\documentclass[titlepage,12pt,openany]{book}
\usepackage{amssymb,mssymb}
\usepackage{graphicx}

\newcommand{\mycomment}[1]{}

\newcommand{\nz}{\mathbb{N}}
\newcommand{\rz}{\mathbb{N}}
\newcommand{\cz}{\mathbb{N}}

\newcommand\oben[1]{\begin{center}\begin{minipage}{#1}\hrule\medskip}
\newcommand\unten  {\hrule \end{minipage}\end{center}}

\def\rz{\ifmmode{I\hskip -3pt R}
    \else{\hbox{$I\hskip -3pt R$}}\fi} 
\def\nz{\ifmmode{I\hskip -3pt N}
    \else{\hbox{$I\hskip -3pt N$}}\fi} 
\def\gz{\ifmmode{Z\hskip -4.8pt Z}
    \else{\hbox{$Z\hskip -4.8pt Z$}}\fi} 
\def\cz{\ifmmode{C\hskip -5pt \line(0,1){8} \hskip 5pt}
    \else{\hbox{$C\hskip -5pt \line(0,1){8} \hskip 5pt$}}\fi}
\def\qz{\ifmmode{Q\hskip-5.0pt\vrule height6.0pt depth 0pt\hskip6pt}
    \else{\hbox{$Q\hskip-5.0pt\vrule height6.0pt depth 0pt\hskip6pt$}}\fi}

\def\Rz{\ifmmode{I\hskip -3.3pt R}
    \else{\hbox{$I\hskip -3.3pt R$}}\fi} 
\def\Nz{\ifmmode{I\hskip -3.3pt N}
    \else{\hbox{$I\hskip -3.3pt N$}}\fi} 
\def\Gz{\ifmmode{Z\hskip -5.3pt Z}
    \else{\hbox{$Z\hskip -5.3pt Z$}}\fi} 
\def\Cz{\ifmmode{C\hskip-5.5pt\vrule height6.5pt\hskip6.9pt}
    \else{\hbox{$C\hskip-5.5pt\vrule height6.5pt\hskip6.9pt$}}\fi}
\def\Qz{\ifmmode{Q\hskip-5.6pt\vrule height6.6pt depth 0pt\hskip6.6pt}
    \else{\hbox{$Q\hskip-5.6pt\vrule height6.6pt depth 0pt\hskip6.6pt$}}\fi}

\newtheorem{Theorem}{Theorem}
\newtheorem{Definition}{Definition}

\newtheorem{Corollary}{Corollary}
\newtheorem{Proposition}{Proposition}
\newtheorem{Remark}{Remark}
\newtheorem{Lemma}{Lemma}
\newtheorem{Example}{Example}

\newcommand{\nin}{{n \in {\Bbb N}_0}}

\newcommand{\supp}{\mbox{supp }}

\newenvironment{proof}{{\bf Proof}. $\;\;$}{\hspace*{\fill} $\diamond$}
\newenvironment{Proof}{{\bf Proof}. $\;\;$}{\hspace*{\fill}
$\diamond$}

\begin{document}
\thispagestyle{empty}

\begin{center} 

\vspace{2cm}

 {\Large\bf Positive Definite Kernels and Random
Sequences Connected to Polynomial Hypergroups}

\vspace{1cm}

{\Large \bf Volker H\"osel}

\vspace{3cm}

\includegraphics[width=8cm]{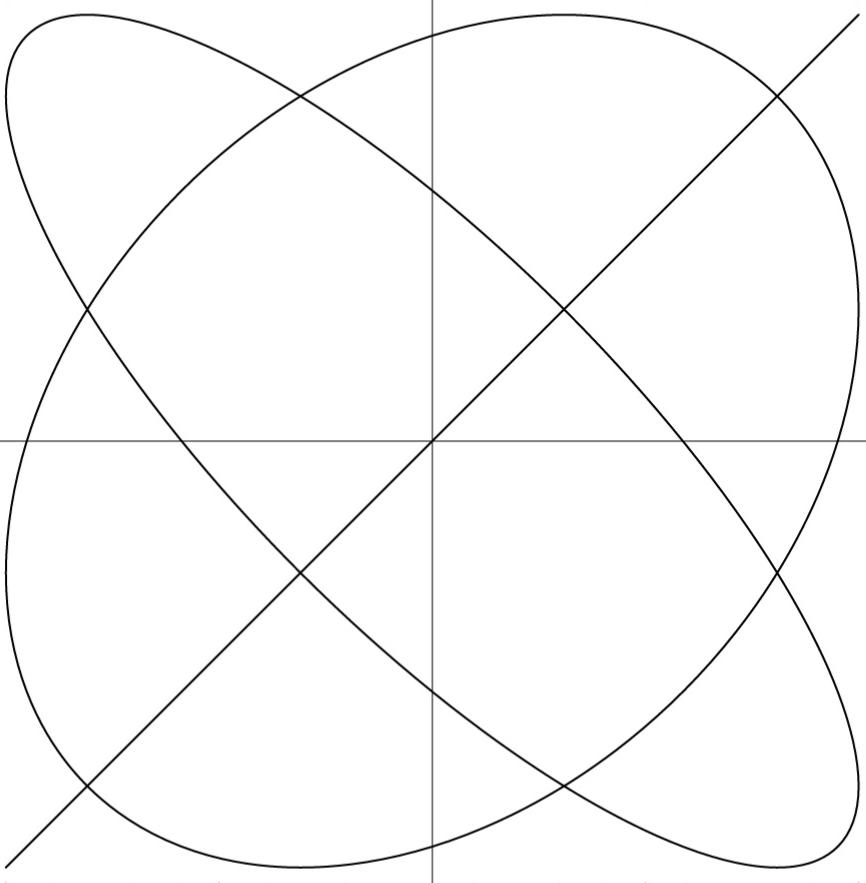}

\vspace{3cm}

{\large \bf November 2024}

\end{center}
\newpage

\setcounter{page}{1}

\vspace{6cm}

\tableofcontents

\addcontentsline{toc}{chapter}{Introduction}

\chapter*{Introduction}

\markboth{\MakeUppercase{Introduction}}{}

The development of stochastic time series modelling advanced significantly with Khintchine's 1934 introduction of weakly stationary stochastic processes \cite{Khin}. Since then, the field has evolved substantially through contributions from researchers like Bochner \cite{Bo}, Cramér \cite{Cra}, and Kolmogorov \cite{Kol}, finding important applications across various domains. A major advantage of weak stationarity is the range of effective strategies it offers for parameter estimation and prediction, as demonstrated by Wiener and Masani \cite{WiMa1}, \cite{WiMa2}.

However, even weak stationarity is often too restrictive for real-world applications. To relax this assumption researchers developed tools to address specific types of nonstationary processes. Treatments of these typically rely on concrete assumptions about the covariance structure, as seen in the works of Getoor \cite{Geto}, Gladyshev \cite{Glad}, Kampé de Fériet \cite{Kam}, and Miamee \cite{Mia}. Further studies, such as those by Rao \cite{Cha}, \cite{RaoHa}, \cite{RaoView} and Hurd \cite{HuRep}, \cite{HuCor}, have approached these processes from an abstract mathematical perspective. On the other hand, there remain limited results on concrete algorithms estimating the covariance structure and making predictions from realizations of such processes.

Traditional weakly stationary processes, indexed by locally compact abelian groups, could serve as a guideline to construct such estimation algorithms. Here, the covariance function is positive definite with respect to group multiplication and Bochner and Cram\'{e}r theorems yield representations of the covariance function and the random variables themselves. Many of the classical results are based upon such representations and thorough investigations of this field may be found in the books of Kakihara \cite{Kak}, Mandrekar and Salehi \cite{Man}, Priestley, \cite{Pri} or Yaglom \cite{Yag}.

In this work, we develop strategies for covariance estimation and prediction for processes indexed by hypergroups, leveraging the positive definiteness of covariance functions. Our analysis extends beyond hypergroup-stationary processes, exploring generalizations such as hypergroup-harmonizable, hypergroup-cyclostationary, and hypergroup-asymptotically stationary processes. For comparison, similar concepts in the classical case are found in the work of Dehay \cite{Deh}, Hanin and Schreiber \cite{Han}, Houdré \cite{Hou}, Miamee and Salehi \cite{MiaSal}, Niemi \cite{Nie}, and Schreiber \cite{SchAs}.

To construct specific algorithms, we focus on discrete stochastic sequences related to orthogonal polynomials. Some results apply to polynomials without necessarily inducing hypergroup structures. But, in general we need Bochner's theorem and thus require that a polynomial hypergroup is involved. In hypergroup settings, estimation strategies are more intricate than in traditional cases. The covariance structure often involves mixtures of random variables that are challenging to disentangle.

Throughout this work, we assume the random sequences comprise of square integrable variables $X_n$, which are elements of a Hilbert space $L^2(\Omega,P)$. For our investigations, only the covariance structure is relevant. Therefore, we assume further that the variables are centred ($E(X_n)=0$). According to Aronszajn \cite{Aron} every positive definite kernel induces an inner product in an appropriate Hilbert space (and vice versa). Consequently, covariances of square integrable processes and positive definite kernels can be identified and studied without regressing to the associated stochastic process.

Following this concept, the first chapter focuses on various types of positive definite kernels. We examine kernels related to specific stochastic sequences, starting with classical structures like cyclostationary and harmonizable kernels, then introducing analogous structures in hypergroup settings. For example, we characterise $P_n$-harmonizable kernels that are T-cyclostationary or prove Wiener type theorems for spectral measures of harmonizable kernels. 

In Chapter 2, we discuss random sequences associated with these kernels, offering examples of processes where they occur.

Chapter 3 introduces our results on estimating the covariance structures of $P_n$-weakly stationary sequences. We derive several estimators and establish mean square consistency under reasonable conditions. Additionally, we introduce periodograms to detect periodicities in $P_n$--harmonizable sequences and propose Fej\'{e}r-based methods for estimating spectral densities.

In Chapter 4, we present a comprehensive theory for predicting $P_n$-weakly stationary sequences, including best linear predictors and prediction errors for various polynomial hypergroups. We also introduce a new concept of prediction with additional information for $MA(q)$ sequences.

We conclude with a Levinson-type algorithm for the coefficients of the best linear predictor. This efficient algorithm, with a complexity of $O(n^2)$, is designed for positive matrices defined by polynomials, providing a faster alternative to the classical Cholesky decomposition, which has a complexity of $O(n^3)$.

\chapter{Positive definite kernels}
\section{Classical structures}

In this section, we present examples of positive definite kernels
which are associated to important classes of stochastic processes.
These classes are the weakly stationary, cyclostationary,
harmonizable and asymptotically stationary sequences. First we
state basic facts on positive definite kernels and positive
definite functions.

\subsection{Kernels on groups and semigroups}

Positive definite kernels can be found in a wide range of
scientific areas. Their appearance can be traced back at least to
1909, when Mercer \cite{Mercer} presented a paper on integral
equations. They are studied for example in a more abstract
framework in connection with linear operators on Banach spaces or
are used in applied fields like machine learning. A main reference
is the book of Berg, Christensen and Ressel \cite{Ressel}. We
consider in this text positive definite kernels as covariances of
stochastic sequences. First we give the basic definition and state
some general properties.
\newpage
\begin{Definition}
Let $A$ be any set. A complex valued function $K(x,y)$ on $A
\times A$ is called hermitian kernel if
\mbox{$K(x,y)=\overline{K(y,x)}$} for all $x,y \in A$. A hermitian
kernel is called positive definite kernel if for any positive
integer $n$ and every choice of elements $x_0, \dots, x_n$ in $A$
and complex numbers $\lambda_0 \ldots \lambda_n$ the inequality
\begin{eqnarray}
\sum_{k,l=0}^n \lambda_k \overline{\lambda_l} K(x_k,x_l) \geq 0
\end{eqnarray}
holds. If the inequality is strict one calls the kernel strictly
positive definite.
\end{Definition}
The theorem of Aronszajn \cite{Aron} implies that every positive
definite kernel can be seen as covariance of a set of square
integrable random variables indexed by $A$. We follow the proof in
Parthasarathy \cite{Part}.
\begin{Theorem}(Aronszajn)
Let $A$ be any set and $K(x,y)$ a positive definite kernel on $A
\times A$. Then there exists a complex Hilbert space $H$ and a
mapping $\phi$ from $A$ into $H$ such that the set $\{\phi(x),x
\in A\}$ span $H$ and
\begin{eqnarray}
K(x,y)= \langle \phi(x), \phi(y) \rangle
\end{eqnarray}
where $\langle , \rangle$ is the scalar product of $H$.
\end{Theorem}
\begin{proof}
For any finite set of elements $x_0, \dots, x_n$ in $A$ the matrix
with entries $K(x_k,x_l)$ for $0 \leq k,l \leq n$ is positive
definite. Therefore there exist n dimensional (complex) Gaussian
probability distributions $\mu(x_1,\ldots,x_n) $ with mean zero
having this matrix as covariance function (compare for example
\cite{Shir}, Chap.II $\S 8$). The probability distributions
$\mu(x_1,\ldots,x_n)$ are consistent by construction. Let $X$ be
the Borel space of all complex valued functions on $A$ with the
smallest $\sigma$-field such that all projections $\pi_x: f
\rightarrow f(x)$ from $X$ to $\mathbb{C}$ are measurable.
Kolmogorov's extension theorem guarantees the existence of a
measure $\mu$ on $X$ such that the joint distribution of
$(f(x_1),\ldots,f(x_n))$ is $\mu(x_1,\ldots,x_n)$. Now defining
$\phi(x):=f(x)$ for $x \in A$ and $f \in X$, then $\phi(x)$ is a
function on $X$ which is an element of the Hilbert space
$L_2(\mu)$. Further, one has $<\phi(x),\phi(y)>=\int_X
f(x)\overline{f(y)}d\mu(f)=K(x,y).$ Now, defining as $H$ the
closure of the linear span of all $\phi(x)$ with $x \in A$ in
$L_2(\mu)$ finishes the proof.

\end{proof}

\begin{Remark}
There exist different approaches to Aronszajn's theorem. Another
utilizes reproducing kernel Hilbert spaces not resorting to
stochastics (for example \cite{Ressel}). The proof given above
demonstrates also the equivalence of positive definite kernels and
covariance matrices of centered Gaussian processes. This
equivalence will we have in mind, when we consider kernels in the
sequel. The investigation of such kernels is worthwhile in itself.
But, in a later chapter we also look at stochastic sequences
having such kernels as covariance. There, it is of especial
interest if such kernels have a spectral representation or induce
positive definite functions (compare \cite{Sas}).
\end{Remark}

\begin{Definition}
A complex valued function $\phi$ on a semigroup ($S,\diamond$)
with involution $*$ is called positive definite if $s \times t
\mapsto \phi(s^*\diamond t)$ is a positive definite kernel on $S
\times S$, i.e. if one has for any positive integer $n$ and every
choice of elements $s_0, \dots, s_n$ in $S$ and complex numbers
$\lambda_0 \ldots \lambda_n$
\begin{equation}\label{posFunction}
    \sum_{k,l =0}^n \lambda_k \overline{\lambda_l} \phi(s_k^* \diamond s_l) \geq
    0.
\end{equation}
\end{Definition}
For bounded positive definite functions $\phi$ on abelian
semigroups the theorem of Lindahl and Maserick \cite{Lindahl}
proves the spectral representation
\begin{equation}\label{LindahlSpec}
    \phi(s)=\int_{\hat{S}} \rho(s)d\mu(\rho) \mbox{  for } s \in S.
\end{equation}
Here, $\mu$ is a uniquely determined positive measure on the
restricted dual semigroup $\hat{S}=\{\rho \in S^*\;|\;|\rho(s)|
\leq 1 \mbox{ for all } s \in S \}$. $S^*$ is the set of
semicharacters on $S$. These are complex functions on $S$
satisfying with the neutral element $0$ and the composition $+$:
$\rho(0)=1$, $\rho(s+t)=\rho(s)\rho(t) \mbox{ for }$ and
$\rho(s^*)=\overline{\rho(s)} \mbox{ for } s,t \in S$. For
discrete abelian groups ($s^*=-s$) the theorem reduces to the well
known Bochner theorem. The analogue Bochner theorem for discrete
abelian hypergroups will be use later in this text. A historical
survey of positive definite functions and generalizations may be
found in \cite{Stewart}.

\subsection{Weakly stationary and harmonizable kernels}

There is a well developed theory about stationary kernels on a
variety of groups $G$ of different type. For abelian groups with
group operation written additively the decisive definition reads
as

$$K(s+t,s)=K(t,0) \mbox{ for all } s,t \in G.$$

The early investigations started with the works of Khintchine
\cite{Khin} on the correlation theory of stationary processes.
Spectral theory proved crucial for all the developments. The
function $R(t):= K(t,0)$ is positive definite and thus Bochners
\cite{Bochner} theorem guarantees the existence of a finite
positive measure $\mu$ on the dual group $\hat{G}$ such that $R$
is the Fourier-Stieltjes transform of $\mu$:

$$R(t)=\int_{\hat{G}}\gamma(t)d \mu(\gamma).$$

The spectral representation of the associated process was given in
the works of Cram\'{e}r \cite{Cra} and Kolmogorov \cite{Kol}. A
theorem of Stone \cite{Stone} on the spectral representation of
groups of unitary operators paved the way. If, for example, $X_n,
\; n \in \Bbb Z$ is a e sequence of square integrable random
variables, then there is a representation
$$X_n=\int_T e^{int}dZ(t)$$
with $Z$ an orthogonal stochastic measure on the Torus $T$.

For multidimensional stationary processes spectral theory is
connected with names like Cram\'{e}r \cite{Cra}, Kolmogorov
\cite{Kol}, Rozanov \cite{Roza}, Wiener and Masani
\cite{WiMa1},\cite{WiMa2}.

Generalization exist in several directions. So the index set has
been relaxed and for example certain non abelian groups of type I
have been considered. On the other hand the domain of the process
can be generalized from Hilbert spaces to Banach spaces or even
more general structures (for a collection of concepts see
\cite{Kak}).

From a theoretical point of view the subject stationarity has
matured satisfyingly and relevant topics for applications like
estimation of the spectral measure or prediction of the process
are considered as completed now.

Nevertheless, stationarity is rarely met in real world phenomena.
Structures more suited to wider applications have been proposed.
There is a trade off between more general definitions which do not
allow for a satisfactory theory and too specific ones like in the
stationary case.

We present some of the most important generalizations which will
be the guideline for our study in connection with $P_n$-weakly
stationary sequences.

A concept using spectral representation of the kernel is
harmonizability. This means that one has in the discrete case a
representation of the positive definite kernel $K$:
$$K(s,t)=\int_{[- \pi,\pi]\times [- \pi,\pi]} \! \! exp(i s x- it y) \;
d \mu(x,y).$$ Here, $\mu$ is a (finite) complex measure on $[-
\pi,\pi]\times [- \pi,\pi]$ which satisfies for every $n \in \Bbb
N$, $a_1,...,a_n$ and Borel sets $A,B, A_1...A_n$ the conditions
\begin{eqnarray*}
a) &\mu(A,B)=\overline{\mu(B,A)}\\
 b) &\quad \sum_{k,l=1}^n a_k \;
\overline{a_l} \; \mu(A_k,A_l)\; \geq 0.
\end{eqnarray*}

This definition was given by Lo\`{e}ve end of the 1940's
\cite{loeve} and is nowadays referred to as strongly harmonizable.
More general concepts have been introduced by Bochner \cite{Bo}
called V-boundedness and by Rozanov which he also called
harmonizable. The latter definitions are equivalent and are now
denoted as weakly harmonizable (\cite{ChaRao},\cite{Niemi}). Here,
a) and b) are still required, but $\mu$ needs only to be a
bimeasure, meaning, that for every fixed Borel set $A$ the set
function $\mu(A,.)$ is a measure with respect to the second set
variable and conversely. Further readings are \cite{Hou},
\cite{Houdre}.

\subsection{Cyclostationarity}
Another important generalization of weakly stationary processes
are periodically correlated (or cyclostationary) processes and
almost periodically correlated processes. The definitions were
introduced by Gladyshev \cite{Glad}. His work on multivariate
stationary processes \cite{Glady} led to cyclostationarity: There
is a bijection between second order periodic correlated sequences
with period T and T-dimensional stationary vector sequences. Many
researchers worked on these sequences and their continuous time
counterparts. We especially mention the contributions of Hurd
\cite{HuRep}, \cite{HuCor}, \cite{Hurd}. A recent survey with a
large collection of literature is \cite{Gard}.

For our purpose, we consider the associated spectral theory for
the discrete case. It can be compared to our results for
hypergroup cyclostationarity, derived later.

The definition of cyclostationary kernels according to Gladyshev
is: For all $s,t \in \Bbb Z$ and one positive integer T it is true
that
$$K(s,t)=K(s+T,t+T).$$
As $K(t,t+s)$ is for every $s \in \{0,...,T-1\}$ a periodic
function in t, discrete Fourier techniques can be applied. The
coefficients are
$$B_k(s)= \frac{1}{T} \sum_{t=0}^{T-1} \; K(t+s,t) exp(-2 i \pi k t/T)$$
and the kernel is representable as the finite Fourier sum
$$K(t,t+s)= \sum_{k=0}^{T-1} B_k(s) \; exp(2 i \pi k t/T).$$

A set of discrete Fourier transforms obtained from a periodic
function $K(t,t+s)$ as above needs to have special properties to
guarantee that $K$ is a positive definite kernel. Gladyshev's
characterization is:

$K$ is a positive definite cyclostionary kernel if and only if for
every choice of $N$, complex numbers $\lambda_n$, and integers
$k_n \in \{0,...,T-1\}$
$$\sum_{n,m=0}^N \lambda_n \overline{\lambda_m} \; \; \beta_{k_n k_m}(n-m) \geq 0$$
holds. Here, one defines $\beta_{k l}(s):=B_{l-k}(s) \; exp(2 i
\pi k s/T)$.

Further, Gladyshev showed for the discrete case that periodically
correlated kernels are (strongly) harmonizable and that the
spectral measure is concentrated on the $2T-1$ diagonal lines
$$\{(x,y)|\, y=x-2 \pi k/T, \; k=-(T-1),...,T+1 \}.$$ This is
illustrated for $T=5$:

\begin{center}
\includegraphics[width=8cm]{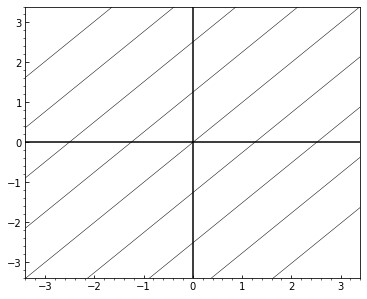}
\end{center}

\newpage

\subsection{Asymptotic stationarity}
Asymptotic stationarity is a concept which was developed in the
works of Kamp\'{e} de F\'{e}riet and Frenkiel \cite{Kam},
\cite{KampFer}, Parzen \cite{Parzen} and Rozanov \cite{Rozanov}.
Regarding positive definite kernels $K$ the definition in the
discrete case  is as follows:

$K$ is called asymptotically stationary if the limit of averages
$$\tilde{K}(s)= \lim_{n \rightarrow \infty} \frac{1}{n}\sum_{k=1}^n K(k,k+s)$$
exists for all $s \in \Bbb Z$ and is positive definite. In this
case \~{K} defined by
$$\tilde{K}(n,m):= \tilde{K}(n-m)$$
is a positive definite kernel which is weakly stationary. It turns
out, that strong harmonizable and cyclostationary kernels are
asymptotic stationary. An overview and generalizations are
presented in Schreiber \cite{SchAsRel}, \cite{SchAs}.

\section{Kernels related to polynomial hypergroups}

The concepts leading to hypergroups have a long history. They date
back to the rise of group theory around 1900. In the mid thirties
hypergroups were studied in connection with cosets of nonabelian
groups (Marty \cite{Marty}, Wall \cite{Wall}). The formal
definition involves algebraic and topological conditions.

In our context, we consider exclusively polynomial hypergroups.
These are discrete and abelian structures which do not require to
introduce hypergroups in full generality. In the following, we
present the basic facts necessary for our investigations.

Let $(P_n(x))_{\nin}$ be a polynomial sequence with degree$(P_n)
=n $, that is orthogonal with respect to a probability measure
$\pi \in M^1({\Bbb R}).$ We always assume that $\pi$ has infinite
support: $|\supp \pi| =  \infty$. Assume that $P_n(1) =1 \;$ and
that all the linearization coefficients $g(m,n,k),$ defined by
\begin{equation}\label{Linearization}
P_m(x) \; P_n(x) \; = \; \sum_{k=|n-m|}^{n+m} g(m,n,k) \; P_k(x)
\end{equation}
are nonnegative. Then we denote the convex combination of point
measures $\epsilon_k$ on ${\Bbb N}_0$ by
$$ \varepsilon_m \ast \varepsilon_n := \;
\sum_{k=|n-m|}^{n+m} g(m,n,k)\; \varepsilon_k $$ and call
$\varepsilon_m \ast \varepsilon_n$ the convolution of $m$ and
$n.\;$ With $\ast$ as convolution, the identity map as involution
and 0 as unit element, ${\Bbb N}_0$ becomes a commutative
hypergroup, that is called a polynomial hypergroup, see Bloom and
Heyer \cite{Bloo}. Section 3.3 of their book contains a long list
of polynomial hypergroups. Every character on the polynomial
hypergroup ${\Bbb N}_0$ is given by ${\alpha_x:{\Bbb N}_0
\rightarrow {\Bbb R},\; \; \alpha_x(n) = P_n(x),\;}$ where $x \in
D_s$ with
\begin{equation}
D_s \; = \; \left\{ x \in {\Bbb R}: \; \{ P_n(x) : \nin \} \;
\mbox{ is bounded} \right\}.
\end{equation}
The character space $\widehat{{\Bbb N}_0}$ is homeomorphic
to $D_s,$ (see \cite{Bloo}). \\
In particular $D_s$ is a compact subset of ${\Bbb R}$ and
$|P_n(x)| \leq 1$ for every $\nin$ and every $x \in D_s$.

The orthogonal polynomials $P_n(x)$ are determined by their
three-term recurrence relation
$$ P_0(x) \; =1,\qquad P_1(x) \; = \; \frac{1}{a_0} \; (x-b_0) \qquad \mbox{ and for } \nin $$
\begin{equation} \label{threeTerm}
P_1(x) \; P_n(x) \; = \; a_n P_{n+1} (x) \; + \; b_n P_n(x) \; +
\; c_n P_{n-1} (x)
\end{equation}
where $(a_n)_\nin,  (b_n)_\nin,  (c_n)_\nin \;$ are three
real-valued sequences with ${a_n,c_n >0, \; b_n \geq 0 \;}$ and
$a_n +b_n +c_n \; =1.\;$ For $a_0,b_0 \in {\Bbb R}$ we assume
$a_0>0$ and $a_0+b_0 =1.\;$

The linearization coefficients $g(m,n,k)$ can be directly
calculated from the recurrence coefficients $a_n,b_n,c_n,\;$ see
Lasser \cite{la01}.

The Haar measure $h$ on the (discrete) polynomial hypergroup is
given by $h(n)\; = \; g(n,n,0)^{-1}$ and also determined by the
weights
\begin{equation}\label{HaarProd}
h(0) \; =1, \quad h(1) \; = \; 1/c_1, \quad h(n) \; = \;
\prod_{k=1}^{n-1} a_k / \prod_{k=1}^n c_k \qquad n=2,3,...
\end{equation}

The support of the orthogonalization measure $\pi$ is contained in
$D_s$. In fact we have
\begin{equation}
\supp \pi \; \subseteq D_s \; \subseteq \; [1-2a_0,1]
\end{equation}
and $\pi$ is the Plancherel measure on $D_s$, i.e.
$$ \int_{D_s} P_n(x) \; P_m(x) \; d \pi (x) \; = \; \left\{ \begin{array}{ll}
   0 & \mbox{ if } n \neq m \\
\displaystyle\frac{1}{h(n)} & \mbox{ if } n=m.
\end{array} \right. $$
We essentially use a Bochner theorem for polynomial hypergroups
(see Theorem 4.1.6 of \cite{Bloo}), which states that for each
bounded positive definite sequence $(d(n))_{\nin}$ there exists a
unique $\mu \in M^+ (D_s)$, called spectral measure, such that
\begin{equation}
d(n) \; = \; \int_{D_s} P_n(x) \; d\mu (x). \label{spec1}
\end{equation}
Hereby positive definiteness stands for
\begin{equation}\label{Hyperposdef}
\sum_{i,j=1}^n \lambda_i \; \overline{\lambda_j} \; \;
\varepsilon_{m_i}\ast\varepsilon_{m_j} (d) \; \geq 0
\end{equation}
for all $\lambda_1,...,\lambda_n \in {\Bbb C}$ and $m_1,...,m_n
\in {\Bbb N}_0.$

To clarify we note that
$$ \varepsilon_{l}\ast\varepsilon_{m} (d) \; = \; \sum_{k=|m-l|}^{m+l} g(l,m,k) \;
d(k).$$

\begin{Example}
The paradigm of our investigations is the hypergroup induced by
the Chebyshev polynomials  of the first kind $T_n(x)=\cos(n t)$,
$x=\cos(t)$. These polynomials are orthogonal on $[-1,1]$ with
respect to (the Plancherel measure)
$$d\pi(x)=\frac{d x}{\pi \sqrt{1-x^2}}.$$
Here, one has $\supp \pi = D_s$ and the Haar measure is $h(0)=1$
and $h(n)=2$ for $n>0$. Due to the product formula for cosine one
has the linearization
\begin{equation}\label{ChebLin}
    T_n(x) T_m(x)= \frac{1}{2}T_{|n-m|}(x)+ \frac{1}{2}T_{n+m}(x).
\end{equation}
\end{Example}

The growth conditions of the Haar measure play an important
r\^{o}le in our investigations. The preceding example shows one of
only a few known cases with bounded Haar measure. We need the
weaker growth conditions (H) which seems to be essential in the
theory of polynomial hypergroups:
\begin{equation}\label{H}
    (H) \qquad \lim_{n \rightarrow \infty}
    \frac{h(n)}{\sum_{k=0}^n h(k)} = 0.
\end{equation}
Condition (H) will proof helpful in the later applications of the
Christoffel-Darboux identity (\cite{laro},(12)), for $x \neq y$
\begin{equation}\label{Christoffel1}
    \sum_{k=0}^n P_k(x) P_k(y) h(k) = a_0 a_n h(n)
    \frac{P_{n+1}(x)P_{n}(y)-P_{n}(x)P_{n+1}(y)}{x-y}.
\end{equation}

\subsection{$P_n$-weakly stationarity}

The definition of $P_n$-stationary kernels is guided by the
classical weakly stationarity $K(n,m)=K(|n-m|,0)$. As there is no
inverse operation in hypergroups $H$ one has to use the
convolution alone. We give the definition for polynomial
hypergroups only. For definitions and investigations related to
more general hypergroups see for example \cite{lalest}.

\begin{Definition}
Let $(P_n)_{\nin}$ induce a polynomial hypergroup with convolution
$*$. A positive definite kernel $K$ is called $P_n$-weakly
stationary if it is bounded and satisfies for all $n,m \in {\Bbb
N}_0$
\begin{equation}\label{WeakStat}
    K(n,m)=K(n*m,0).
\end{equation}
\end{Definition}
Then, $P_n$-weakly stationarity explicitly means
\begin{equation}\label{KernelWeak}
    K(m,n) \; = \; \sum_{k=|m-n|}^{m+n} g(m,n,k)\; K(k,0).
\end{equation}

The definition implies that $(K(n,0))_{\nin}$ is a bounded
positive definite sequence on the hypergroup.  Thus, Bochner's
theorem yields a spectral representation of the kernel $K$ as well
(compare \ref{spec1}):
\begin{equation}\label{SpecK}
    K(n,m)=\int_{D_s} P_n(x) P_m(x) d\mu(x).
\end{equation}
For later applications, we note the equalities (valid for all
$s,n,m \in \Bbb N_0$)
\begin{equation}\label{SpecKConv}
     K(n*s,m)=\int_{D_s} P_{n*s}(x) P_m(x) d\mu(x)
     = \int_{D_s} P_{s}(x) P_n(x) P_m(x) d\mu(x).
\end{equation}


\subsection{$P_n$-cyclostationarity}
The structure on the kernels corresponding to cyclostationary
processes
$$K(x,y)=K(x+T,y+T) \quad \hbox{for all } x,y \in \Bbb Z \quad\hbox{and
some integer  } T>0$$ seems to have the canonical equivalence
\mbox{$K(n,m)=K(n*T,m*T)\quad\hbox{for all } n,m \in {\Bbb N}_0.$}
But this does not extend hypergroup stationarity since
$K(n*m,0)=K(n*T,m*T)$ does not hold in general for any fixed
$T>0$:

$K(n*m,0)$ is a function of $K(|n-m|,0),\ldots,K(n+m,0)$, but
$K(n*T,m*T)$ depends on $K(n+m+2T,0)$ since the corresponding
coefficient
$$g(n,T,n+T)g(m,T,m+T)g(n+T,m+T,n+m,2T)$$
is not zero. This holds since $g(n,m,n+m) \neq 0$ for all $n,m \in
\Bbb N_0$ (compare \ref{Linearization}). An elementary observation
guides to a more appropriate definition of cyclostationarity for
hypergroups:

\begin{Proposition}
For a positive definite kernel $K(x,y)$ defined on $\mathbb{Z}$
the following conditions are equivalent
\begin{eqnarray*}
C1:  &&  K(x,y)= K(x+T,y+T)\\
C2:  &&  K(x,y)= K(x-y+y \mbox{ mod T},y \mbox{ mod T})\\
C3: &&  K(x+T,y)= K(x,y-T)
\end{eqnarray*}

where $\mbox{ mod T}$ means the remainder of dividing y by T and
all equations have to be valid for all $x,y \in \Bbb Z$ and some
integer $T>0$.
\end{Proposition}

\begin{Proof}\\
$C1$ equivalent to $C3$ is trivial: substitute $y$ by $y+T$ in
$C3$.

$C1$ implies $C2$: This is evident from the chain of equalities
\begin{eqnarray*}
K(x,y)= && K(x-T,y -T) = \ldots = K(x-N_y T, y \mbox{ mod T})\\
= && K(x-y+y \mbox{ mod T},y \mbox{ mod T})
\end{eqnarray*}
with $N_y$ the quotient of y divided by $T$.

$C2$ implies $C1$: The substitution $\tilde{x}=x+T$,
$\tilde{y}=y+T$ in the equation C2 gives with $y \mbox{ mod T}
=\tilde{y} \mbox{ mod T}$
$$K(\tilde{x}-T,\tilde{y}-T)=K(\tilde{x}-\tilde{y}+\tilde{y} \mbox{ mod T},
\tilde{y} \mbox{ mod T})$$ Renaming $\tilde{x}, \tilde{y}$ to
$x,y$ again and comparing it to $C2$ yields $K(x,y)=K(x-T,y-T)$
and thus $C1$.

$C1$ equivalent to $C3$ is clear from the shifts  $\tilde{x}=x+T$
and $\tilde{x}=x-T$ with $\tilde{y}=y$.

\end{Proof}

\begin{Remark}
For $T=1$ we have stationarity. The analogous conditions $C1$ and
$C3$ for abelian groups are obviously also equivalent. $C2$ and
its equivalence to the other conditions could be shown for other
free abelian groups than $\mathbb{Z}$. In the case of more general
structures like monoids or hypergroups the lack of inverse
elements require different conditions which are in general not
equivalent.
\end{Remark}

For polynomial hypergroups we consider
\begin{eqnarray*}
HC1: &&  K(n,m)= K(n*T,m*T)\\
HC2: &&  K(n,m)= K(n*m*m \mbox{ mod T},m \mbox{ mod T})\\
HC3: &&  K(n*T,m)= K(n,m*T)
\end{eqnarray*}
Each condition has to be valid for all $n,m \in {\Bbb N}_0$ and
some integer $T>0$.

$HC1$ induces interesting structures on hypergroups and could be
used for semigroups as well. As we have seen before, it does not
extend the definition of stationarity. So we do not discuss it
here.

$HC2$ indeed includes stationarity ($T=1$) and is applicable to
polynomial hypergroups. It has no equivalence in semigroups or
more general hypergroups and should be investigated elsewhere.

$HC3$ is applicable to hypergroups and semigroups. We show that
for $T=1$ the condition is equivalent to stationarity.

\begin{Remark}
Another condition including stationarity $(T=0)$ is
 $$K(n,m)=K(n*m*T,T).$$
Here, as in $HC2$ the second component $m$ of $K$ is shifted to
the first component. This structure has been investigated in
\cite{Hof}. The condition implies that for $T>0$ the entire kernel
$K(n,m), \; n,m \in \Bbb N_0$, is given  by the $2T$ values
$K(T,0),\ldots,K(2T-1,0)$ only. Therefore we will not use this
restrictive condition.
\end{Remark}

For the proof that $K(n*1,m)=K(n,m*1)$ is equivalent to
stationarity we use the following lemma.
\begin{Lemma}
$K(n*T,m)=K(n,m*T)$ for all $n,m \in {\Bbb N}_0$ and a fixed $T
\in \Bbb N$ implies
\begin{equation}\label{convo}
    K(n*T^{*s},m)=K(n,m*T^{*s})
\end{equation}
for all $s \nin$. Here we define $T^{*s}$ as the $s$ time
convolution $T*\ldots*T$ of $T$ with itself  and set $T^{*0}:=0$
and $T^{*1}:=T$.
\end{Lemma}

\begin{Proof}
We use induction on $s$. The case $s=0$ is trivial and $s=1$ gives
the assumption. Now assume $K(n*T^{*t},m)=K(n,m*T^{*t})$ holds for
$t<s$ and all $n,m \in {\Bbb N}_0$. With the associativity of $*$
one has
\begin{eqnarray*}
K(n*T^{*s},m)&=&K(n*T*T^{*(s-1)},m)=\sum_{k=|n-T|}^{n+T}g(n,T,k)\,K(k*T^{*(s-1)},m)\\
&=& \sum_{k=|n-T|}^{n+T}g(n,T,k)\,K(k,m*T^{*(s-1)})=
K(n*T,m*T^{*(s-1)}) \\
&=&\sum_{k=0}^{m+(s-1)\cdot T} a(k)\,
K(n*T,k)=\sum_{k=0}^{m+(s-1)\cdot T} a(k)\, K(n,k*T) =
K(n,m*T^{*s}).
\end{eqnarray*}
The positive numbers $a(k)$ result from the respective
convolutions.
\end{Proof}

\begin{Theorem} \label{nStar1}
Let $K(n,m)$ be a positive definite kernel defined on a polynomial
hypergroup. Equivalent are:
\begin{eqnarray*}
a) && K(n,m)=K(n*m,0) \quad \mbox{(stationarity)}\\
b) && K(n*1,m)=K(n,m*1) \\
c) && K(n*t,m)=K(n,m*t) \quad \mbox{for all } t \in {\Bbb N}_0
\end{eqnarray*}
Each condition is meant to be valid for all $n,m \in {\Bbb N}_0$.
\end{Theorem}

\begin{proof}
a) implies c) and thus b) :
\begin{eqnarray*}
K(n*t,m)&=&\sum_{k=|n-t|}^{n+t}g(n,t,k)\,K(k,m)=
\sum_{k=|n-t|}^{n+t}g(n,t,k)\,K(k*m,0)=K(n*t*m,0)\\
&=&\sum_{k=|m-t|}^{m+t}g(t,m,k)\,K(n*k,0)
=\sum_{k=|m-t|}^{m+t}g(t,m,k)\,K(n,k) =K(n,m*t)
\end{eqnarray*}
Here, we use associativity and commutativity of $*$.

Now, we show by induction that b) implies c).

Assume that $K(n*t,m)=K(n,m*t)$ is valid for $t<s$ and all
$m,\nin$. Then the equations
\begin{eqnarray*}
K(n*1^{*s},m)&=& \sum_{k=0}^{s} a(k) \, K(n*k,m)\\
K(n,m*1^{*s})&=& \sum_{k=0}^{s} a(k) \, K(n,m*k)
\end{eqnarray*}
yield together with the lemma: $a(s) \, K(n*s,m)= a(s) \,
K(n,m*s)$. Since $a(s)$ is a product of factors of the form
$g(1,k,k+1)$ it is not zero and thus we have  $K(n*s,m)=
K(n,m*s)$.

Finally, setting $m=0$ in c) implies a).
\end{proof}

From the discussion above the following definition is well
motivated.

\begin{Definition}
A positive definite kernel $K$ defined on a polynomial hypergroup
is T-cyclostationary if there is a minimal $T> 0$ with
\begin{equation}\label{T-cyclo}
    K(n*T,m)=K(n,m*T) \quad \mbox{for all } n,m \in {\Bbb N}_0
\end{equation}
where * is the convolution of the hypergroup.
\end{Definition}

The theorem states the equivalence of 1-cyclostationarity and
stationarity for polynomial hypergroups.


\begin{Example}\label{cycloExam}

Consider for some real $C \neq 0$ the kernel on $\Bbb N_0$ given
by
\begin{equation}\label{CycExamp}
K(n,m)=\cases{(C-1)^2, & for n and m odd;\cr (C+1)^2,& for n and m
even \cr (C-1)(C+1),& else.\cr}
\end{equation}

$K$ is positive definite:
\begin{eqnarray*}
\sum_{n,m}\lambda_n \overline{\lambda_m}K(n,m)&=&(C-1)^2\sum_{n,m
\,odd} \lambda_n \overline{\lambda_m}+(C+1)^2 \sum_{n,m \,
even}\lambda_n
\overline{\lambda_m}\\
&&+(C-1)(C+1) \sum_{n+m\, odd}\lambda_n \overline{\lambda_m}\\
&=& \biggl|(C-1)\sum_{n\, odd}\lambda_n + (C+1) \sum_{n\, even}
\lambda_n \biggr|^2 \geq 0.
\end{eqnarray*}

$K$ is 2-cyclostationary with respect to Chebyshev polynomials of
the first kind $T_n$:
\begin{eqnarray*}
K(n*2,m)&=&\frac{1}{2}K(|n-2|,m)+ \frac{1}{2} K(n+2,m)\\
&=& \frac{1}{2} K(n, |m-2|) + \frac{1}{2} K(n,m+2)= K(n,m*2)
\end{eqnarray*}
is for all $n,m \in {\Bbb N}_0$ obviously satisfied.

$K$ is not stationary:
$$K(1*1,0)=\frac{1}{2}K(0,0)+ \frac{1}{2}
K(2,0)=(C+1)^2 \neq (C-1)^2 = K(1,1).$$

\end{Example}

\subsection{$P_n$-asymptotic stationarity}
For polynomial hypergroups, the concept for asymptotic
stationarity, paralleling the classical situation, could be:

\begin{Definition}

A positive definite kernel $K$ defined on a polynomial hypergroup
given by $(P_n)_{n\nin}$ is called M-asymptotically $P_n$-weakly
stationary if
\begin{equation}\label{Masymp}
    M(s):=\lim_{n \to \infty} \frac{1}{\sum_{k=0}^n h(k)} \sum_{k=0}^n
K(k*s,k)\; h(k)
\end{equation}
is bounded and exists for all $s \in {\Bbb N}_0$.

\end{Definition}

This definition is given in \cite{Hof}. It is not yet clear, if
this is in general an extension of $P_n$-weakly stationarity:

If a $P_n$-weakly stationary kernel is M-asymptotically
$P_n$-weakly stationary then (compare \ref{SpecKConv})
\begin{eqnarray*} M(s)&=& \lim_{n \to \infty} \frac{1}{\sum_{k=0}^n h(k)} \sum_{k=0}^n
\int_{D_s}P_{k*s}(x) P_k(x)h(k)d\mu(x)\\
 &=& \lim_{n \to \infty}
 \int_{D_s}P_s(x) \frac{\sum_{k=0}^n P^2_{k}(x)h(k)}
 {\sum_{k=0}^n h(k)}d\mu(x).
\end{eqnarray*}

This implies that for a spectral measure $\mu=\delta_{x}$ (one
atom only) with any $x \in D_s$
$$C(x):=\lim_{n \to \infty} \frac{\sum_{k=0}^n P^2_{k}(x) h(k)}{\sum_{k=0}^n h(k)}$$
has to exist. Though, true for the known examples, it is still an
unsolved problem if this holds in general.

On the other hand, if $C(x)$ exists for $\mu$ almost all $x \in
\supp \mu$, then by the theorem of dominated convergence (for
example (\cite{rudre},1.34))
$$M(s)= \int_{D_s} P_s(x)C(x) d\mu(x).$$
Now, $M(s)$ is positive definite and the associated spectral
measure is $\tilde{\mu}=C(x) d\mu.$

We look at a specific class:

\begin{Example}
Jacobi polynomials $P_n^{(\alpha,\beta)}$:

These classical polynomials are orthogonal on $[-1,1]$ with
respect to $d \pi(x) = f_{\alpha,\beta}(x) dx$ ($\alpha,\beta
>-1$) with
$$f_{\alpha,\beta}(x)= c_{\alpha,\beta}(1-x)^{\alpha}(1+x)^{\beta}.$$
The constant $c_{\alpha,\beta}$ is calculated to make $\pi$ a
probability measure. According to Gasper \cite{Gasper} these
polynomials induce polynomial hypergroups if and only if
$(\alpha,\beta) \in V$:
$$V=\{(\alpha,\beta)\;|\; \alpha \geq \beta >-1
\mbox{ and }a(a+5)(a+3)^2 \geq (a^2-7a-24)\; b^2 \}.$$ Here,
$a=\alpha+\beta+1$ and $b=\alpha-\beta$. A more accessible subset,
almost covering V is $$W=\{(\alpha,\beta)\;|\; \alpha \geq \beta
>-1 \mbox{ and } \alpha+\beta+1 \geq 0 \}.$$ The asymptotic of
Christoffel functions \cite{Mate} gives for $x$ in the open
interval $(-1,1)$

\begin{equation}\label{Christoffel}
    \lim_{n \rightarrow \infty}\frac{\sum_{k=0}^n
(P_k^{(\alpha,\beta)}(x))^2 h_{\alpha,\beta}(k)}{n} = \frac{1}{\pi
f_{\alpha,\beta}(x) \sqrt{1-x^2}}.
\end{equation}

Since, in the big O notation $h_{\alpha,\beta}= O(n^{2 \alpha+1})$
(see \cite{la01}) the equation \ref{Christoffel} yields for
$\alpha> -1/2$
$$C(x)=0 \mbox{  for  } x \in (-1,1).$$
Now, $C(1)=1$ and also $C(-1)=1$ for the symmetric case
$\alpha=\beta > -1/2$. Application of (\cite{Sz}, 4.1.4) gives
$C(-1)=0$ in case $\alpha > \beta, \alpha > -1/2 $.

Chebyshev polynomials of the first kind $T_n$ are the remaining
case $\alpha=\beta = -1/2$. Here, the asymptotic
(\ref{Christoffel}) directly shows
$$C(x)=1 \mbox{ for  } x \in [-1,1].$$

Altogether, we find (with an application of dominated
convergence):

\begin{eqnarray*}
M(s)= \left\{ \begin{array}{lll}
 \mu(\{1\})  & \mbox{ if  } &  \alpha > \beta, \alpha > -\frac{1}{2} \\
 \mu(\{1\})+ (-1)^s \mu(\{-1\}) & \mbox{ if  } & \alpha=\beta > -\frac{1}{2} \\
  d(s)  & \mbox{ if  } &  \alpha=\beta = - \frac{1}{2} \\
\end{array}\right.
\end{eqnarray*}

with $d(s)=\int_{D_s}T_s(x) d\mu(x)$.

\end{Example}

The example indicates that the function $M(s)$ contains little
structural information. We, therefore, present a second concept,
which is closer to the philosophy of $P_n$-weakly stationarity.


\begin{Definition}
A positive definite kernel $K$ defined on a polynomial hypergroup
given by $(P_n)_{\nin}$ is called H-asymptotically $P_n$-weakly
stationary if
$$\lim_{n \to \infty} \frac{1}{\sum_{k=0}^n h(k)} \sum_{k=0}^n
(K(k*s,0)-K(k,s))\;h(k) =0  \quad \mbox{for all  }s \in {\Bbb
N}_0.$$
\end{Definition}

We will give later a condition for a H-asymptotically $P_n$-weakly
stationary kernel to be in a class defined in the next section.

\subsection{$P_n$-harmonizability and Wiener theorems}
We now introduce kernels which are strongly harmonizable with
respect to polynomial hypergroups. Wiener type theorems will be
proved, that show connections of discrete parts of the spectral
measure with its Fourier transform.

\begin{Definition} A positive definite kernel on a polynomial hypergroup
given by $(P_n)_{\nin}$ is called $P_n$-harmonizable if $K(m,n)$
admits a representation
\begin{equation}\label{hypHarm}
    K(m,n)=\int_{D_s \times D_s}P_m(x)P_n(y)d\mu(x,y)
\end{equation}
for all $m,n\in\mathbb{N}_0$ where $\mu$ (spectral measure) is a
complex Borel measure on $D_s\times D_s$.
\end{Definition}
Note, that the positive definiteness of $K$ imposes additional
conditions on $\mu$.
\begin{Proposition} A $P_n$-weakly stationary kernel $K$ is also
$P_n$-harmonizable.
\end{Proposition}

\begin{Proof}
Let $K(n,m)=\int_{D_s}P_m(x)P_n(x)d\mu(x)$ be the spectral
representation of the $P_n$-weakly stationary kernel $K$. For
$A\subseteq D_s$ define \mbox{$\Delta(A):=\{(x,x)\in D_s\times
D_s:x\in A\}$} and a measure $\tilde{\mu}$ which is concentrated
on $\Delta(D_s)$ with the property
\[\tilde{\mu}(\Delta(A))=\mu(A)\hspace{1cm} \mbox{ for all Borel sets } A\subseteq D_s.\]
Then
$$K(m,n)=\int_{D_s }P_m(x)P_n(x)d\mu(x)=\int_{D_s \times D_s}P_m(x)P_n(y)d\tilde{\mu}(x,y),$$
and thus $K$ is $P_n$-harmonizable.
\end{Proof}

The classical theorem of N. Wiener characterizes the discrete part
of a complex Borel measure $\mu \in M(T)$ on the torus group T in
the formula
\[ \sum_{z \in T} |\mu({z})|^2 = \displaystyle \lim_{n \to \infty} \frac{1}{2n+1} \sum_{k=-n}^n
|\hat{\mu}(k)|^2 ,\] where $\hat{\mu}$ denotes the
Fourier-Stieltjes transform of $\mu$; see, e.g., (\cite{Graham},
p. 415). This implies in particular $\{z \in T \;| \; |\mu(\{z\})|
\} = \varnothing$ if and only if
$$\displaystyle \lim_{N \to \infty} \frac{1}{2n+1} \sum_{k=-n}^n |\hat{\mu}(k)|^2 = 0 .$$

Now, let $(p_n)_{\nin}$ be an orthonormal polynomial sequence on
the real line with respect to a probability measure $\pi$ with
$|\supp \pi| = \infty$. The polynomials are assumed to be real
valued with $\mbox{deg }(p_n)=n$ and positive leading coefficients
$\gamma_n$. Note, that for the following theorem we do not require
that the polynomials induce a hypergroup.

Also, let $A$ be an arbitrary compact subset of $\Bbb R$ and $A^N$
its N-fold cartesian product. The following notions depend on $A$,
but we omit the prefix $A$. Denote $m(n)=\max \{|p_n(x)|\;|\;x \in
A \}$ and set
\begin{equation}\label{normalization}
    P_n(x)=\frac{p_n(x)}{m(n)}.
\end{equation}

For a complex Borel measure $\mu_N$ on $A^N$ define the
Fourier-Stieltjes transform by
\begin{equation}\label{muHat}
\widehat{\mu_N}
(k_1,\ldots,k_N)=\displaystyle{\frac{1}{m(k_1)\cdot \ldots \cdot
m(k_N)} \int_{A^N} p_{k_1}(x_1) \cdot \ldots \cdot p_{k_N}(x_N)
 d\mu_n(x_1,...,x_N)}
\end{equation}
for all $k_1,\ldots,k_N \in {\Bbb N}_0$. We set for $x,y \in A$
\begin{equation}\label{Tfun}
t_n(x,y):=\displaystyle{\frac{1}{\sum_{k=0}^n m^2(k)}
\sum_{k=0}^{n} p_k (x) p_k (y)},
\end{equation}
\begin{equation}
S_N:=\{x_1,\ldots, x_N \in A \;|\; \limsup_{n \to \infty}\; t_n(x_1,x_1)\cdot \ldots \cdot t_n(x_N,x_N)\; > 0\}
\end{equation}
and define the condition

\hspace{1 cm}(G) \hspace{2.5
cm}$\displaystyle{\lim_{n \to \infty} \;
\frac{\gamma_n}{\gamma_{n+1}} \; \frac{m(n) \; m(n \! + \!
1)}{\sum_{k=0}^n m^2(k)}} \; = \; 0$.

Now, we can state the following theorems.

\begin{Theorem} \label{WienerGeneral}
Let $(p_n)_{\nin}$ be a orthogonal polynomial
system and $A$ a compact real set such that (G) holds. Further let
$\mu_N$ be a complex Borel measure on $A^N$ with discrete part
$\mu_{N d}$. Then the following conditions are equivalent:
\begin{eqnarray*}
(i) &  S_N \cap \mbox{supp }\mu_{Nd}=\emptyset \\
(ii) & \displaystyle{\lim_{n \to \infty} \frac{1}{(\sum_{k=0}^n
m^2(k))^N} \sum_{k_1,\ldots,k_N=0}^{n} |
\widehat{\mu_N}(k_1,...,k_N) m(k_1) \cdot... \cdot m(k_N)|^2}=0
\end{eqnarray*}
\end{Theorem}

\begin{Proof}
Denote in $A^N \! \! \times \!\! A^N$ the set
$$\underline{\Delta}_N:=\{(x_1,...,x_N) \times (x_1,...,x_N)|
\; x_1,...,x_N \in A \}.$$ From the definitions we get
\begin{eqnarray*}
&\mbox{ }&\displaystyle{\frac{1}{(\sum_{k=0}^n m^2(k))^N}
\sum_{k_1,...,k_N=0}^{n}| \widehat{\mu_N}(k_1,...,k_N)\cdot
m(k_1) \cdot ... \cdot m(k_N)|^2}\\\\
&=& \displaystyle{ \int_{A^N \! \times \! A^N}
\sum_{k_1,...,k_N=0}^{n}\frac{ p_{k_1}(x_1)p_{k_1}(y_1)...
p_{k_N}(x_N) p_{k_N}(y_N)}{(\sum_{k=0}^n m^2(k))^N}\; d \mu_N \!
\times \! \overline{\mu_N}
(x_1,...,x_N,y_1,...,y_N)}\\\\
&=& \displaystyle{ \int_{A^N \! \times \! A^N \setminus
\underline{\Delta}_N} t_n(x_1,y_1)... t_n(x_N,y_N)\;
d \mu \times \overline{\mu_N}(x_1,...,x_N,y_1,...,y_N)}\\\\
& &\displaystyle{+\int_{\underline{\Delta}_N}t_n(x_1,y_1)\cdot ...
\cdot t_n(x_N,y_N) \; d \mu_N \! \times \!\!
\overline{\mu_N}(x_1,...,x_N,y_1,...,y_N)}.
\end{eqnarray*}

For $x \neq y$ the identity of Christoffel-Darboux gives
$$t_n(x,y)=\frac{\gamma_n}{\gamma_{n+1}}\;\frac{m(n+1)m(n)}
{\sum_{k=0}^n m^2(k)}\;
\frac{P_{n+1}(x)P_n(y)-P_n(x)P_{n+1}(y)}{x-y}.$$

We have $|t_n(x_1,y_1)\cdot ... \cdot t_n(x_N,y_N)| \leq 1$ for
any $x_1,...,x_N,y_1,...,y_N \in A$ and \\
$\lim_{n \to \infty}|t_n(x_1,y_1)\cdot ... \cdot t_n(x_N,y_N)|=0 $
for $(x_1,...,x_N,y_1,...,y_N )\in A^N \! \times \! A^N \setminus
\underline{\Delta}_N$ by condition (G).

Thus, the first integral of the last equality tends to zero with
growing n by the dominated convergence theorem. For the second
integral we obtain with Fubini's theorem
$$\displaystyle{ \int_{\underline{\Delta}_N}t_n(x_1,y_1)\cdot ...
\cdot t_n(x_N,y_N) \; d \mu_N \! \times \!\!
\overline{\mu_N}(x_1,...,x_N,y_1,...,y_N)}$$
$$ = \displaystyle{
\int_{A^N}t_n(x_1,x_1)...t_n(x_N,x_N)
\overline{\mu}(\{(x_1,...,x_N)\})d\mu(x_1,...,x_N)}.$$

Hence, for n to infinity one has
$$\displaystyle{\frac{1}{(\sum_{k=0}^n m^2(k))^N}
\sum_{k_1,...,k_N=0}^{n}| \widehat{\mu_N}(k_1,...,k_N)\cdot m(k_1)
\cdot ... \cdot m(k_N)|^2}$$
$$= \displaystyle{ \sum_{(x_1,...,x_N) \in supp \, \mu_{N d}} t_n(x_1,x_1)...t_n(x_N,x_N) |\mu
\{(x_1,...,x_N)\}|^2+o(1)}$$
Since $t_n(x_1,x_1)...t_n(x_N,x_N) \;
\geq 0$ condition (ii) implies (i).

Conversely, Lebesgue's theorem of dominated convergence applied to
the discrete measure
$$\nu_N:= \sum_{(x_1,...,x_N) \in supp \, \mu_{N d}} |\mu_N
\{(x_1,...,x_N)\}|^2 \delta _{(x_1,...,x_N)}$$

shows that (i) implies (ii).

\end{Proof}

The above theorem fits for $N=2$ the case of $P_n$-harmonizable
kernels $K(n,m)$. Here $\hat{\mu}_2(n,m)= K(n,m)$, $\mu_2$ is the
spectral measure and $A^2$ is the set $D_s \times D_s$.

\begin{Remark} \label{SWiener}
The case $N=1$ has been investigated in \cite{wie}.

There, we examined the set $S=S_1$ and the growth condition (G)
for several classical polynomial systems. For example, we found
for Jacobi polynomials $P_n^{(\alpha,\beta)}$ that (G) is always
valid and $S$ is given by

\begin{eqnarray*}
S= \left\{ \begin{array}{lll}
 ]-1,1[  & \mbox{ if  } &  -1< \alpha,\beta < - \frac{1}{2} \\
 ]-1,1]  & \mbox{ if  } & -1< \beta < \alpha = - \frac{1}{2} \\
  \mbox{$[-1,1[$} & \mbox{ if  } &  -1< \alpha < \beta = - \frac{1}{2} \\
  \mbox{$[-1,1]$}  & \mbox{ if  } &  \alpha=\beta = - \frac{1}{2} \\
\{-1,1\} &  \mbox{ if  }  & \alpha=\beta > -\frac{1}{2} \\
\{1\}    &  \mbox{ if  }  & \alpha > \beta, \alpha > -\frac{1}{2} \\
\{-1\}   &  \mbox{ if  }   & \alpha < \beta, \beta > -\frac{1}{2}\\
\end{array}\right.
\end{eqnarray*}
\end{Remark}

In the situation of $P_n$-weakly stationary kernels we give an
extended version of the theorem. With $m(n)^2=h(n)$, the Haar
measure, the following holds true:

\begin{Theorem} \label{satz21}
Let $K(x,y)$ be weakly stationary with respect to $(P_{n})_{n \in
\nz_0}$ with spectral measure $\mu$ having discrete part $\mu_d$.
Assume that $(P_{n})_{n \in \nz_0}$ define a polynomial hypergroup
with condition $(H)$. Then, the following is equivalent:

\renewcommand{\labelenumi}{\alph{enumi})}

\begin{enumerate}
\item $S \cap supp \mu_d = \emptyset  $;\\[-0.4cm]
\item $\displaystyle \lim_{n\to \infty}
 \frac{1}{\sum_{k=0}^{n}h(k)}\sum_{k=0}^n {|K(k,0)| h(k)}=0  $;\\
\item $\displaystyle \lim_{n\to \infty}
 \frac{1}{\sum_{k=0}^{n}h(k)}\sum_{k=0}^n {|K(k,0)|^2 h(k)}=0  $;\\
\item $\displaystyle \lim_{n\to \infty}
 \frac{1}{( \sum_{k=0}^{n}h(k))^2}\sum_{k,l=0}^n {|K(k,l)|^2 h(k) h(l)}=0$.\\
\end{enumerate}
\end{Theorem}

\begin{Proof}

The equivalence of a) and c) has already been shown.

a) equivalent to d):

$\displaystyle \frac{1}{(\sum_{k=0}^n h(k))^2} \sum_{k,l=0}^n |K(k,l)|^2 h(k) h(l)\\[0.3cm]
\hspace*{0.4cm} = \int_{D_s} \int_{D_s} \Bigl(
\frac{1}{(\sum_{k=0}^n h(k))^2} \sum_{k,l=0}^n
 P_k(x) P_l(x) P_k(y) P_l(y) h(k) h(l) \Bigr) d \mu (x) d \overline{\mu} (y) \\[0.3cm]
\hspace*{0.4cm} =\int_{D_s} \int_{D_s} \Bigl(
\frac{1}{(\sum_{k=0}^n h(k))^2} (\sum_{k=0}^n
 P_k(x) P_k(y) h(k))^2 \Bigr)  d \mu (x) d \mu (y) \\[0.3cm]
\hspace*{0.4cm} = \int_{D_s \times D_s} t^2_n (x,y) d \mu \times  \overline{\mu }(x,y)$.\\[0.3cm]

From here one proceeds with the same arguments as in the preceding
theorem with $t_n^2$ instead of $t_n$.

The equivalence of b) and c) follows from the inequalities\\[-0.5cm]
\begin{eqnarray}
\Bigl(\frac{1}{\sum_{k=0}^n h(k)} \sum_{k=0}^n |K(k,0)| h(k)
\Bigr)^2 & \leq &
\frac{1}{\sum_{k=0}^n h(k)} \sum_{k=0}^n |K(k,0)|^2 h(k) \nonumber\\[0.2cm]
& \leq& \frac{K(0,0)}{\sum_{k=0}^n h(k)} \sum_{k=0}^n | K(k,0))|
h(k).\hspace{3cm} \label{nr3}
\end{eqnarray}
The first inequality is derived from the Cauchy-Schwarz inequality\\
\begin{equation} \label{nr4}
 \sum_{k=0}^n | K(k,0)| h(k) \leq (\sum_{k=0}^n K(k,0)^2 h(k))^{\frac{1}{2}}
( \sum_{k=0}^n h(k))^{\frac{1}{2}},
\end{equation}
by squaring both sides of (\ref{nr4}) and dividing by
\mbox{$\displaystyle( \sum_{k=0}^n h(k))^2 $}. \\[0.2cm]
The second inequality in (\ref{nr3}) follows from \quad \mbox{$
K(0,0) \geq | K(k,0) |,$} \quad \mbox{$\forall \: k \in \nz_0$},
which is valid for the function $K(.,0)$, being positive definite
with respect to the hypergroup convolution (compare also
\cite{jew}).

\end{Proof}

For the proof that b) implies c) we used the positive definiteness
of $K(.,0)$. All other implications could be verified without
assuming hypergroup structure. Here, the condition (H) is the
substitute for condition (G).

If $\mu$ has no discrete part, then, of course, all statements of
the equivalence are true.

The following corollary for Chebyshev polynomials of first kind
$T_n$ is required in a later chapter, where we give consistent
estimators of covariance structures.

\begin{Corollary}\label{WienerCheby}
If the kernel $K$ in Theorem \ref{satz21} is weakly stationary
with respect to Chebyshev polynomials of first kind $T_n$ then the
following statements are equivalent:

\[
\begin{array}{ll}
\begin{array}{lcl}
&  &  \\
a)\displaystyle{\mbox{ }\mu \mbox{ has no atoms}}, &  &  \\[0.4cm]
b)\displaystyle{\lim_{n\to \infty }\frac
1n\sum_{k=0}^{n-1}|K(k,0)|^2(k)=0}, &  &
\\[0.2cm]
&  &
\end{array}
&
\begin{array}{lcl}
c)\displaystyle{\lim_{n\to \infty }\frac
1{n^2}\sum_{k,l=0}^{n-1}|K(k,l)|^2 =0
}, &  &  \\[0.2cm]
d)\displaystyle{\lim_{n\to \infty }\frac
1n\sum_{k=0}^{n-1}|K(k,0)|=0}. &  &
\\[0.2cm]
&  &
\end{array}
\end{array}
\]
\end{Corollary}

\begin{Proof}
This is immediate from $h(n)=2$ for all $n >0$ and $S=D_s=[-1,1]$
(see \cite{Bloo}).
\end{Proof}

\subsection{Harmonizability and cyclostationarity}
In the classical situation, Gladyshev showed that discrete
cyclostationary sequences are always harmonizable. He essentially
used Fourier techniques which are not available in our situation.
What we will do instead is to characterize cyclostationary
processes that are harmonizable.

First, recall our elementary example of cyclostationarity:

$$K(n,m)=\cases{(C-1)^2, & for n and m odd;\cr
(C+1)^2,& for n and m even \cr (C-1)(C+1),& else.\cr}$$

$K$ is harmonizable: The crucial equation
$$K(n,m)=\int_{D_s \times D_s}T_n(x)T_m(y)d\mu(x,y)$$
holds with the spectral measure
$$\mu(x,y)=\bigl(C \delta_{1}(x)+ \delta_{-1}(x)\bigr) \times \bigl (C \delta_{1}(y)+ \delta_{-1}(y) \bigr).$$
This is easily verified with the symmetries $T_n(1)=1$ and
$T_n(-1)=(-1)^n$ for all $\nin$.

Our main theorem about harmonizable and T-cyclostationary kernels
requires a feature of complex (regular) Borel measures which is
also used in the classical theory of harmonizable processes. Since
a proof is not easily found in literature, we provide one tailored
to our needs in some detail.
\begin{Lemma}
Let $(P_n)_{\nin}$ be a sequence of real polynomials with $
\mbox{deg }(P_n) =n$ and $\mu$ a complex Borel measure on a
compact set $K \subset \Bbb R^2$. For any real polynomial $S$ the
conditions
$$ \int_K P_n(x)P_m(y) S(x,y) d\mu(x,y)=0 \mbox{  for all }n,m \in {\Bbb N}_0
\qquad \mbox{and}$$
$$ \mbox{supp } \mu \subseteq \{(x,y)\in K \;|\; S(x,y)=0)\} $$
are equivalent.
\end{Lemma}
\begin{Proof}
Since $deg (P_n)=n$ any monomial $x^k y^l$ and thus any complex
polynomial $Q$ can be represented as a (complex) linear
combination of terms of the form $P_n(x)P_m(y)$. So,
$$ \int_K P_n(x)P_m(y) S(x,y) d\mu(x,y)=0 \mbox{  for all }n,m \in {\Bbb N}_0$$
implies for any complex polynomial $Q$
$$ \int_K Q(x,y) S(x,y) d\mu(x,y)=0.$$
From Stone-Weierstra{\ss} we know that the complex polynomials are
dense in $C(K,\Bbb C)$, the space of continuous complex valued
functions on K, with respect to the supremum norm $\|.\|_\infty$.
As the integrals of an uniform convergent sequence of continuous
functions on the compact $K$ converge, we get by
$$L(f):=\int_K f(x,y) S(x,y) d\mu(x,y) \mbox{  for } f \in C(K,\Bbb
C)$$ a continuous linear functional on $C(K,\Bbb C)$, which is
zero. The Riesz representation theorem provides a unique measure
$\sigma$ for the representation $$L(f)=\int_K f(x,y)
d\sigma(x,y).$$ Then, $\sigma=S\;\mu$ has to be the zero measure.

It remains to show that the restriction of $\mu$ to $K \cap \{
(x,y) \;|\; S(x,y)\neq 0 \}$ gives the zero measure. As $\mu$ is a
complex measure we need some additional work.

The complex measure $\mu$ has the decomposition $\mu = \mu_R+i
\mu_I$ with (finite) signed measures $\mu_R$ and $\mu_I$.
Therefore, $S \mu_R$ and $S \mu_I$ have to be zero measures. For
this reason, it is enough to do the proof for signed measures.

Let now $\mu$ be a signed measure. The Hahn-Jordan decomposion
$$\mu=\mu^+ - \mu^-$$
provides positive measures $\mu^+$ and $\mu^-$ together with a
decomposition $K=A \cup B$ of $K$ into disjoint measurable sets
$A$ and $B$ such that $\mu^+(E)=\mu (E \cap A)$ and
$\mu^-(E)=-\mu(E \cap B)$ hold for any measurable set $E \subset
K$.

The equation
$$0= \sigma(E) = \int_{A \cap E} S(x,y) d\mu^+(x,y) - \int_{B \cap E} S(x,y) d\mu^-(x,y)$$
then shows ($\sigma(A \cap E)= \sigma(B \cap E)=0$) that both
integrals on the right hand side are zero and we can further
restrict the work to positive measures.

Now we are almost done:

The zero set of $S$ partitions $K$ in connected components with
constant sign of $S$. Let M be a component where $S$ is nowhere
zero and $K_M$ a compact subset of M. The continuous function
$|f|$ has a minimum $\epsilon >0$ in $K_M$.

One has for $S$ positive in $K_M$
$$0=\int_{K_M}S(x,y)d\mu^+(x,y) \geq \epsilon \int_{K_M}d\mu^+(x,y)$$
and likewise
$$0=\int_{K_M}S(x,y)d\mu^+(x,y) \leq - \epsilon \int_{K_M}d\mu^+(x,y)$$
for $S$ negative in $K_M$. In any case $\mu^+(K_M)=0$ follows and
from the (inner) regularity of $\mu^+$ also $\mu^+(E_M)=0$ for all
measurable sets $E_M$ in $M$. This is valid for any component M
and in the same manner for the measures $\mu^-$, $\mu_R$, $\mu_I$.
Finally we get $\mu (E) =0$ for any measurable set in $$K \cap \{
(x,y) \;|\; S(x,y)\neq 0 \}.$$ This proves one direction.

The other direction is trivial.

Note, that we used for the support of the complex measure $\mu$
the definition $\supp \mu =\supp \mu_R^+ \cup \supp \mu_R^- \cup
\supp \mu_I^+ \cup \supp \mu_I^- $
\end{Proof}
\newpage
\begin{Theorem} \label{harmcyclo}
Let $K$ be a harmonizable kernel with respect to a hypergroup with
polynomial system $(P_n)_{\nin}$. $K$ is T-cyclostationary if and
only if
$$\mbox{supp } \mu \subseteq \{(x,y) \in D_s \times D_s \;|\;
P_T(x)=P_T(y)\}.$$
\end{Theorem}

\begin{Proof}
If $K$ is harmonizable one has in general for all $n,m,l \in {\Bbb
N}_0$
\begin{eqnarray*}
K(n*l,m)&=& \sum_{k=|n-l|}^{n+l} g(n,l,k) K(k,m)
=\sum_{k=|n-l|}^{n+l} g(n,l,k) \int_{D_s \times
D_s}P_k(x)P_m(y)d\mu(x,y)\\
&=& \int_{D_s \times D_s}P_n(x)P_l(x)P_m(y)d\mu(x,y).
\end{eqnarray*}

Thus, for $K$ also T-cyclostationary
\begin{eqnarray*}
\int_{D_s \times D_s}P_T(x)P_n(x)P_m(y)d\mu(x,y)-\int_{D_s \times
D_s}P_T(y)P_n(x)P_m(y)d\mu(x,y)=0
\end{eqnarray*}
holds for all $n,m \in {\Bbb N}_0$. The above lemma now shows
$$ \mbox{supp } \mu \subseteq \{(x,y) \in D_s \times D_s \;|\; P_T(x)=P_T(y)\}.$$

Conversely, if the spectral measure $\mu$ of a harmonizable kernel
$K$ has the mentioned support,
$$K(n*T,m)-K(n,m*T)=\int_{D_s \times D_s}P_n(x)P_m(y) \bigl (P_T(x)-P_T(y) \bigr)d\mu(x,y)=0 $$
is valid for all $n,m \in {\Bbb N}_0$.
\end{Proof}

\newpage
As in the classical case, we show the maximal support of the
 spectral measure for kernels being $5-cylostationary$ with respect to Chebyshev polynomials of the first kind.

\begin{center}
\includegraphics[width=8cm,keepaspectratio]{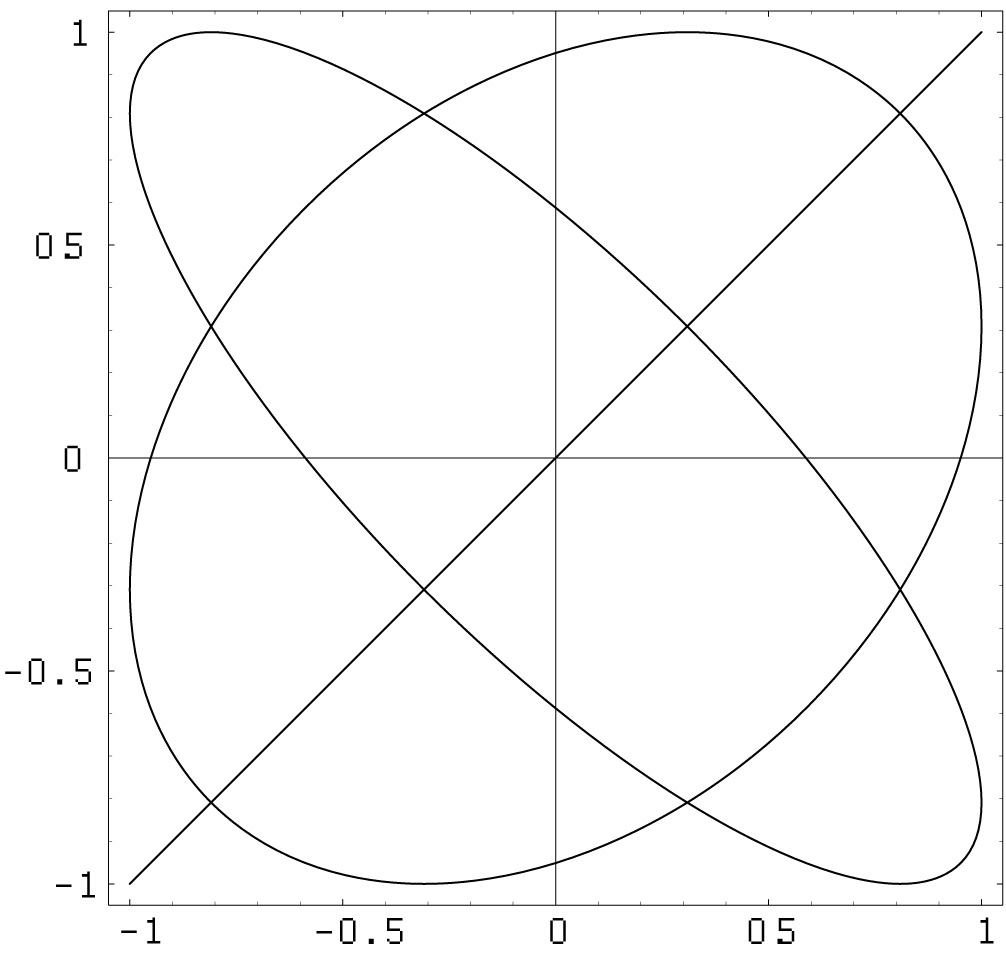}
\end{center}

\subsection{Harmonizability and asymptotic stationarity}
In the classical case all strongly harmonizable kernels are
asymptotically stationary \cite{Rozanov}. In the scenario of
polynomial hypergroups with Haar measure satisfying condition (H)
we can identify the $P_n$-harmonizable kernel that are
H-asymptotically $P_n$-weakly stationary.

\begin{Theorem} \label{harmasymp}
Let the Haar weights of a polynomial hypergroup satisfy condition
(H) and let a kernel $K$ be $P_n$-harmonizable with respect to
this hypergroup. Then, $K$ is H-asymptotically $P_n$-weakly
stationary if and only if the spectral measure $\mu$ satisfies for
all $s \in {\Bbb N}_0$
$$\int_{\widetilde{\Delta}} P_s(y) d\mu(x,y)=\mu(\widetilde{\Delta})$$
with $\widetilde{\Delta}:= \{(1,y)\;| \;y \in D_s, y \neq 1 \}.$
\end{Theorem}

\begin{Proof}
\begin{eqnarray*}
M_s(n)&:=&\frac{1}{\sum_{k=0}^n h(k)} \sum_{k=0}^n
(K(k*s,0)-K(k,s))h(k)\\\\
&=& \frac{1}{\sum_{k=0}^n h(k)} \sum_{k=0}^n \Bigr(\int_{D_s
\times
D_s} (P_k(x)P_s(x)-P_k(x)P_s(y)) \;d \mu(x,y)\Bigl)h(k)\\\\
&=&  \int_{D_s \times D_s}\frac{\sum_{k=0}^n
P_k(x)h(k)}{\sum_{k=0}^n h(k)} (P_s(x)-P_s(y))d\mu(x,y)\\\\
&=&  \int_{D_s \times D_s \;\backslash
\widetilde{\Delta}}\frac{\sum_{k=0}^n P_k(x)h(k)}{\sum_{k=0}^n
h(k)} (P_s(x)-P_s(y))d\mu(x,y)+
\int_{\widetilde{\Delta}}(P_s(x)-P_s(y))d\mu(x,y)\\\\
&=&\int_{\widetilde{\Delta}}(P_s(x)-P_s(y))d\mu(x,y)+ o(1)
\end{eqnarray*}

The last equality can be proved as follows:

For $x \neq y $ the Christoffel-Darboux formula
(\ref{Christoffel1}) gives
\begin{equation}
    \frac{\sum_{k=0}^n P_k(x)P_k(y) h(k)}{\sum_{k=0}^n
h(k)}=\frac{h(n)}{\sum_{k=0}^n h(k)}a_n a_0 \frac{P_{n+1}(x)
P_n(y)-P_n(x)P_{n+1}(y)}{x-y}
\end{equation}
Thus, for $x \neq 1$ and $y=1$ ($P_n(1)=1$), condition (H) implies
\begin{equation}\label{ChristConf}
    \lim_{n \to \infty}\; \frac{\sum_{k=0}^n
P_k(x)h(k)}{\sum_{k=0}^n h(k)}=0
\end{equation}
and by the dominated convergence theorem  one gets
$$\lim_{n \to \infty} \; \int_{D_s \times D_s \;\backslash
\widetilde{\Delta}}\frac{\sum_{k=0}^n P_k(x)h(k)}{\sum_{k=0}^n
h(k)} (P_s(x)-P_s(y))d\mu(x,y)=0.$$ Note, that at $(x,y)=(1,1)$
the integrand is always zero. So we have the equation
\begin{equation}\label{harmasym}
    \lim_{n \to \infty} M_s(n)=\int_{\widetilde{\Delta}}(1
-P_s(y))d\mu(x,y),
\end{equation}
showing that $\displaystyle{\lim_{n \to \infty} M_s(n)=0}$ is
valid if and only if $\int_{\widetilde{\Delta}} P_s(y)
d\mu(x,y)=\mu(\widetilde{\Delta})$ holds.


\end{Proof}

\begin{Corollary} $ $
\begin{enumerate}
\item If $\mu$ has no mass in $\widetilde{\Delta}$, K is
H-asymptotically stationary.

\item If K is H-asymptotically stationary and
$\mu(\widetilde{\Delta})=0$, $\mu$ has no mass in
$\widetilde{\Delta}$.
\end{enumerate}
\end{Corollary}

\begin{Proof}
\begin{enumerate}
\item If $\mu$ has no mass in $\widetilde{\Delta}$, the condition
$\int_{\widetilde{\Delta}}(1 -P_s(y))d\mu(x,y)=0$ is obviously
satisfied for all $s \in {\Bbb N}_0$.

\item $\int_{\widetilde{\Delta}}P_s(y) d\mu(x,y)=0$ for all $s \in
{\Bbb N}_0$, implies by the Riesz representation theorem, that on
$\widetilde{\Delta}$, $\mu$ is the zero measure.
\end{enumerate}
\end{Proof}
\begin{Remark}\label{RemAsHarm} $ $
\begin{enumerate}
    \item The spectral measure $\mu$ can well have mass on $\widetilde{\Delta}$
    and K still be H-asymptotically stationary. Thereto, the
    modified momentum problem
     $$\int_{\widetilde{\Delta}} P_s(y) d\mu(x,y)=\mu(\widetilde{\Delta})$$
     for all $s \in {\Bbb N}_0$, has to be solved. If $D_s \subseteqq [-1,1]$
     and $\mu(\widetilde{\Delta})>0$, there is a
     positive Borel measure solving the problem, as has been
     shown in \cite{LaMM}.
    \item On the other hand, if $\mu(1,y)=\delta_{y_0}$, with one
    atom $y_0 \neq 1$ only, $P_s(y_0)=1$ cannot be satisfied for
    all  $s \in {\Bbb N}_0$, as the unit character is given by $1 \in
    D_s$. Thus the associated kernel is not H-asymptotically
    stationary.
    \item Note, that because of the symmetry
    $d\mu(x,y)=d\overline{\mu(y,x)}$,
    conditions for $\widetilde{\Delta}$ also apply to $\{(x,1)\;| \;x \in D_s, x \neq 1 \}$.
\end{enumerate}
\end{Remark}

\chapter{Stochastic sequences}
The preceding discussion of positive definite kernels will now be
linked to stochastic sequences. Aronszajn's theorem shows that
every positive definite kernel can be regarded as covariance
kernel of some second order stochastic process. We consider
stochastic sequences having kernels as discussed before. Thereto,
we assume that the random variables $X_n$ of the stochastic
sequence are elements of $L_0^2(P)$, the Hilbert space of
(complex) centered square integrable random variables on a given
probability space $(\Omega, \Sigma, P)$ and denote the sequence
according to the respective kernel. To get consistency with the
existing literature, we use in the following the notations
$d(n,m):=K(n,m)$ and $d(n):=K(n,0)$.

\section{Spectral representations and dilations}
In this section we present some general properties of
$P_n$-harmonizable and $P_n$-weakly stationary sequences. Two
classes of harmonizable sequences which need not to be
$P_n$-weakly stationary are introduced. One class has discrete
spectral measure, the other emerges from truncated sequences and
has spectral measures continuous with respect to $d\pi \times
d\pi$.

\subsection{Harmonizable sequences}
The following theorems on spectral representations and dilations
hold  for all weakly harmonizable processes indexed by abelian
hypergroups (Rao \cite{RaoBiH} and Leitner \cite{Leihyp}).

\begin{Theorem}\label{harmDil}
Let $(X_n)_{\nin} \subset H = L_0^2(P)$ be a $P_n$-harmonizable
sequence with spectral measure $\mu$, then the following is valid:
\begin{enumerate}
\item[(i)] There exists a unique stochastic measure Z, mapping
Borel sets of $D_s$ to elements of $L_0^2(P)$, such that
\begin{equation}\label{HarmSequSpec}
    X_n = \int_{D_s} P_n(x) dZ(x) \quad \mbox{ for all } \nin,
\end{equation}
where $E(Z(A) \overline{Z(B)}) = \mu(A,B).$
 \item[(ii)] $(X_n)_{\nin}$ is a norm bounded sequence and
\begin{equation}\label{normbounded}
    \parallel \sum_{k=0}^{\infty}f(k)X_k h(k) \parallel \leq M \parallel \hat{f}
    \parallel,\qquad f \in l^1(h)
\end{equation}
with some constant M and $\hat{f}(x)=\sum_{k=0}^{\infty}f(k)P_k(x)
h(k)$ for $x \in D_s$.
 \item[(iii)] There exists an extension Hilbert space $\tilde{H}=L_0^2(\tilde{P})
 \supseteqq H$ on an enlarged probability space
 $(\tilde{\Omega},\tilde{\Sigma},\tilde{P})$ and a $P_n$-weakly stationary sequence
 $(\tilde{X}_n)_{\nin} \subset \tilde{H}$ such that
 \begin{equation}\label{dilation}
    X_n = pr_H \tilde{X}_n, \quad \nin,
\end{equation}
where $pr_H: \tilde{H} \rightarrow H$ is the orthogonal projection
from  $\tilde{H}$ to H.
\end{enumerate}
\end{Theorem}

\begin{Proof}
The proof is rather involved and works in a more general
situation. We will not replicate it here and refer instead to Rao
\cite{RaoBiH}, sec 2.5, Th.3, Leitner \cite{Leihyp}, Th.1 and also
Niemi \cite{Niemi}.
\end{Proof}

\begin{Remark}
Actually, it is proven (\cite{Leihyp},Theorem 1), that the
statements of the above theorem are all equivalent to the random
sequence being weakly harmonizable. The sequence
$(\tilde{X}_n)_{\nin} \subset \tilde{H}$ is called dilation of
$(X_n)_{\nin}$. The existence of dilations can be used to derive
some sort of shift operators for $P_n$-harmonizable sequences. Let
us recall the classical situation: For a (classical) weakly
stationary sequence $(Y_n)_{n \in \Bbb Z} $ there exists a unitary
shift operator $U$ with $U Y_n= Y_{n+1}$ and $Y_n=U^n Y_0$. For
$P_n$-weakly stationary sequences $(X_n)_{\nin}$ shift operators
 are defined by (compare \cite{lalest})
$$\displaystyle T_n(X_m) := X_{n*m}=\sum_{s=|m-n|}^{n+m} g(n,m,s) X_s.$$
Linear continuation of $T_n$ in $H$, the closed linear span of
$(X_n)_{\nin}$ in $L^2(P)$, yields operators with $\parallel T_n
\parallel \leq 1$, which commute, are symmetric and satisfy
\begin{equation}\label{shifts}
    T_n T_m = T_{n*m}=\sum_{s=|n-m|}^{n+m} g(n,m,s) T_s.
\end{equation}
The sequence $(X_n)_{\nin}$ is again determined by $X_0$: $X_n=
T_n X_0$. Further, any start vector $X_0 \in H$ induces in this
way a $P_n$-weakly stationary sequence.

The harmonizable case is more involved. In general, there exist no
canonical shift operators. But, one can assign sequences of
contractive operators, defined on an enlarged Hilbert space, which
induce the random sequence from a start value. Details are given
in the following theorem, which also holds for weakly harmonizable
processes indexed by abelian hypergroups (\cite{Leihyp}, Th. 3).
\end{Remark}

\begin{Theorem}\label{OpShift}
Let $(X_n)_{\nin} \subset H = L_0^2(P)$ be a $P_n$-harmonizable
sequence. There exists an extension Hilbert space
$\tilde{H}=L_0^2(\tilde{P}) \supseteqq H$ on an enlarged
probability space $(\tilde{\Omega},\tilde{\Sigma},\tilde{P})$, a
random variable $\tilde{X}_0 \in \tilde{H}$ and a sequence
$(V_n)_\nin$ of bounded linear operators from $\tilde{H}$ to $H$
with
$$X_n= V_n \tilde{X}_0 \qquad \mbox{ for all } \nin .$$
The restricted sequence $(S_n)_\nin$ with $S_n=V_n |_H$ satisfies
$S_0=id$, $S^*_n=S_n$ and
\begin{equation}\label{operatorshift}
    0 \! \leq \! \sum_{m, \nin} \sum_t E(S_t W_n, W_m)
\varepsilon_m \! \ast \! \varepsilon_d \! \ast \! \varepsilon_d \!
\ast \! \varepsilon_n (t) \! \leq \! \sum_{m, \nin} \sum_t E(S_t
W_n, W_m) \varepsilon_m \! \ast \! \varepsilon_n (t),
\end{equation}
for all $d \in {\Bbb N}_0$ and every sequence $(W_n)_{\nin}
\subset H$ with $W_n=0$ for almost all $\nin$.
\end{Theorem}

\begin{Proof}
According to Theorem \ref{harmDil}(iii), there exists a Hilbert
space $\tilde{H} \supseteqq H$ and a $P_n$-weakly stationary
sequence
 $(\tilde{X}_n)_{\nin} \subset \tilde{H}$ such that
    $X_n = pr_H \tilde{X}_n$ for all $\nin$. If $(T_n)_{\nin}$ are the
    translation operators of $(\tilde{X}_n)_{\nin}$, define
    $V_n:= pr_H T_n$ for all $\nin$. Then $X_n=pr_H \tilde{X}_n= pr_H T_n
    \tilde{X}_0 = V_n \tilde{X}_0$ and $S_0=pr_H T_0=id$, since $T_0=id$.
    Further, one has for all $W,W' \in H$
    $$E(S_n W, W')= E(T_n W, pr_H W')=E(pr_H W, T_n W')= E(W, S_n W')$$
    and thus $S^*_n=S_n$. Here we used $T_n=T^*_n$ and that orthogonal
    projections are self adjoint. Finally take a sequence
    $(W_n)_{\nin} \subset H$ with $W_n=0$ for almost all $\nin$. Then, one has
    ´with the properties of $T_n$:

 $ \displaystyle \sum_{m, \nin} \! \sum_t E(S_t W_n, W_m)
\varepsilon_m  \ast  \varepsilon_d  \ast  \varepsilon_d  \ast
\varepsilon_n (t) \!=\!\! \sum_{m, \nin}\!\! E(T_m T_d T_d T_n
W_n,W_m) =
\parallel \sum_{\nin} T_d T_n W_n \parallel \\
\leq \parallel \sum_{\nin} T_n W_n \parallel = \sum_{m, \nin}
E(T_m T_n W_n,W_m)= \! \sum_{m, \nin} \sum_t E(S_t W_n, W_m)
\varepsilon_m \! \ast \! \varepsilon_n (t)$.

This finishes the proof.
\end{Proof}

\begin{Remark}
Note, that one can not avoid extension of $H$, by working with
$S_n$ alone. It might not be possible to have $\tilde{X}_0$ in $H$
with $\tilde{X}_0= X_0$ (compare example \cite{RaoHa} p.301).

The conditions for $S_n$ correspond in the classical group case to
operator families of positive type. Here, the important dilation
theorem of Sz.-Nagy applies ( \cite{Nagy}): Operator families of
positive type can be lifted to families of unitary operators in
extended Hilbert spaces (see also Sz.-Nagy's short historical
review on the topic \cite{Czech}). Therefrom, one gets a converse
of Theorem \ref{OpShift}: A family of positive type operators
defines for any start vector $\tilde{X}_0 \in H$ a weakly
harmonizable process (\cite{RaoHa}, Th.6.3).

For hypergroups an extension of Sz.-Nagy's theorem is valid, which
allows to lift operators of type $S_n$ to shift operators $T_n$ of
$P_n$-weakly stationary sequences (\cite{Leihyp}, Th. 2). The
proof depends on the main dilation theorem for certain operators
indexed by *-semigroups (\cite{RiSz}). As a consequence, one gets,
analogously to the classical case: Operators, with the properties
of $S_n$, define, for any $\tilde{X}_0$ in $H$, hypergroup weakly
harmonizable processes.

\end{Remark}

We now leave the structure theory and come to concrete examples of
sequences.

 The first class of $P_n$-harmonizable sequences, that we present,
has spectral measures consisting of a finite number of atoms.
These are the so called periods of the associated stochastic
sequence. In a later chapter we derive periodograms which allow to
detect this periods from realizations of the stochastic sequence.

\begin{Example}

Consider the sequence
$$X_n := \sum_{k=1}^s P_n (x_k) A_k  \qquad \mbox{with fixed} \quad x_k \in D_s$$
and $A_k$ arbitrary square integrable random variables. The
covariance is
$$EX_n \bar{X}_m=\sum_{k,l=1}^s \sigma_{k l} P_n(x_k) P_m(x_l) \mbox{  with  } \sigma_{k l}:= E A_k \bar{A}_l$$
and the associated spectral measure
$$\mu(x,y) = \sum_{k,l=1}^s \sigma_{k l}\; \delta_{(x_k,x_l)}
(x,y).$$

From the preceding discussions it follows that $(X_n)_{\nin}$

\begin{itemize}
    \item is T-cyclostationary  if and only if
$$ \{ (x_k,y_l)\; | \; k,l=1,...,s\} \subseteq \{(x,y) \in D_s \times
D_s \; | \;P_T(x)=P_T(y)\},$$
    \item is H-asymptotically $P_n$-weakly stationary if
$$ \{ (x_k,y_l)\; | \; k,l=1,...,s \} \subseteq \{(x,y) \in D_s \times
D_s \; | \; x=y=1 \mbox{ or }(x \neq 1 \mbox{ and } y \neq 1 \}.$$
    \item is $P_n$-weakly stationary  if and only if
$$ \{(x_k,y_l) \; | \; k,l=1,...,s \} \subseteq \{(x,x) \; | \; x \in D_s \}.$$
\end{itemize}

\end{Example}

\begin{Remark}
This characterization also yields examples of sequences being
T-cyclostationary with respect to $(P_n)_{\nin}$ but not
H-asymptotically $P_n$-weakly stationary: Any case with
$P_T(1)=P_T(x)$ for $x \neq 1$ could serve as an example (compare
Remark \ref{RemAsHarm}(3)) .

\end{Remark}

The next example is given in Leitner \cite{Leihyp} within the more
general framework of weakly harmonizability. It is an example of
$P_n$-harmonizability as long as the mentioned set $A$ is finite.

\begin{Example}\label{trunc}

Let $(X_n)_{\nin}$ be $P_n$-weakly stationary and $A$ a finite set
in $\Bbb N_0$. Then, the truncated process
\begin{eqnarray*}
\tilde{X}_n = \left\{ \begin{array}{lll}
 X_n  & \mbox{ if  } &  n \in A \\
0 &  else & \\
\end{array}\right.
\end{eqnarray*}
is $P_n$-harmonizable with covariance
$$ \tilde{d}(n,m) = \textbf{1}_{A \times A}(n,m) d(n,m).$$
Orthogonality of the polynomials with respect to the measure $\pi$
shows
\begin{eqnarray*}
\textbf{1}_{A \times A}(n,m) d(n,m)\! \! \! &=&\! \! \!
d(n,m)\!\int_{D_s}\!\!P_n(x)\sum_{s \in A}
        P_s(x) h(s) d\pi(x)\! \int_{D_s}\!\!P_m(y)\sum_{s \in A} P_s(y) h(s) d\pi(y) \\
&=&\int_{D_s \times D_s} \sum_{s,t \in
A}P_s(x)P_t(y)h(s)h(t)d(s,t) d \pi(x)d\pi(y).
\end{eqnarray*}

 Thus, one has the representation
$$\tilde{d}(n,m) =\int_{D_s \times D_s} P_n(x) P_m(y) f(x,y) d\pi(x)d\pi(y),$$
with the density function
$$f(x,y)= \sum_{s,t\in A} d(s,t) h(s) h(t) P_s(x) P_t(y).$$

\end{Example}

\subsection{$P_n$-weakly stationary sequences}
The best investigated class, connected to polynomial hypergroups
are $P_n$-weakly stationary sequences.

For a $P_n$-weakly stationary sequence $(X_n)_{\nin}$ with
spectral measure $\mu \in M^+(D_s)$ one has for all $n,m \in {\Bbb
N}_0$ (\ref{SpecK})
\begin{equation}
E(X_m \overline{X_n}) \; = \; d(m,n) \; = \; \int_{D_s} P_m(x) \;
P_n (x) \; d\mu (x).
\end{equation}
Therefore, $(X_n)_{\nin}$ belongs to the Karhunen class
((\cite{ChaRao}, p.59). We know from the harmonizable case that a
Cram\'{e}r representation theorem (see \cite{CraHilb}) holds.
Here, one gets an \textbf{orthogonal} stochastic measure $ Z :
{\cal B} \rightarrow L^2(P),\;\; {\cal B}$ being the Borel
$\sigma$-algebra on $D_s$ such that
\begin{equation}
X_k \; = \; \int_{D_s} P_k(x) \; dZ(x),
\end{equation}
where $\parallel Z(A)\parallel_2^2 \; = \mu (A) \;$ for all $A \in
{\cal B}.$

In fact, the following characterizations of weakly stationarity
are valid.

\begin{Theorem}\label{charHypostation}
For $(X_n)_{\nin} \subset L^2(\Omega,P)$ the conditions are
equivalent:
\begin{enumerate}
\item[(i)] $(X_n)_{\nin}$ is weakly stationary on the polynomial
hypergroup ${\Bbb N}_0$ induced by $(P_n)_{\nin}.$ \item[(ii)]
$E(X_m \overline{X_n})\;=\;\displaystyle\int_{D_s}
P_m(x)\;P_n(x)\;d\mu (x) \quad$ for every $n,m \in {\Bbb N}_0,\;$
where $\mu$ is a bounded positive Borel measure on $D_s.$
\item[(iii)] $X_k=\displaystyle\int_{D_s} P_k(x) \; dZ(x) \quad$
for every $k \in {\Bbb N}_0,\;$ where $Z$ is an orthogonal
stochastic measure on $D_s$.
\end{enumerate}
\end{Theorem}

\begin{proof}
$(i) \Rightarrow (ii)\;$ is already shown. $(ii) \Rightarrow (iii)
\;$ is based on the isometric isomorphism $\Phi$ between
$L^2(D_s,\mu)$ and $H = \overline {span \{ X_n: \nin \}} \subseteq
L^2(\Omega,P) \;$ determined by $\; \Phi (P_n) = X_n.\;$ The
stochastic measure is defined by $Z(A)=\Phi(\chi_A),\;\;A $ being
a Borel subset of $D_s.\;$ The construction parallels (like in the
harmonizable situation) the classical case (see Shiryayev
(\cite{Shir}, 395-403)).

Finally, assuming $(iii)$ we get for every $m,n \in {\Bbb N}_0$
$$ E(X_m \overline{X_n}) \;=\; \int_{D_s} P_m(x) \;P_n(x) \; d\mu(x),$$
where $\mu$ is defined by $\mu(A) \;= \parallel
Z(A)\parallel_2^2,\;\; A \in {\cal B}.\;$ The linearization of
$P_m(x)P_n(x)\;$ gives
$$ E(X_m \overline{X_n}) \;=\; \sum_{k=|n-m|}^{m+n} g(m,n,k)\; E(X_k\overline{X_0}).$$
\end{proof}

We look at some special classes of $P_n$-weakly stationary
processes, having their analogue in classical weakly stationary
processes \cite{lalest}:
\renewcommand{\labelenumi}{\arabic{enumi})}
\begin{enumerate}
\item White noise with respect to $(P_n)_{n \in \nz_0}$:

For all $n,m \in \nz_0$ it holds $EZ_m \overline Z_n =
\delta_{m,n} g(m,n,0)$.

\item Moving average-process MA($\infty$) or MA(p) for $p \in
\nz$:

 There is a sequence of complex numbers
$(a_n)_{n \in \nz_0} \in l^2 (h)$, a white noise $(Z_k)_{k \in
\nz_0}$ and associated shift operators $T_k$ (compare
\ref{shifts}) with:
\begin{equation}\label{MovAverage}
    X_n := \sum_{k=0}^\infty a_k T_n Z_k h(k), \quad n \in \nz_0.
\end{equation}
The spectral measure is with the Plancherel measure $\pi$ given
by:
\begin{equation}\label{MovAverageSpec}
    |\hat a |^2 \pi = | \sum_{k=0}^\infty a_k P_k h(k)|^2 \pi.
\end{equation}

\item Autoregressive processes AR(q):

There are $a_k \in \cz, \: k=1,...,q$, such that:
\[
 X_n + a_1 T_1 X_n + ... + a_q T_q X_n = Z_n
\]
\end{enumerate}



\section{Occurrence of $P_n$-weakly stationary sequences}
The concept of $P_n$-weakly stationary processes is motivated by
the idea that the observation of real events is often a result of
averaging "simpler" events which themselves may be modelled as
weakly stationary processes. For example, one may regard the
sequence of arithmetic means $ (1/(2n+1)\sum_{k=-n}^nY_k)_{n \in
\mathbb{N}_0}$ as an "integrated" observation of a non-observable
weakly stationary process $(Y_k)_{k\in \mathbb{Z}}$. The problem
is that even for the simplest rules of averaging the resulting
sequence is no longer weakly stationary. But certain classes of
such averaged processes are $P_n$-weakly stationary.

\subsection{Averaging classical weakly stationary processes}

$P_n$-weakly stationary processes occur from averaging of weakly
stationary processes as follows. Let $(P_n)_{n\in {\mathbb{N}_0}}$
be a sequence of polynomials which define a hypergroup as above.
Then these polynomials have for $x\in [-1,1]$ and $x=cos(t)$ a
trigonometric representation
\[
P_n(cos(t))=\sum_{k=-n}^n\alpha _{nk}e^{ikt}, \qquad t\in {]-\pi
,\pi ]}.
\]
If $(Y_n)_{n\in {\mathbb{Z}}}$ is a classical weakly stationary
sequence with mean $M$, covariance function $d_Y$, and spectral
measure $\mu _Y$, then
\[
X_n:=\sum_{k=-n}^n\alpha _{nk}Y_k
\]
defines an unbiased estimator of $M$. It can be shown (compare
\cite{Lasser}) that \noindent $(X_n)_{n\in {\mathbb{N}_0}}$ is
$P_n$-weakly stationary with mean $M$ and covariance function
\[
d(n)=\sum_{k=-n}^n\alpha _{nk}d_Y(k).
\]
The spectral measure is $\mu =T(\mu _Y)$, the image measure of
$\mu _Y$ with
respect to the map $T:\,]-\pi ,\pi ]\rightarrow [-1,1],\;\;T(t)=cos(t)$.\\

Important examples of estimators of the mean can be obtained in
the described way.

For example, the sequence of classical arithmetic mean estimators
is given by $\alpha_{n k} = 1/(2n+1)$ and is weakly stationary
with respect to the Jacobi polynomials $P_n^{(1/2,-1/2)}$.

More generally, the weighted mean estimators $M_n^{\lambda}$
connected to $P_n^{( \lambda-1/2,-1/2)}$ provide asympotically
optimal estimators of the mean for ARMA($p,q$) processes
\cite{laro}.

Another class of weighted mean estimators is given by
ultraspherical polynomials $P_n^{(\alpha, \alpha)}$. In \cite
{laro}, the coefficients $\alpha_{nk}$ have been calculated for
these polynomials which define hypergroups, if and only if $\alpha
\geq -1/2 $.\\

The lower two, of the plots below, show realizations of $P_n$-
weakly stationary sequences $X_n$, occurring from averaging:
\begin{center}
\includegraphics[width=12cm]{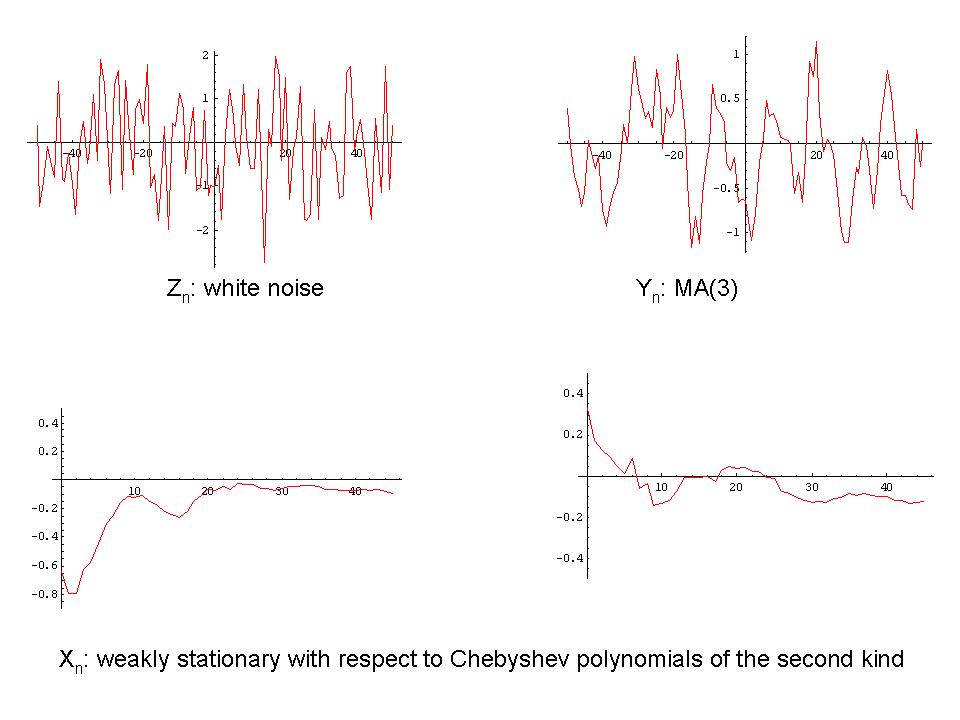}
\end{center}
Starting with a white noise sequence $(Z_k)_{k \in \Bbb Z}$ of
independent $N(0,1)$ distributed random variables, we define the
moving average ($MA(3)$) sequence:
$$Y_n :=1/4(Z_n+Z_{n-1}+Z_{n-2}+Z_{n-3}) \quad n \in \Bbb Z.$$
The Chebyshev polynomials of the second kind $P^{(1/2,1/2)}_n$
have for all $\nin$ the trigonometric representation (compare
\cite{Sz} (4.1.7)):
$$P^{(1/2,1/2)}_n (cos(t))=\frac{1}{n+1}\sum_{k=0}^{n} e^{i t (n-2k)},$$
and therefore, $$X_n := \frac{1}{n+1}\sum_{k=0}^{n} Y_{n-2k}$$
defines a sequence which is weakly stationary with respect to
$P^{(1/2,1/2)}_n$.

\subsection{Real and imaginary parts}

Let $(Y_n)_{n \in {\Bbb Z}}$ be a weakly stationary complex-valued
process with symmetry,  that is $Y_{-n} = \overline{Y_n}.\;$ The
random sequence $(U_n)_{\nin}$ of the real parts $U_n= $ Re $ Y_n
\; = \frac{1}{2} (Y_n+Y_{-n}),\;\; \nin \;$ is no longer weakly
stationary in the usual sense. However, one can easily check that
$$E(U_mU_n) \;=\; \frac{1}{2} \; E(U_{n+m}U_0) \;+\;\frac{1}{2} \; E(U_{|n-m|}U_0).$$
That means $(U_n)_{\nin}$ is a weakly stationary random sequence
on the polynomial hypergroup ${\Bbb N}_0$ induced by the Chebyshev
polynomials $T_n(x)$ of the first kind.

Denote the imaginary part of $Y_n$ by
$$ V_n \;=\; \mbox{Im } Y_n \;=\; \frac{1}{2i} \;(Y_n-Y_{-n}).$$
Since $V_0=0$ we have
$$ E(U_nU_0) \;+ \;iE(V_nU_0) \;=\; E(Y_n \overline{Y_0})\;=\;
\int_{-\pi}^\pi \cos (nt)\; d\mu (t) \;+\; i \int_{-\pi}^\pi \sin
(nt) \; d\mu (t),$$ where $\mu \in M^+ (\; ]-\pi,\pi]\,)\;$ is the
spectral measure of $(Y_n)_{n \in {\Bbb Z}}.$

The random sequence $(U_n)_{\nin}$ is weakly stationary with
respect to $T_n(x) = \cos (nt),\;$ where $x=\cos t \;$ for $x \in
[-1,1], \;\; t \in [0,\pi].\;$ Hence we have a unique spectral
representation
$$ E(U_nU_0)\;=\; \int_0^\pi \cos (nt) \; d \nu (t),$$
$\nu \in M^+([0,\pi]).\;$ Now it is clear that $\nu \;=
\mu|[0,\pi] \;+ \mu_1|\;]0,\pi[,\;$ where $\mu_1$ is the image
measure of $\mu$ under the mapping $t \rightarrow -t.\;$

Further we observe that
$$ E(U_mU_n)\;+\;E(V_mV_n)\;=\; \mbox{Re } E(Y_m\overline{Y_n}) $$
$$ =\;\int_{-\pi}^\pi \cos ((m-n)t)\; d\mu(t) \;=\; \int_0^\pi \cos ((m-n)t)\; d\nu (t) $$
$$ =\; \int_0^\pi \cos (mt) \;\cos (nt) \;
d\nu (t) \;+\; \int_0^\pi \sin (mt) \; \sin (nt) \; d\nu (t).$$
Since $E(U_mU_n)\;=\displaystyle\int_0^\pi \cos(mt)\; \cos(nt) \;
d\nu (t),\;$ it follows that
$$ E(V_mV_n) \;=\; \int_0^\pi \sin (mt) \; \sin (nt)\; d\nu (t). $$
Define for $\nin$:
$$ X_n:= \;\frac{1}{n+1} \; V_{n+1} \;=\; \frac{1}{n+1} \; \mbox{Im } Y_{n+1}.$$
Then
$$ E(X_mX_n) \;=\; \int_0^\pi \frac{\sin ((m+1)t)}{(m+1) \sin t} \;
\frac{\sin ((n+1)t)}{(n+1) \sin t} \; (\sin t)^2 \; d \nu (t).$$

By Theorem \ref{charHypostation} we see that $(X_n)_{\nin}$ is a
weakly stationary random sequence on the polynomial hypergroup
${\Bbb N}_0$ induced by the Chebyshev polynomials $U_n(x)$ of the
second kind. Notice that $U_n(x)\;= P_n^{(1/2,1/2)} (x) \;=
\displaystyle\frac{\sin ((n+1)t)}{(n+1) \sin t} ,\quad x= \cos t.$

\subsection{Stationary increments}

A random sequence $(Y_n)_{n \in {\Bbb Z}}$ is called a sequence
with stationary increments if all $E(Y_{n+k}-Y_k) $ depend only on
$n \in {\Bbb Z}$ and $E\left(
(Y_{n_1+k}-Y_k)\;(\overline{Y}_{n_2+k} - \overline{Y_k}) \right)
\;$ depend only on $n_1,n_2 \in {\Bbb Z},$ see Yaglom (\cite{Yag},
sec. 23).

It is readily seen, compare (\cite{Yag}, (4.227)) that
$$ Y_{n+k} -Y_k \; = \; \int_{-\pi}^\pi e^{ikt} \; \frac{e^{int} -1}{e^{it}-1} \; dZ(t),$$
where $Z$ is an orthogonal stochastic measure on $\;]-\pi,\pi].\;$
Hence we get for arbitrary $\nin$
$$ Y_{n+1} -Y_{-n} \;=\; \int_0^\pi \frac{ \sin \left( (2n+1)\frac{t}{2} \right)}
{\sin \left( \frac{t}{2} \right) }\; d \tilde{Z} (t),$$ where
$\tilde{Z}$ is the orthogonal stochastic measure on $[0,\pi]$
defined by $\tilde{Z} \;= Z|[0,\pi]\;+Z_1,\;\; Z_1(\;]a,b]) \;=
Z([-b,-a[\;) \;$ for $\;]a,b] \subseteq \;]0,\pi]\;$ and $Z_1
(\{0\} )=0.\;$

Putting
$$X_n := \; \frac{1}{2n+1} \; (Y_{n+1} -Y_{-n})$$
we have a weakly stationary sequence $(X_n)_{\nin}$ on the
polynomial hypergroup ${\Bbb N}_0$ induced by the Jacobi
polynomials $\;P_n^{(1/2,-1/2)} (x)=\; \displaystyle\frac{\sin
\left( (2n+1) \frac{t}{2} \right)}{(2n+1) \sin \left( \frac{t}{2}
\right)},\quad x=\cos t.$

\subsection{Random orthogonal expansions for density estimation}

Suppose that the distribution of a random variable $X$ is
absolutely continuous with respect to a positive Borel measure
$\pi$ on the interval $[-1,1],\;$ i.e. $P(X\in A) \;=
\displaystyle\int_A f(x) \;d\pi (x),$ where $f \in L^1(
[-1,1],\pi), \; f \geq 0.\;$ Consider the sequence $(p_n)_{\nin}$
of polynomials that are orthonormal with respect to $\pi$, and let
as before $P_n(x) = p_n(x)/p_n(1).\;$ Further assume that
$(P_n)_{\nin}$ induces a polynomial hypergroup on ${\Bbb N}_0.\;$
Given independent random variables $X_1,X_2,...,X_N$ equally
distributed as $X$, the unknown density function $f(x)$ can be
estimated by the random orthogonal expansion
$$ f_N(\omega;x):= \; \sum_{k=0}^{q(N)} a_{N,k} c_{N,k} (\omega) \; p_k(x).$$
$q(N)$ is the truncation point, and $a_{N,k}$ are numerical
coefficients to be chosen in an appropriate manner, see Lasser et
al. \cite{laob} and Devroye and Gy\"orfi \cite{Dev}. The random
coefficients are given by
$$ c_{N,k}(\omega) := \; \frac{1}{N} \; \sum_{j=1}^N p_k\left(X_j(\omega)\right).$$
Define $C_k:= c_{N,k} /p_k(1)\;$ for each $k \in {\Bbb N}_0.\;$
The random sequence $(C_k)_{k \in {\Bbb N}_0}$ is weakly
stationary on the polynomial hypergroup ${\Bbb N}_0$ induced by
$(P_k)_{k \in {\Bbb N}_o}.\;$ In fact we have $C_n =
\displaystyle\frac{1}{N} \sum\limits_{j=1}^N P_n(X_j)\;$ and hence
$$ E(C_mC_n) \;=\; \frac{1}{N^2} \; \sum_{i,j=1}^N E\left( P_m(X_i)\;P_n(X_j)\right) $$
$$ =\; \frac{1}{N^2} \; \sum_{i,j=1}^N \int_{-1}^1 P_m(x) \;P_n(x) \;f(x) \; d\pi(x) \;=
\; \int_{-1}^1 P_m(x) \;P_n(x) \; f(x) \; d\pi(x) $$
$$ = \; \sum_{k=|n-m|}^{n+m} g(m,n,k)\;\int_{-1}^1 P_k(x) \; P_0(x) \; f(x) \; d\pi(x) \;=
\; \sum_{k=|n-m|}^{n+m} g(n,m,k) \; E(C_kC_0).$$ Estimating the
density function $f$ by random orthogonal expansion is hence
strongly related to estimating the spectral measure $f \pi$ of the
random sequence $(C_k)_{k \in {\Bbb N}_0}, \;$ where $(C_k)_{k \in
{\Bbb N}_0}$ is weakly stationary with respect to $P_k \; =
p_k/p_k(1).$

\subsection{Stationary radial stochastic processes on homogeneous trees}

We denote by $T$ a homogeneous tree of degree $q \geq 1$ with
metric $d$. Let $G$ be the isometry group of $T,\;\; t_0 \in T$ an
arbitrary but fixed knot of $T$ and let $H$ be the stabilizer of
$t_0$ in $G.\;$ We identify $T$ with the coset space $G/H$ and
call a mapping on $T$ radial if it depends only on $|t|\;=
d(t,t_0).\;$ A square integrable stochastic process $(X_t)_{t \in
T}$ is called stationary, if there exists a function $\phi:{\Bbb
N}_0 \rightarrow {\Bbb R} \;$ such that
$$ E(X_s \overline{X_t})\;=\; \phi(d(s,t))$$
for all $s$ and $t$ in $T$, compare Arnaud \cite{Arn}. It is
known, see also \cite{Arn}, that the following spectral
representation is true
$$ E(X_s \overline{X_t}) \; = \; \int_{-1}^1 P_{d(s,t)} (x) \; d \mu (x),$$
where $(P_n (x))_{\nin}$ are the orthogonal polynomials (also
called Cartier-Dunau polynomials) corresponding to homogeneous
trees of degree $q$. The polynomials $P_n(x)$ are determined by
$$ P_0(x) \; =1, \qquad P_1(x) \; =x \qquad \qquad \mbox{and}$$
$$ P_1(x) \;P_n(x) \; = \; \frac{q}{q+1} \; P_{n+1} (x) \; + \;
\frac{1}{q+1} \; P_{n-1} (x),\qquad \nin.$$ They induce a
polynomial hypergroup on ${\Bbb N}_0$, see Lasser \cite{la01}.

We assume moreover that the stationary stochastic process
$(X_t)_{t \in T}$ is radial, that is $X_t = X_s$ if $ |t| =
|s|.\;$ Putting $X_n:= X_t$ whenever $|t|=n \;$ we get a
well-defined random sequence. We have
$$ \int_H P_{d(h(s),t)} (x) \; d \beta (h) \; = \; P_{|s|} (x) \; P_{|t|} (x),$$
where $\beta$ is the Haar measure on the compact stabiliser $H$,
and $h(s)$ is the action of $h \in H$ on $T=G/H,\;$ compare e.g.
Cowling et al \cite{Cow}. From the above spectral representation
we can derive
$$E(X_m\overline{X_n}) \; = \; \int_{-1}^1 P_m(x) \; P_n(x) \; d \mu (x).$$
Again, by Theorem \ref{charHypostation} we obtain that
$(X_n)_{\nin}$ is weakly stationary on the polynomial hypergroup
${\Bbb N}_0$ induced by the Cartier-Dunau polynomials
$(P_n)_{\nin}.$

\chapter{Finding Structure from Realizations}
\section{Covariance estimation}

In this chapter we present mean square consistent estimators for
all $P_n$-weakly stationary processes under reasonable conditions.

Processes which are weakly stationary with respect to the
Chebyshev polynomials of the first kind $T_n(x)$ will serve as a
germ of our estimation theory. With a covariance satisfying the
equation
\begin{equation}
d(n,m)=\frac 12(d(|n-m|)+d(n+m)),  \label{coveq}
\end{equation}
these $T_n$-weakly stationary processes are for large $n$, $m$
``close" to classical weakly stationary processes, if $d(k)$ tends
to zero with growing $k$. Yet, one has to keep in mind that the
associated dual space is not the torus but the real interval
$[-1,1]$. We will first construct estimators for the special case
of $T_n$-weakly stationary processes and then derive therefrom
solutions for the general case by linear transformations.

\subsection{Estimators for $T_n$-weakly stationary processes }

The estimators $\tilde{d}_N(s)$ and $\stackrel{\approx
\qquad}{d_N(s)}$ introduced in this section exploit the indicated
similarity of $T_n$-weakly stationary processes and weakly
stationary processes. The estimator
\begin{equation} \label{schatz1}
\tilde{d}_N(s):=\frac 2{N-s+1}\sum_{k=0}^{N-s}X_k\overline{X}_{k+s},
 \qquad 0 \leq s \leq N
\end{equation}
equals, except for the factor 2, the classical estimator, but is
now no longer unbiased. Its definition and (\ref{coveq}) show for
the mean of $\tilde{d}_N(s)$:
\begin{eqnarray}\label{bias}
E\tilde{d}_N(s) & = & \frac 1{N-s+1} \sum_{k=0}^{N-s} (d(s) + d(2k+s)) \\
                & = & d(s) + \frac 1{N-s+1} \sum_{k=0}^{N-s} d(2k+s).  \nonumber
\end{eqnarray}
However, not too restrictive conditions on the spectral measure guarantee
asymptotic unbiasedness, that is
\[
\lim_{N\rightarrow \infty }E\tilde{d}_N(s)=d(s).
\]
\vspace{0.1 em}
\begin{Theorem} \label{s1}

Let the spectral measure $\mu $ of a centered $T_n$-weakly stationary \\
process $(X_n)_{n\in {\mathbb{N}_0}}$ be of the form $\mu =f\pi +\mu _d$ with $f\in
L^1(\pi )$, and a discrete measure $\mu _d$ having no atoms or limit points
of atoms in $\{-1,+1\}$. Then $\tilde{d}_N(s)$ is for fixed $s\in \mathbb{N}_0$
an asymptotically unbiased estimator of $d(s)$.
\end{Theorem}

\begin{Proof}
The spectral representation of $d(s)$ and the symmetry
$ T_n(x)=(-1)^nT_n(-x)$ yield (note that $D_s=[-1,1]$):
\begin{eqnarray}
E\tilde{d}_N(s)-d(s) &=&\frac 1{N-s+1}\int_{-1}^1\frac
12\sum_{k=s}^{2N-s}(T_k(x)+(-1)^sT_k(-x))f(x)d\pi (x)
 \nonumber \\
& &+ \ \frac 1{N-s+1}\int_{-1}^1\frac
12\sum_{k=s}^{2N-s}(T_k(x)+(-1)^sT_k(-x))d\mu _d(x).\label{bias1}
\end{eqnarray}
According to the Riemann-Lebesgue lemma $\int_{D_s}T_n(x)f(x)d\pi
(x)\rightarrow 0$ for \mbox{$n \rightarrow \infty$} holds, and
therefore the first summand on the right side of (\ref{bias1})
also vanishes asymptotically. The Christoffel Darboux identity for
Chebyshev polynomials with $y=1$ and $x\neq 1$
(\ref{Christoffel1}):
\begin{equation}
T_0(x)+2\sum_{k=1}^NT_k(x)=\frac{T_{N+1}(x)-T_N(x)}{x-1}
\end{equation}
and $|T_n(x)|\leq 1$ on $D_s$ imply for $s\in \mathbb{N}_0$:
\begin{equation}
(\quad \int_{-1}^1 \frac 12 | \sum_{k=s}^{2N-s} ( T_k(x) +
(-1)^sT_k(-x) )|
 d\mu _d(x) \leq C + \int_{-1}^1 \frac 1{1-x^2} d\mu _d(x) \qquad
\nonumber \label{bias2}
\end{equation}
with a constant $C$ depending only on $s$. The existence of the
last integral is guaranteed by the assumptions on $\mu _d$. The
second summand in (\ref{bias1}) is thus $O(1/N)$, and the bias of
$\tilde{d}_N(s)$ tends to zero with growing $N$.
\end{Proof}

If $d$ has the property $\displaystyle{|\sum_{k=0}^\infty
d(2k+s)|<\infty }$, then of course, the bias of $\tilde{d}_N(s)$
is $O(1/N)$. But this condition is not always satisfied, even if
$\mu =f\pi $ with continuous $f$. On the contrary, the next
theorem demonstrates that such an asymptotic decrease is valid for
the bias of the Fej\'{e}r weighted estimator
\begin{equation} \label{schatz2}
\qquad \stackrel{\approx \qquad}{d_N(s)}:= \frac
4{N-s+1}\sum_{k=0}^{N-s}(1-\frac k{N-s+1})X_k \overline{X}_{k+s},
\quad 0 \leq s \leq N \qquad \nonumber
\end{equation}
under rather general assumptions on the Radon-Nikodym derivative $f$ of $\mu$.
As usual, we call \mbox{$t_0 \in ]-\pi ,\pi ]$} a Lebesgue point of a function
\mbox{$g \in L^1(]-\pi ,\pi ],dt)$}, $dt$ being the Lebesgue measure, if there exists
a constant $c$ with
\[
\lim_{h\rightarrow 0}\frac 1{2h}\int_{-h}^h|g(t_0+t)-c|dt=0,
\]
and define $g(t_0):=c$.

\begin{Theorem}
Let the spectral measure $\mu $ of a centered $T_n$-weakly
stationary process $(X_n)_{n\in {\mathbb{N}_0}}$ be as in Theorem
\ref{s1}, and let $0$ and $ \pi $ be Lebesgue points of
$f(cos(.))$. Then $\stackrel{\approx \qquad}{d_N(s)}$ is for fixed
$s\in \mathbb{N}_0$ an asymptotically unbiased estimator of $d(s)$
with
\[
E\stackrel{\approx \qquad}{d_N(s)}-d(s)=O(1/N).
\]
\end{Theorem}

\begin{Proof}
We calculate the bias of $\stackrel{\approx
\qquad}{d_N(s)}$ to be
\begin{equation}\label{bias2Tilde}
\quad E\stackrel{\approx \qquad}{d_N(s)} - d(s) = \frac 1{N-s+1}
(d(s)+ 2\sum_{k=0}^{N-s}(1- \frac k{N-s+1})d(2k+s)).\qquad
\nonumber \label{bs1}
\end{equation}
Denoting the classical Fourier transform on $]-\pi ,\pi ]$ by
$^{\wedge }$, and defining \\\mbox{$f \circ \cos (.) := f(\cos(.))
$}, we find with an elementary computation that
\[ \int_{-1}^1T_k(x)f(x)d\pi(x) = (f \circ \cos)^{\wedge}(k) \]
and obtain
\begin{eqnarray}
\sum_{k=0}^{N-s}(1-\frac k{N-s+1})d(2k+s) &=&\sum_{k=0}^{N-s}(1-\frac
k{N-s+1})(f\circ \cos)^{\wedge }(2k+s)  \nonumber  \label{fe1} \\
& +& \sum_{k=0}^{N-s}\!(1\!-\!\frac
k{N-s+1})\!\int_{-1}^1\!T_{2k+s}(x)d\mu_d(x).
\end{eqnarray}
As in the proof of Theorem \ref{s1}, we conclude the boundedness
of the last sum as a consequence of the choice of $\mu _d$. The
first sum on the right side of (\ref{fe1}) is now considered
further. The symmetry of $f\circ cos$ yields the identity
\[
(f\circ cos(\pi- .))^{\wedge }(l)=(-1)^l(f\circ cos)^{\wedge }(l),
\]
and consequently we calculate the Fourier coefficients of
\[
g_s(t):=e^{ist}(f\circ cos(t)+(-1)^sf\circ cos(\pi -t))
\]
as $g_s^{\wedge }(l)=(1+(-1)^l)f\circ cos(l+s)$. We get the equation
\begin{equation}
\; \; \sum_{k=0}^{N-s}(1 \! - \! \frac k{N-s+1}) (f \circ cos
)^{\wedge } (2k \! + \! s ) = \frac12 \sum_{l=0}^{2N-2s+1} \! \!
(1 \! - \! \frac l{2(N-s+1)}) g_s^{\wedge }(l)  \nonumber
\label{devil}
\end{equation}

by utilizing $g_s^{\wedge }(2l+1)=0$ and by substituting $l=2k$.
The inequality
\begin{eqnarray}\label{onTorus}
|g_s(t)-g_s(0)| &\leq &|f\circ cos(t)-f(1)|+|f\circ cos(\pi -t)-f(-1)|
\nonumber \\
& &+ \ |e^{ist}-1|\;\;|f(1)+(-1)^sf(-1)|  \nonumber
\end{eqnarray}
and the assumption that $f\circ cos$ has Lebesgue points at $t=0$ and $t=\pi
$ clearly imply that $g_s$ has a Lebesgue point at $t=0$:
\[
\lim_{h\rightarrow 0}\frac 1{2h}\int_{-h}^h|g_s(t)-g_s(0)|dt=0.
\]
The required boundedness in (\ref{fe1}) to prove the theorem is
now guaranteed by the convergence of the Fej\'{e}r series
(\ref{devil}) in Lebesgue points.

\end{Proof}

\subsection{The least square estimator{\normalsize {\rm
\ $d_N^{LS}$}}}
For classical weakly stationary processes $(Y_k)_{k\in Z}$ with
$EY_k=0$, \mbox{$ E|Y_k|^4<\infty $}, and covariance function
$R(s)$, the variables $Y_k \overline{Y}_{k+s}$,
\mbox{$s=0,...,N$}, \mbox{$k=0,...,N-s$} define a linear model
with the equations $Y_k\overline{Y}_{k+s}=R(s)$, for example
And\v{e}l (\cite{andl}, Ch.12). The least square estimator of this
model is just the classical unweighted \mbox{estimator} of $R(s)$:
\begin{equation}
\hat{R}_N(s)=\frac 1{N-s+1}\sum_{k=0}^{N-s}Y_k\overline{Y}_{k+s}.  \label{kl}
\end{equation}
To get an unbiased estimator for $P_n$-weakly stationary processes in the
same way, one has to choose the equations of the linear model. We define the
estimator $d_N^{LS}(s)$ to be the least square solution based on the
realizations of $X_k\overline{X_l}$ with $(k,l)\in I_N$, where
\[
I_N:=\{ (k,l)| \ k,l \in \mathbb{N}_0, \ k+l \leq N, \ 0 \leq k \leq l \leq N \}.
\]
\newpage

\begin{Proposition}\label{EstimationProp}
Let $(X_n)_{n\in \mathbb{N}_0}$ be a $P_n$-weakly stationary process with
\mbox{$EX_n=0$} and $E|X_n|^4<\infty $ for all $n\in \mathbb{N}_0$. Then:\\

a) The least square estimator $d_N^{LS}(s)$ associated to the equations
\[
EX_k\overline{X}_l=\sum_{s=|k-l|}^{k+l}g(k,l,s)d(s),\qquad (k,l)\in I_N
\]
is uniquely determined.

 b) With a lexicographical ordering on
$I_N$ and the notations
\[
\begin{array}{ll}
\begin{array}{lcl}
{\bf A}_N & := & (g(k,l,s))_{((k,l),s)\in I_N\times \{0,\ldots ,N\}}, \\
[0.2cm]
{\bf {\hat{d}}}_N & := & (d_N^{LS}(0),\ldots ,d_N^{LS}(N))^T, \\[0.2cm]
&  &
\end{array}
&
\begin{array}{lcl}
{\bf X}_N & := & (X_k\overline{X}_l)_{(k,l)\in I_N}, \\[0.2cm]
{\bf G}_N & := & Cov({\bf X}_N,{\bf X}_N) \\[0.2cm]
&  &
\end{array}
\end{array}
\]
the following equations are valid:

$$\hat{{\bf d}}_N=({\bf A}_N^T{\bf A}_N)^{-1}{\bf A}_N^T{\bf
X}_N$$
$$Cov(\hat{{\bf d}}_N,\hat{{\bf d}}_N)=({\bf A}_N^T{\bf A}_N)^{-1}{\bf A}_N^T
{\bf G}_N{\bf A}_N({\bf A}_N^T{\bf A}_N)^{-1}.$$
\end{Proposition}

\begin{Proof}
$(g(0,l,s))_{0\leq l,s \leq N}$ is a submatrix of ${\bf A}_N$. Due
to $\displaystyle g(0,l,s)=\delta _{ls}$ the matrix ${\bf A}_N$
has the full rank $N+1$, and all statements follow as general
properties of linear models (for example \cite{andl}, $12.2$).
\end{Proof}

\newpage
The above representation of $d_N^{LS}(s)$ and its covariance is
not well suited for further investigations. To perform consistency
analysis we focus again on the case of $T_n$-weakly stationary
processes where we can give an explicit representation of
$d_N^{LS}(s)$. For clarity, we will denote this estimator as
$d_{TN}^{LS}(s)$.

\begin{Theorem}\label{EstimationMain}
Let $(X_n)_{n\in \mathbb{N}_0}$ be a $T_n$-weakly stationary process with $EX_n=0$
and $E|X_n|^4<\infty $ for all $n\in \mathbb{N}_0$. Then the least square
estimator $d_{TN}^{LS}(s)$ has for {\bf odd N} the representation
\begin{equation}
\quad d_{TN}^{LS}(s) \! = \! \frac 8{(\! N \! + \! 4 )^2 \! - \!
1} ((\! N \! + \! 2 ) X_0\overline{X}_s + \! \frac{N \! + \! 1}
{2} \! \! \! \! \sum_{(k,l)\in I_1^N(s)} \! \! X_k\overline{X}_l-
\! \! \! \! \sum_{(k,l)\in I_2^N(s)} \! \! \! X_k \overline{X}_l)
\nonumber \label{schatz3}
\end{equation}
with the sums running over:\\
\noindent $I_1^N(s):= I_N \cap \{(k,l)| \ 0 < k,\ |k \pm l| =
s\},$
\\$I_2^N(s):= I_N \cap
\{(k,l)| \ k+l \neq s,\ l-k \neq s, \ k+l+s: \mbox{ even}\}$.\\
\begin{flushleft}
The estimator for {\bf even N} is given by:

\end{flushleft}
${d_{TN}^{LS}(s)}=\left\{
\begin{array}{cl}
d_{TN+1}^{LS}(s) & \mbox{ for s: even} \\[0.1cm]
d_{TN-1}^{LS}(s) & \mbox{ for s: odd .}
\end{array}
\right. $
\end{Theorem}

The proof of this theorem is rather technical, and in order not to
disturb the main course we postpone it to the end of the chapter.

For illustration we show the domains $I_N$, $I_1^N(s)$ and
$I_2^N(s)$  for the both estimators:

\newpage

\begin{eqnarray*}
d_{T10}^{LS}(5) &=&d_{T9}^{LS}(5)=\frac 1{21}(11 X_0\overline{X}
_5+5\sum_{(k,l)\in I_1^9(5)}X_k\overline{X}_l-\sum_{(k,l)\in I_2^9(5)}X_k
\overline{X}_l), \\
d_{T10}^{LS}(6) &=&d_{T11}^{LS}(6)=\frac 1{28}(13 X_0\overline{X}
_6+6\sum_{(k,l)\in I_1^{11}(6)}X_k\overline{X}_l-\sum_{(k,l)\in
I_2^{11}(6)}X_k\overline{X}_l).
\end{eqnarray*}

\begin{center}
\includegraphics[width=8cm,keepaspectratio]{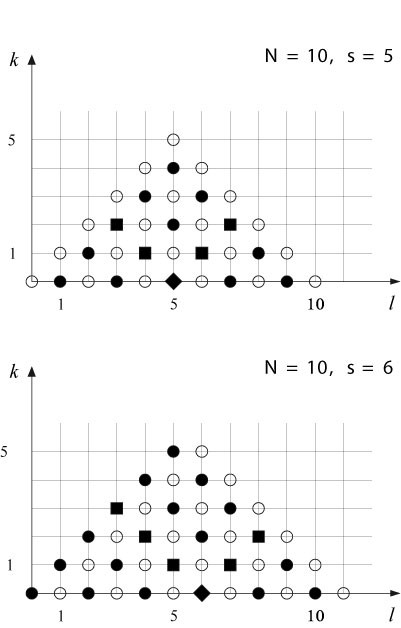}
\end{center}

FIG.: The upper and lower diagram show the cases $N=10$, $s=5$
resp. $N=10$, $s=6$. The sets \mbox{$I_1^{10}(5)=I_1^9(5)$} resp.
\mbox{$I_1^{11}(6)=I_1^{10}(6)$} are marked by $\rule{2mm}{2mm}$
and \mbox{$I_2^{10}(5)=I_2^9(5)$} resp.
\mbox{$I_2^{11}(6)=I_2^{10}(6)$} are marked by $\bullet $. All
marked lattice points together, including $\circ $ and $\triangle
$, represent $I_{10}$ resp. $I_{11}$.

\vspace{2cm}

\subsection{Mean square consistency}

In this section we  prove the consistency of $\tilde{d}_N$,
$\stackrel{\approx}{d}_N$ and $d_{TN}^{LS}$ for real Gaussian
$T_n$-weakly stationary processes having spectral measures without
atoms and show how these results can be used to define mean square
consistent estimators for general $P_n$-weakly stationary
processes.

An application of the Wiener type Corollary \ref{WienerCheby} to
(\ref{bias}) and (\ref{bias2Tilde}) yields:

\begin{Corollary} \label{c1}
 Let the spectral measure $\mu $ of a centered $T_n$-weakly
stationary process have no atoms. Then the estimators
$\tilde{d}_N$ and $\stackrel{\approx}{d}_N$ are asymptotically
unbiased.
\end{Corollary}

Remark: Note, that the corollary is not a consequence of Theorem
\ref{s1}, since $\mu $ is more general here.

\begin{Theorem} \label{t1}
Let $(X_n)_{n\in {\mathbb{N}_0}}$ be a centered real Gaussian
$T_n$-weakly stationary process with a spectral measure $\mu $ which has no
atoms. Then each estimator of $d(s)$ of the form
\[
\hat{d}_N(s)=\sum_{0\leq k,l\leq N}a_{kl}^{Ns}X_kX_l
\]
which is asymptotically unbiased and has coefficients
$a_{kl}^{Ns}\in {\mathbb{R}} $ which satisfy for all but finite
$N$ the inequalities
\begin{equation}\label{crazy}
\sum_{k=0}^N(|a_{kl}^{Ns}|+|a_{lk}^{Ns}|)\leq \frac {C(s)}{N+1},
\qquad 0 \leq l \leq N,
\end{equation}
where $C(s)$ depends only on $s$,

is mean square consistent:
\[
\lim_{N\rightarrow \infty }E(\hat{d}_N(s)-d(s))^2=0.
\]
\end{Theorem}

\begin{Proof}
The process satisfies the Isserlis identity \cite{Isserlis}
\[
E(X_mX_nX_rX_t)=d(m,n)d(r,t)+d(m,r)d(n,t)+d(m,t)d(n,r),
\]
and therefore, elementary computations show for the covariances:
\begin{equation}
Cov(X_mX_n,X_rX_t)=d(m,r)d(n,t)+d(m,t)d(n,r).  \label{zz}
\end{equation}
This allows us to calculate upper bounds for the variance of $\hat{d}_N(s)$.
The \linebreak assumptions on $a_{kl}^{Ns}$ and the Cauchy Schwarz
inequality now yield:
\begin{eqnarray*}
\lefteqn{E(\hat{d}_N(s)-E\hat{d}_N(s))^2 = \sum_{0\leq m,n,r,t\leq
N}a_{mn}^{Ns}a_{rt}^{Ns}Cov(X_mX_n,X_rX_t)} \\
&\leq &\sum_{0\leq m,n,r,t\leq
N}|a_{mn}^{Ns}|(|a_{rt}^{Ns}|+|a_{tr}^{Ns}|)|d(m,r)d(n,t)| \\
&\leq &(\sum_{0\leq m,r\leq
N}(\sum_{n=0}^N|a_{mn}^{Ns}|
\sum_{t=0}^N(|a_{rt}^{Ns}|+|a_{tr}^{Ns}|))d(m,r)^2)^{1/2} \\
&&\mbox{ $\cdot$ }(\sum_{0\leq n,t\leq
N}(\sum_{m=0}^N|a_{mn}^{Ns}|
\sum_{r=0}^N(|a_{rt}^{Ns}|+|a_{tr}^{Ns}|))d(n,t)^2)^{1/2} \\
&\leq &\frac{C(s)^2}{(N+1)^2}\sum_{k,l=0}^Nd(k,l)^2.
\end{eqnarray*}

The consistency of $\hat{d}_N(s)$ now follows from the assumption,
Corollary \ref{WienerCheby} and the decomposition
\begin{equation}
E(\hat{d}_N(s)-d(s))^2=(E\hat{d}_N(s)-d(s))^2+E(\hat{d}_N(s)-E\hat{d}
_N(s))^2. \quad \Box  \qquad \nonumber \label{zerl}
\end{equation}
\end{Proof}

\begin{Remark}
Note that the summation variable in (\ref{crazy}) is
$k$ and that $l$ and $s$ are fixed. The boundedness condition is
especially satisfied for the estimators $\tilde{d}_N(s)$,
$\stackrel{\approx \qquad}{d_N(s)}$ and $d_{TN}^{LS}(s)$ which is
used in the following consistency theorem.
\end{Remark}

\begin{Theorem}\label{t2}
Let $(X_n)_{n \in {\mathbb{N}_0}}$ be a centered real Gaussian $T_n$-weakly
stationary process with a spectral measure $\mu $ which has no atoms. Then
the estimators $\tilde{d}_N(s)$, $\stackrel{\approx \qquad}{d_N(s)}$ and $
d_{TN}^{LS}(s)$ are mean square consistent for all $s\in \mathbb{N}_0$.
\end{Theorem}

\begin{Proof}
The explicit representations (\ref{schatz1}) and (\ref{schatz2})
of the estimators $\tilde{d}_N(s)$ and $\stackrel{\approx
\qquad}{d_N(s)}$ show that one gets $(|a_{kl}^{Ns}|+|a_{lk}^{Ns}|)
\neq 0 $ only for $k=l-s$ or $k=l+s$ and that $N \geq s$ implies
\[\sum_{k=0}^N(|a_{kl}^{Ns}|+|a_{lk}^{Ns}|) \leq \frac {8(s+1)}{N+1} \]
for both estimators. In the case of $d_{TN}^{LS}(s)$ it is clear
from (\ref{crazy}) that for fixed $l$ at most two of the
$a_{kl}^{Ns}$ are of order $O(1/N)$. The remaining $O(N)$
coefficients which are not equal to zero are related to $I_2^N(s)$
and of order $O(1/N^2)$. These observations and Corollary \ref{c1}
assure that all three estimators satisfy the assumptions of
Theorem \ref{t1}.
\end{Proof}

\begin{Remark}
In the situation of Theorem \ref{t2} one can easily
calculate the variances of $\tilde{d}_N(0)$ and $\stackrel{\approx
\qquad}{d_N(0)}$ using again (\ref{zz}).
\begin{eqnarray*}
Var(\tilde{d}_N(0)) &=&\frac 8{(N+1)^2}\sum_{k,l=0}^Nd(k,l)^2,
\end{eqnarray*}
\\[-0.6cm]
\begin{eqnarray*}
Var(\stackrel{\approx \qquad}{d_N(0)}) &=&\frac{32}{(N+1)^2}\sum_{k,l=0}^N(1-\frac
k{N+1})(1-\frac l{N+1})d(k,l)^2 \\
&\geq &\frac 8{(N+1)^2}\sum_{k,l=0}^{N/2}d(k,l)^2.
\end{eqnarray*}
Now, again Corollary \ref{WienerCheby} shows that for
$\tilde{d}_N(0)$ and $\stackrel{\approx \qquad}{d_N(0)}$ to be
mean square consistent it is necessary and sufficient that $ \mu $
has no atoms. We do not know whether this result, which likewise
holds for the classical estimator $\hat{R}_N(0)$, is also valid
for $ d_{TN}^{LS}(0)$.
\end{Remark}

From explicit estimators for the $T_n$-weakly stationary
case we will now construct estimators for the case of general
$P_n$-weakly stationary processes. Our tool is the linear
transformation defined by the connection coefficients $c_{TP}(n,k)$ of
the representation
\begin{equation}\label{Connect}
    T_n(x)=\sum_{k=0}^nc_{TP}(n,k)P_k(x).
\end{equation}
The reverse representation is
\begin{equation}\label{ConnectRev}
    P_n(x)=\sum_{k=0}^nc_{PT}(n,k)T_k(x).
\end{equation}
\begin{Proposition}

 Let $(X_n)_{n\in \mathbb{N}_0}$ be a
centered $P_n$-weakly stationary process with spectral measure
$\mu $ satisfying $supp(\mu )\subseteq [-1,1]$. Then
\[
X_n^T:=\sum_{k=0}^nc_{TP}(n,k)X_k
\]
defines a centered $T_n$-weakly stationary process with the same spectral
measure.
\end{Proposition}

\begin{Proof}
The normalization $P_k(1)=T_k(1)=1$ for $k\in {\mathbb{N}}_0$
causes \mbox{$\sum_{k=0}^nc_{TP}(n,k)=1$} and thus $(X_n^T)_{n\in
\mathbb{N}_0}$ is centered. The isometric isomorphism between
$L^2(D_s,\mu )$ and $ \overline {span \{ X_n: \nin \}}$ mentioned
in the proof of Theorem \ref{charHypostation} yields for the
covariance of $(X_n^T)_{n\in \mathbb{N}_0}$:
\[
d^T(n,m)=E(X_n^T\overline{X_m^T})=\int_{D_s}T_n(x)T_m(x)d\mu (x),
\]
and linearization of $T_n(x)T_m(x)$ in the integral shows the equation
\[
d^T(n,m)=\frac 12(d^T(|n-m|)+d^T(n+m)).
\]
The assumption $supp(\mu )\subseteq [-1,1]$ allows a spectral representation
\begin{equation}
d^T(n)=\int_{-1}^1T_n(x)d\mu _T(x)  \label{dTn}
\end{equation}
with spectral measure $\mu _T=\mu $. The required boundedness of $d^T$ is
now obvious:
\[
|d^T(n)|\leq \mu (D_s).
\]
\end{Proof}

\begin{Remark}
The assumption $supp(\mu )\subseteq D_s\subseteq [-1,1]$ does not
impose any restrictions on the class of $P_n$-weakly stationary
processes. With $a_0>0$ and $a_0P_1(x)=x-b_0$, the inclusion
$D_s\subseteq [1-a_0,1]$ holds in general \cite{la01}, and
therefore one can choose a linear transformation $x^{\prime
}=ax+(1-a)$ which maps $D_s$ into $[-1,1]$ without changing the
induced hypergroup or covariance structure; the transformed
polynomials have the same linearization coefficients.
\end{Remark}

Each consistent estimator of $d^T$ can now be used to define a
consistent estimator of $d$. We thereto use the connection
coefficients $c_{PT}$ occuring in the representation of $P_n(x)$
with regard to $T_k(x)$, $ k=0,...,n$.

\begin{Theorem} \label{t3}
Let $(X_n)_{n\in \mathbb{N}_0}$ be a centered
$P_n$-weakly stationary process with $ supp(\mu )\subseteq
[-1,1]$, and $\hat{d}_N^T$ a mean square consistent estimator of
the covariance function of the associated process $ (X_n^T)_{n\in
\mathbb{N}_0}$. Then
\[
[\hat{d}^T]_N(s):=\sum_{k=0}^sc_{PT}(s,k)\hat{d}_N^T(k)
\]
defines a mean square consistent estimator of $d(s)$.
\end{Theorem}

\begin{Proof}
First we observe that (see (\ref{dTn})):
\begin{eqnarray*}
d(s) = \int_{D_s}P_s(x)d\mu(x)
&=&\sum_{k=0}^sc_{PT}(s,k)\int_{-1}^1T_k(x)d\mu _T(x) \\
&=&\sum_{k=0}^sc_{PT}(s,k)d^T(k).
\end{eqnarray*}
Then we perform the decomposition (compare (\ref{zerl})):
\begin{eqnarray*}
E([\hat{d}^T]_N(s)-d(s))^2 = (\sum_{k=0}^sc_{PT}(s,k)(E\hat{d}_N^T(k)-d^T(k)))^2 \\
+ \sum_{0\leq k,l\leq s}c_{PT}(s,k)c_{PT}(s,l)
             Cov(\hat{d}_N^T(k),\hat{d}_N^T(l)). \quad \label{err}
\end{eqnarray*}
\end{Proof}

The consistency of $\hat{d}_N^T(k)$, \mbox{$k=0,...,s$} now
implies that all covariances and bias terms in (\ref{err})
disappear asymptotically. $\Box $

\begin{Corollary}
Let $(X_n)_{n\in \mathbb{N}_0}$ be a $P_n$-weakly stationary real
Gaussian process with spectral measure $\mu $ having no atoms and
$supp(\mu )\subseteq [-1,1]$. Then the estimators
$[\tilde{d}]_N(s)$, $[\stackrel{\approx }{d}]_N(s)$ and
$[d_T^{LS}]_N(s)$ are mean square consistent for all $s\in
\mathbb{N}_0$.
\end{Corollary}

\begin{Proof}
This follows at once from Theorem \ref{t2} and Theorem \ref{t3}.
One has only to observe that $(X_n^T)_{n\in \mathbb{N}_0}$ is
again Gaussian and has the same spectral measure as $(X_n)_{n\in
\mathbb{N}_0}$.
\end{Proof}

\begin{Remark}
Estimation of the $d(0),...,d(M)$ requires the coefficients
$c_{TP}(l,k)$, $c_{PT}(l,k)$, $l=0,...,M$, $k=0,...,l$. In some
cases these are explicitly known, for example Askey (\cite{AsII},
Lect.7). In general, connection coefficients can be calculated
from the recursion coefficients of the involved polynomial systems
(compare \cite{Ask}). In any case, the advantage of using
$[\hat{d}]_N(s)$ with any explicit estimator $\hat{d}$ for $d^T$
instead of applying $d_N^{LS}(s)$ is obvious: the calculation of
$d_N^{LS}(s)$ needs a matrix inversion for each $N$.
\end{Remark}

\subsection{Coefficients of associated Chebyshev processes}

Covariance estimations involving associated Chebyshev processes
require the connection coefficients $c_{T P}$ as well as $c_{P
T}$. The prediction algorithm, presented in the last chapter,
calculates the coefficients $c_{\phi P}$ for monic polynomials
$\phi$. Here, we consider polynomials $P_n$ and $Q_n$, both having
normalization $P_n(1)=Q_n(1)=1$. Specializing to the polynomials
of the considered hypergroup  and to Chebyshev polynomials of
first kind, respectively, yields the coefficients $c_{T P}$ and
$c_{T P}$.

Now, let $ (P_n)_{n \in \nz_0}$ und $(Q_n)_{n \in \nz_0} $ be
orthogonal polynomial systems in hypergroup normalization
satisfying the three term recursions

\[
P_0(x)=1; \; P_1(x)= \frac{x-b_0}{a_0}; \qquad P_1 P_n = a_n
P_{n+1} + b_n P_n
 +c_n P_{n-1};\]
\[
Q_0(x)=1; \; Q_1(x)= \frac{x- \tilde b_0}{\tilde a_0}; \qquad Q_1
Q_n = \tilde a_n Q_{n+1} + \tilde b_n Q_n + \tilde  c_n Q_{n-1}\]
with\\[0.3cm]
$a_n, \tilde a_n, c_n, \tilde c_n >0; \qquad b_n, \tilde b_n \geq
0; \qquad \tilde a_n + \tilde b_n + \tilde c_n = a_n + b_n + c_n =
1 \quad
\mbox{ for } n \in \nz$\\[0.2cm]
and: $ a_0 + b_0 = \tilde a_0 + \tilde b_0 = 1 $
\quad with $ a_0, \tilde a_0 >0$.\\[0.2cm]

\begin{Theorem} \label{s511}
The connection coefficients $c_{QP}(n,k), \; n \in \nz_0, \: k=0,
\ldots, n $, satisfy the recursion:

$\displaystyle c_{QP}(0,0) =1; \qquad
 c_{QP}(1,0) = \frac{b_0 - \tilde b_0}{\tilde a_0}; \qquad c_{QP}(1,1)=
  \frac{a_0}{\tilde a_0}$

and for $ n \in \nz $:
\begin{eqnarray*}
c_{QP}(n+1,0) & = & \frac{1}{\tilde a_0 \tilde a_n} [(b_0  -
  \tilde a_0 \tilde b_n - \tilde b_0) c_{QP}(n,0)
  +  a_0 c_1 c_{QP}(n,1) -  \tilde a_0 \tilde c_n c_{QP}(n-1,0)];\\[0.2cm]
c_{QP}(n+1,1) & = & \frac{1}{\tilde a_0 \tilde a_n} [a_0
c_{QP}(n,0)
  +  (a_0 b_1 + b_0 - \tilde a_0 \tilde b_n - \tilde b_0) c_{QP}(n,1)\\
  &&\qquad +  a_0 c_2 c_{QP}(n,2) -  \tilde a_0 \tilde c_n c_{QP}(n-1,1)]
\end{eqnarray*}

and $k=2, \ldots, n+1:$
\begin{eqnarray*}
 c_{QP}(n+1,k) & = & \frac{1}{\tilde a_0 \tilde a_n} [a_0 a_{k-1} c_{QP}(n,k-1)
  +  (a_0 b_k + b_0 - \tilde a_0 \tilde b_n - \tilde b_0) c_{QP}(n,k)\\
  &&\qquad +  a_0 c_{k+1} c_{QP}(n,k+1) -  \tilde a_0 \tilde c_n c_{QP}(n-1,k)]
\end{eqnarray*}

with $c_{QP}(s,t):=0$ for $s<t$.
\end{Theorem}

\begin{proof}
(compare \cite{AsII} for monic polynomial systems)\\[0.5cm]
$ a(0,0)=1 $ \quad holds, since \quad $ P_0(x) = Q_0(x) = 1 $.\\[0.3cm]
From $\displaystyle
 \quad  Q_1(x)=\frac{x- \tilde b_0}{\tilde a_0} = c_{QP}(1,1)
  \frac{x-b_0}{a_0} + c_{QP}(1,0)$ \quad one gets by comparing coefficients:\\[0.3cm]
  $\displaystyle \quad  c_{QP}(1,0) = \frac{b_0- \tilde
b_0}{\tilde a_0}
   \mbox{ and }  c_{QP}(1,1) = \frac{a_0}{\tilde a_0}.$

With the three term recursions one gets for $ n \geq 1:$

\begin{eqnarray}
\tilde a_n Q_{n+1} (x) & = & Q_1 (x) Q_n (x) - \tilde  b_n Q_n(x)
  -\tilde c_n Q_{n-1}(x)\nonumber\\[0.2cm]
 & = & (\frac{a_0}{\tilde a_0} P_1(x) + \frac{b_0 - \tilde b_0}{\tilde a_0})
  \sum_{k=0}^n c_{QP}(n,k) P_k (x)\nonumber\\ [0.2cm]
 && - \; \tilde b_n \sum_{k=0}^n c_{QP}(n,k) P_k(x) -
  \tilde c_n \sum_{k=0}^{n-1} c_{QP} (n-1,k) P_k(x)\nonumber\\[0.2cm]
 & = & \frac{a_0}{\tilde a_0} c_{QP} (n,0) P_1(x) + \frac{a_0}{\tilde a_0}
  \sum_{k=1}^n c_{QP} (n,k) [ a_k P_{k+1} (x) + b_k P_k (x) + c_k P_{k-1}(x)]
  \nonumber\\[0.2cm]
 && + \; \frac{b_0 - \tilde b_0 - \tilde a_0 \tilde b_n}{\tilde a_0}
   \sum_{k=0}^n c_{QP} (n,k) P_k(x) - \tilde c_n
   \sum_{k=0}^{n-1} c_{QP} (n-1,k) P_k (x). \label{gl52}
\end{eqnarray}
Setting $n=1$ gives:
\begin{eqnarray*}
Q_2 (x) & = & \frac{a_0 a_1}{\tilde a_0 \tilde a_1} c_{QP} (1,1) P_2 (x)\\[0.2cm]
&& + \; \frac{1}{\tilde a_0 \tilde a_1} [ ( a_0 b_1 + b_0 - \tilde
b_0
  - \tilde a_0 \tilde b_1) c_{QP} (1,1) + a_0 c_{QP} (1,0) ] P_1(x)\\[0.2cm]
&& + \; \frac{1}{\tilde a_0 \tilde a_1} [  a_0 c_1 c_{QP}(1,1) +
(b_0 -
  \tilde b_0 - \tilde a_0 \tilde b_1) c_{QP} (1,0) - \tilde c_1 c_{QP} (0,0)]
  P_0(x).\\[0.2cm]
\end{eqnarray*}
For $n\geq 2$ we sort the summands in equation (\ref{gl52})
according to $P_k$:
\begin{eqnarray*}
Q_{n+1}(x) & = & \frac{1}{\tilde a_0 \tilde a_n}  \Bigl( a_0 a_n
c_{QP} (n,n)
 P_{n+1}(x)\\[0.2cm]
&& + \; [ a_0 a_{n-1} c_{QP}(n,n-1) + (a_0 b_n + b_0
  - \tilde b_0 - \tilde a_0 \tilde b_n) c_{QP} (n,n)] P_n(x)\\[0.2cm]
&& + \; \sum_{k=2}^{n-1} [ a_0 a_{k-1} c_{QP}(n,k-1) + (a_0 b_k +
   b_0 - \tilde b_0 - \tilde a_0 \tilde b_n) c_{QP} (n,k) \\[0.2cm]
&& \hspace{1.1cm} + \; a_0 c_{k+1} c_{QP}(n,k+1) - \tilde a_0
\tilde c_n  c_{QP} (n-1,k)]
 P_k(x)\\[0.2cm]
&& + \; [ a_0 c_{QP}(n,0) + (a_0 b_1 + b_0 - \tilde b_0 -
   \tilde a_0 \tilde b_n) c_{QP} (n,1) \\[0.2cm]
&& \quad + \; a_0 c_2 c_{QP}(n,2) - \tilde a_0 \tilde c_n  c_{QP}
(n-1,1)]
 P_1(x)\\[0.2cm]
&& + \; [( b_0 - \tilde b_0 - \tilde a_0 \tilde b_n) c_{QP} (n,0)
  + a_0 c_1 c_{QP} (n,1)\\[0.2cm]
&& \quad - \; \tilde a_0 \tilde c_n  c_{QP} (n-1,0)]
 P_0(x) \Bigr)\\[0.2cm]
\end{eqnarray*}
This representation proves the theorem.
\end{proof}

\subsection{Proof for estimator $d_{TN}^{LS}$}

It remains to prove Theorem \ref{EstimationMain}.

\begin{Proof}
According to Proposition \ref{EstimationProp}, we have to
calculate the right side of

\[
d_{TN}^{LS}(s)=(({\bf A}_N^T{\bf A}_N)^{-1}{\bf A}_N^T{\bf X}_N)_s
\]
with the matrices of the linear model ${\bf A}_N^T{\bf A}_N=(a_{ij}^N)_{0
\leq i,j\leq N}$ defined by
\[
a_{ij}^N=\sum_{(k,l)\in I_N}g(k,l,i)g(k,l,j),\qquad 0\leq i,j\leq N.
\]
Denoting ${\bf [}x{\bf ]}$ the largest integer less then or equal to $x$, we
get for {\bf $i=j$} :\\[-0.2cm]
\begin{eqnarray*}
a_{ii}^N &=&\sum_{(k,l)\in I_N}g(k,l,i)^2=\sum_{k=0}^{{\bf [}\frac{N-i}2{\bf
]}}g(k+i,k,i)^2+\sum_{l={\bf [}\frac{i+1}2{\bf ]}}^{i-1}g(i-l,l,i)^2 \\
[0.3cm]
&=&1+{\bf [}\frac{N-i}2{\bf ]}\frac 14+(i-{\bf [}\frac{i+1}2{\bf ]})\frac
14=\left\{
\begin{array}{cl}
\frac{N+8}8 & \mbox{ $N$ even, $i$ even,} \\[0.2cm]
\frac{N+6}8 & \mbox{ $N$ even, $i$ odd,} \\[0.2cm]
\frac{N+7}8 & \mbox{ $N$ odd.}
\end{array}
\right.
\end{eqnarray*}
\\[-0.3cm]In the case {\bf $i\neq j$} we find: \\[-0.2cm]
\begin{eqnarray*}
a_{ij}^N &=&\sum_{\left. {\scriptsize \mit
\begin{array}{c}
(k,l)\in I_N \\
|k\pm l|=i \\
|k\pm l|=j
\end{array}
}\right. }g(k,l,i) g(k,l,j) \\[0.2cm]
&=&\left\{
\begin{array}{cl}
g(\frac{|j-i|}2,\frac{j+i}2,i) g(\frac{|j-i|}2,\frac{j+i}2,j)=\frac 14 &
\mbox{ \qquad for $i+j$ even} \\[0.2cm]
0 & \mbox{ \qquad else.}
\end{array}
\right.
\end{eqnarray*}
As easily can be verified, the inverse matrices
\[
({\bf A}_N^T{\bf A}_N)^{-1}=(b_{ij}^N)_{0\leq i,j\leq N}
\]
have the coefficients \\[0.2cm]
$\displaystyle
b_{ij}^N=\left\{
\begin{array}{clll}
C_{N+1}(N+3) & i=j & \mbox{ $i$ even} & \mbox{$N$ even} \\
C_{N-1}(N+1) & i=j & \mbox{ $i$ odd} & \mbox{$N$ even} \\
C_N(N+2) & i=j &  & \mbox{$N$ odd} \\
-C_{N+1} & i\neq j & \mbox{ $i+j$ even, $i$ even} & \mbox{$N$ even} \\
-C_{N-1} & i\neq j & \mbox{ $i+j$ even, $i$ odd} & \mbox{$N$ even} \\
-C_N & i\neq j & \mbox{ $i+j$ even} & \mbox{$N$ odd} \\
0 & i\neq j & \mbox{ $i+j$ odd} &
\end{array}
\right. $
\[
C_N:=\frac 8{(N+4)^2-1}\mbox{ .}
\]
{\bf $N$ odd}:\\
\begin{eqnarray*}
{\bf A}_N^T{\bf X}_N=(\sum_{(k,l)\in I_N}g(k,l,t)X_k\overline{X}_l)_{0\leq
t\leq N}=(X_0\overline{X}_t+\frac 12\sum_{I_1^N(t)}X_k\overline{X}_l)_{0\leq
t\leq N},
\end{eqnarray*}
and the above values of $b_{st}^N$ for odd $N$ yield:
\begin{eqnarray}
\lefteqn{d_{TN}^{LS}(s)=\sum_{t=0}^Nb_{st}^N (X_0\overline{X}_t+\frac
12\sum_{I_1^N(t)}X_k\overline{X}_l)}  \nonumber  \label{dosu} \\
&=&C_N\!(\!(N\!+\!2)(X_0\overline{X}_s+\frac 12\!\!\sum_{I_1^N(s)}X_k
\overline{X}_l) -\!\!\!\!\!\!\!\!\!\!\!\sum_{\left. {\scriptsize \mit
\begin{array}{c}
t=0 \\
t\neq s \\
t+s:\mbox{ even}
\end{array}
}\right. }^N\!\!\!\!\!\!\!\!\!(X_0\overline{X}_t+\frac
12\!\!\sum_{I_1^N(t)}\!\!X_k\overline{X}_l)\!).
\end{eqnarray}
\noindent
Each $(k,l)\in I_N$ lies in $I_1^N(t)$, if and only if $t=l+k$ or $t=l-k$.
Due to $k>0$ the domains $I_1^N(l+k)$ and $I_1^N(l-k)$ are not identical.
Therefore the variables $X_k\overline{X}_l$ with $k+l\neq s$ and $l-k\neq s$
as well as $l+k+s:even$ count exactly twice in
\[
\sum_{\left. {\scriptsize \mit
\begin{array}{c}
t=0 \\
t\neq s\ t+s:\mbox{ even}
\end{array}
}\right. }^N\sum_{I_1^N(t)}X_k\overline{X}_l.
\]
Note, that $t+s:even$ implies $l+k+s:even$. Further, $l+k=s$ or $l-k=s$
causes that $X_k\overline{X}_l$ occurs just once. Now, we can simplify the
double sum:

\begin{equation}
\sum_{\left. {\scriptsize \mit
\begin{array}{c}
t=0 \\
t\neq s \\
t+s:\mbox{ even}
\end{array}
}\right. }^N\sum_{I_1^N(t)}X_k\overline{X}_l=2\sum_{I_2^N(s)\cap
\{(k,l)|k>0\}}X_k\overline{X}_l+\sum_{I_1^N(s)}X_k\overline{X}_l.
\label{doppel}
\end{equation}
This result is also valid for even $N$, and applied to the right side of
(\ref{dosu}) gives us the required representation for $N$ odd.\\[0.4cm]
{\bf $N$ even}:
\\[0.3cm]We use the values of $b_{st}^N$ for even $N$ and apply (\ref{doppel})
again.\\
\noindent
a) {\it s even:}
\begin{eqnarray*}
d_{TN}^{LS}(s) &=&C_{N+1} ((N+3)(X_0\overline{X}_S+\frac 12\sum_{I_1^N(s)}X_k
\overline{X}_l)-\!\!\!\!\!\!\!\sum_{\left. {\scriptsize \mit
\begin{array}{c}
t=0 \\
t\neq s \\
t+s:\mbox{ even}
\end{array}
}\right. }^N\!\!\!\!\!\!\!(X_0\overline{X}_t+\frac 12\sum_{I_1^N(s)}X_k
\overline{X}_l)) \\
&=&C_{N+1}((N+3)(X_0\overline{X}_S+\frac{N+2}2\!\sum_{(k,l)\in I_1^N(s)}\!X_k
\overline{X}_l) -\!\sum_{(k,l)\in I_2^N(s)}\!X_k\overline{X}_l)
\end{eqnarray*}
b) {\it s odd}:
\begin{eqnarray*}
d_{TN}^{LS}(s)=C_{N-1}((N+1)X_0\overline{X}_S+\frac N2\sum_{(k,l)\in
I_1^N(s)}X_k\overline{X}_l-\sum_{(k,l)\in I_2^N(s)}X_k\overline{X}_l)
\end{eqnarray*}
To complete the proof it remains to show that the domains $I_1^N(s)$ and $
I_2^N(s)$ in a) and b) can be replaced by $
I_1^{N+1}(s),I_2^{N+1}(s) $ and $I_1^{N-1}(s),I_2^{N-1}(s)$, respectively. But,
this is a consequence of the equations\\[0.2cm]$I_1^{N+1}(s)=\;I_1^N(s)\cup
\{(k,l)|0<k\leq l,k+l=N+1,l-k=s\},$\\[0.2cm]$I_2^{N+1}(s)=\;I_2^N(s)\cup
\{(k,l)|0\leq k\leq l,k+l=N+1,l-k\neq s,k+l+s\mbox{
even}\}$.\\[0.2cm]N even and s even implies that the right sets in the
unions are empty. Thus,\\[0.1cm]$I_1^N(s)=I_1^{N+1}(s)$ and $
I_2^N(s)=I_2^{N+1}(s)$ hold. For $N$ even, $s$ odd, one gets \\[0.1cm]$
I_1^N(s)=I_1^{N-1}(s)$ and $I_2^N(s)=I_2^{N-1}(s)$, analogously.
\end{Proof}

\section{Spectral estimation}
In the classical time series analysis it is usually assumed that
the spectral measure of a weakly stationary process has the form
$\mu =f \lambda + \mu_d$ with a absolute part $f \lambda$ (here
$\lambda$ denotes the Lebesgue measure on $T$) and a discrete part
$\mu_d$ representing periods of the process.

To detect such periods from realizations $x_0, \ldots, x_{N-1}$
Schuster \cite{schu} introduced 1898 the periodogram

\[
I_N(t) = \frac{1}{2 \pi N} | \sum_{k=0}^{N-1} x_k e^{-ikt} |^2.
\]

If $\mu$ has only a finite number of atoms (periods), $I_N (t)$ is
suited to detect them. If there are no periods ($\mu_d=0$), the
periodogram is for each $t \in T$ an asymptotically unbiased
estimator of the spectral density, which lacks to be consistent.
Pointwise consistency is achieved with weighted periodograms (for
example \cite{Ros}).

The periodogram is a result of approximating the spectral measure
$\mu$ weakly with the Fej\'er series built from the covariance
function $R(n)$. Substituting in this series the canonical
estimator $\hat R(n)$ for $R(n)$ gives the periodogram $I_N(t)$
(for example \cite{Shir}, Chap. VI §4).

\subsection{Periodogram}
For $P_n$-harmonizable sequences with Haar weights satisfying (H),
we can find the periods lying in the set $S_2$ of Theorem
\ref{WienerGeneral}. Actually, the proof of this Wiener theorem
shows that the following periodogram works beyond the hypergroup
case: Replace condition $(H)$ by $(G)$ and h(n) by $m(n)^2$.
\newpage
\begin{Theorem}
Let $(X_n)_{n \in \nz_0}$ be a $P_n$-harmonizable sequence and let
(H) be satisfied. Then the estimator
\[
I_N(x,y):= \frac{1}{(\sum_{k=0}^N h(k))^2} \sum_{k,l=0}^{N} X_k
\bar{X}_l P_k (x)P_l(y) h(k) h(l), \qquad x,y \in D_s \] detects
the periods in $S_2$, in the sense that
\[ \lim_{N \to \infty} E I_N (x,y) = 0 \qquad \mbox{if } \quad (x,y) \notin
S_2 \cap \supp \mu_d;\]
\[ \limsup_{N \to \infty} E I_N (x,y) \neq 0 \qquad \mbox{for }\quad (x,y) \in
S_2 \cap \supp \mu_d.\] Here, $\mu_d$ is the discrete part of
spectral measure $\mu$.
\end{Theorem}

\begin{proof}
From the definitions we have
\begin{eqnarray*}
E I_N (x,y) & = &  \frac{1}{(\sum_{k=0}^N h(k))^2} \int_{D_s
\times D_s}
 \sum_{k,l=0}^{N} P_k(x)P_l(y)P_k(s)P_l(t) h(k) h(l) d\mu(s,t)\\\\
  &=& \int_{D_s \times D_s \backslash \{(x,y)\}} t_N(s,x) t_N(t,y) d\mu(s,t)
  +  t_N(x,x) t_N(y,y)\mu(\{(x,y)\}).
\end{eqnarray*}

As in proof of Theorem \ref{WienerGeneral}, condition (H), the
Christoffel-Darboux formula and the dominated convergence
criterion together show that
\[\int_{D_s \times D_s \backslash \{(x,y)\}} t_N(s,x) t_N(t,y)
d\mu(s,t)\]
  tends to zero for increasing N.
  Thus we have for growing N
  \[ E I_N (x,y)=t_N(x,x)
  t_N(y,y)\mu(\{(x,y)\}) + o(1),\]
  showing the theorem.
\end{proof}

\newpage
We illustrate the theorem for sequences which are harmonizable
with respect to $T_n$, the Chebyshev polynomials of the first
kind. Condition (H) is clearly satisfied and $S=D_s=[-1,1]$
(compare Remark \ref{SWiener}). According to Example 4, we choose
independent $N(0,1)$ distributed random variables $Z_0$ and $Z_1$
and define the harmonic sequence
$$X_n=T_n(-1/3) Z_0 + T_n(1/2) Z_1, \quad \nin.$$
Four realizations of the sequence:
\begin{center}
\includegraphics[width=10cm]{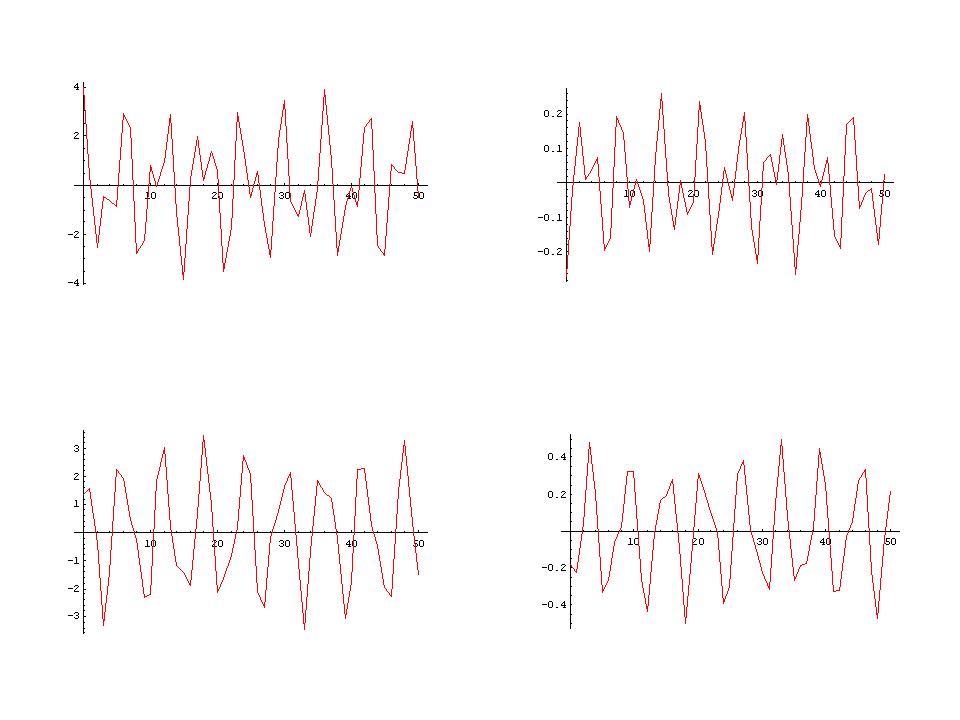}
\end{center}
The next plot shows the absolute values of the periodogram for a
realization of the first 25 variables. Peaks emerge around the
four atoms of the spectral measure:
\begin{center}
\includegraphics[width=8cm]{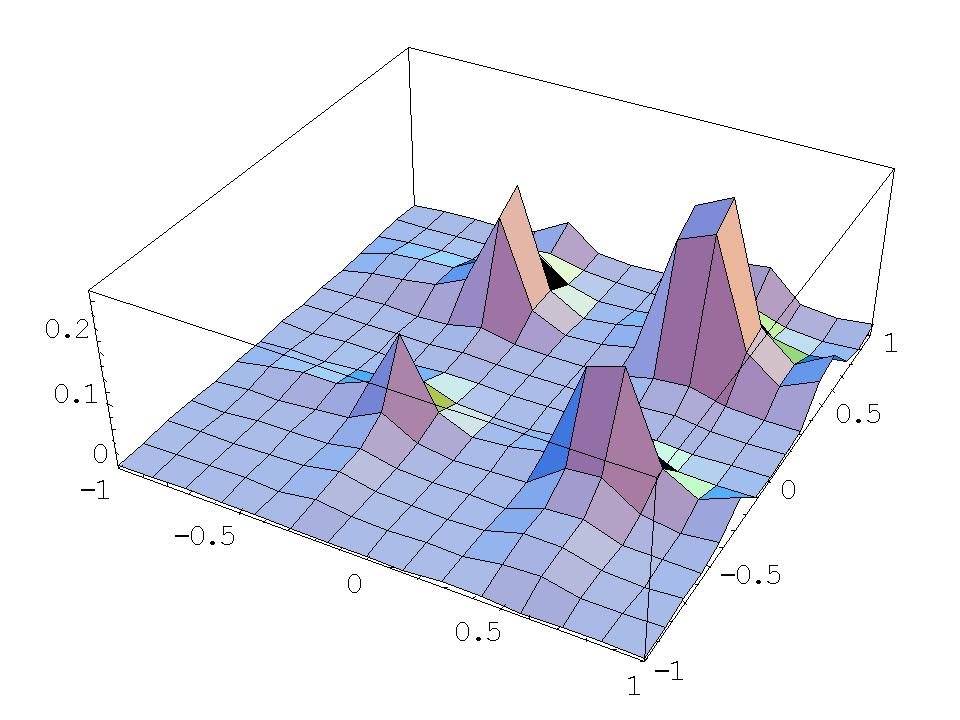}
\end{center}

\subsection{Density estimation }
In this section we discuss the estimation of the spectral density
of $P_n$-weakly stationary sequences. The methods use Fej\'er
series in certain homogeneous Banach spaces $B$ (compare Lasser at
al. \cite{laob}). We restrict our discussion to the $P_n$-weakly
stationary case were we will present estimation strategies.

Let $(X_n)_{n \in \nz_0}$ be a $P_n$-weakly stationary process
with a spectral measure $\mu = f \pi$ which is absolutely
continuous with respect to the orthogonalization measure $\pi$ and
has a spectral density $f$ which lies one of the Banach spaces
$B=L^p(D_s), \quad 1 \leq p< \infty,$ or $\;B= C(D_s)$.

Further, assume that $a_{Ns}$ are complex numbers with $|a_{Ns}|
\leq A$ for some $A>0$ and all $N \in \nz_0,\; s=0,\ldots,N$ and
the property that the series

\[
A_N(g) := \sum_{s=0}^N a_{Ns} \check{g} (s) P_s h(s) \mbox{\quad
with \quad}
 \check{g} (s) = \int_{D_s} g P_s d \pi \mbox{\quad for } g \in B\]
 approximates $g$ in the norm of $B$:\\[0.2cm]
\[\lim_{N \to \infty} || A_N (g) - g ||_B =0.\]

We define for an estimator $\hat d$ of the covariance function of
$(X_n)_{n \in \nz_0} $ and a sequence $(q(N))_{N \in \nz_0}$ with
$q(N) \in \nz_0$ the following function:

\[
\hat f_N (A_N, \hat d, q(N)):= \sum_{s=0}^N a_{Ns} \hat
d_{q(N)}(s) P_s h(s), \qquad N \in \nz_0.\]

In this situation we can prove the following:

\begin{Theorem} \label{s221}

If $(q(N))_{N \in \nz_0} $ is a sequence of positive integers such
that the condition \\
\[ \max_{s \leq N} \{ E| \hat d_{q(N)} (s) - d(s)|^2\} = o
(\frac{1}{(\sum_{k=0}^N h(k))^2}) \]
is satisfied, then it is true that:\\
\[
\lim_{N \to \infty} E || \hat f_N (A_N, \hat d, q(N))- f||_B=0.\]
\end{Theorem}

\begin{proof}

We denote \quad $||\;||:=||\;||_B; \qquad \hat f_N:= \hat f_N
(A_N, \hat d,
q(N))$.\\[0.2cm]
The triangle inequality shows:
\begin{eqnarray*}
 E|| \hat f_N -f|| &\leq& E||\hat f_N-E\hat f_N||+||E\hat f_N-
f||;\\
||E \hat f_N -f|| &\leq & || A_N(f) -f||+ ||\sum_{s=0}^N E (\hat
d_{q(N)} (s) - d(s)) a_{Ns} P_s h(s)||\\
&\leq& || A_N(f)- f||+(\sum_{s=0}^N h(s)) A ||1|| \max_{s \leq N}
\{ | E (\hat d_{q(N)} (s) - d(s))|\}.\\
\end{eqnarray*}

From the Cauchy-Schwarzs inequality we get for all $x \in D_s$:\\
\begin{eqnarray*}
\lefteqn{|\sum_{s=0}^N (\hat d_{q(N)}(s) - E \hat d_{q(N)} (s))
a_{Ns} P_s(x) h(s)|}\\
&&\leq \Bigr(\sum_{s=0}^N (\hat d_{q(N)} (s) - E\hat d_{q(N)}
(s))^2 h(s)\Bigl)^{\frac{1}{2}} \Bigr( \sum_{s=0}^N (a_{Ns}
P_s(x))^2
h(s)\Bigl)^\frac{1}{2}.\\
\end{eqnarray*}

The considered Banach spaces $B$ have a monotone norm.\\
$|g_1| \leq |g_2|$ f\"{u}r $g_1,g_2 \in B$ implies $||g_1||\leq
||g_2||$, and therefore Cauchy-Schwarz and Jensen inequality yield
(\cite{rudre}, Th. 3.3):
\begin{eqnarray*}
\lefteqn{E|| \hat f_N - E \hat f_N|| = E || \sum_{s=0}^N (\hat
d_{q(N)} (s) - E \hat d_{q(N)} (s)) a_{Ns} P_s h(s)||}\\
&& \leq E \Bigl(\sum_{s=0}^N (\hat d_{q(N)} (s) - E \hat d_{q(N)}
(s))^2
h(s)\Bigl)^\frac{1}{2} \, A \,||1||\, \Bigr(\sum_{s=0}^N h(s)\Bigl)^\frac{1}{2}\\
&& \leq \Bigl(\sum_{s=0}^N E (\hat d_{q(N)} - E \hat d_{q(N)}
(s))^2
h(s)\Bigl)^\frac{1}{2} \, A \, ||1||\, \Bigr(\sum_{s=0}^N h(s)\Bigl)^\frac{1}{2}\\
&& \leq \Bigl(\sum_{s=0}^N h(s) \Bigr) \,A \, ||1||\, \max_{s \leq
N} \{ \Bigl( Var \, \hat
d_{q(N)}(s)\Bigr)^\frac{1}{2}\}.\\
\end{eqnarray*}
From the equation \quad $\displaystyle E |\hat d_{q(N)}
(s)-d(s)|^2 = Var (\hat
d_{q(N)}(s)) +(E (\hat d_{q(N)} (s)-d(s))^2$\\[0.2cm]

it follows that\\[0.2cm]
$\displaystyle \max_{s \leq N} \{ | E ( \hat d_{q(N)} (s)-d(s))|\}
+ \max_{s \leq N} \{ ( Var \, \hat d_{q(N)} (s))^\frac{1}{2}\}\\[0.2cm]
\hspace*{2cm} \leq 2 \Bigl( \max_{s \leq N} \{ E |\hat d_{q(N)}
(s)-
d(s)|^2\} \Bigr)^{1/2}$.\\[0.2cm]

Thus, we also have  \\

$\displaystyle \hspace*{1cm} E|| \hat f_N-f|| \leq || A_N(f) -f||
+ ( \sum_{s=0}^N h(s))
2  A  \Bigl( \max_{s \le N} \{ E| \hat d_{q(N)} (s)-d(s)|^2\} \Bigr)^{1/2}$\\

which together with the assumption finishes the proof.
\end{proof}

Apparently, if the estimator $\hat d$ is consistent we can choose
sequences $q(N)$ which satisfy the assumptions of the preceding
theorem. In general one does not know the necessary growth order.
For the praxis of spectral estimation it will be recommended to
work with a $N$ which is "small" with respect to the number of
observations $x_0, \ldots ,x_{q(N)}$. This means, to require for
an estimators $\hat f_N$ "good" estimations $\hat d_{q(N)}$ of
$d$.

\chapter{Prediction Theory}

\section{One and n-step prediction}
We now study the prediction problem. In the sequel assume that the
polynomial sequence $(P_n)_{\nin}$ induces a polynomial
hypergroup. For a $P_n$-weakly stationary random sequence
$(X_n)_{\nin}$ the problem in one-step prediction can be
formulated as follows. Given $\nin$ denote the linear space
generated by $X_0,...,X_n$ by $\; H_n:= \; $ span $ \{ X_0,...,X_n
\} \; \subseteq \; H = \overline {span \{ X_n: \nin \}}.\;$ We
want to characterize the prediction $\hat{X}_{n+1} \in H_n $ of
$X_{n+1}$ with the minimum property
\begin{equation}
\parallel \hat{X}_{n+1} -X_{n+1} \parallel_2 \; =
\, \min \{ \parallel Y-X_{n+1} \parallel_2: \; Y \in H_n \}.
\end{equation}
It is well-known that $\hat{X}_{n+1} \, = P_{H_n} X_{n+1}$, where
$P_{H_n}$ is the orthogonal projection from $L^2(P) $ to $H_n.\;$
The problem is to determine the coefficient $b_{n,k},\;k=0,...,n
\;$ in the representation
\begin{equation}
\hat{X}_{n+1} \; = \; \sum_{k=0}^n b_{n,k} X_k
\end{equation}
and to decide whether the prediction error
\begin{equation}
\delta_n := \; \parallel \hat{X}_{n+1} -X_{n+1} \parallel_2
\end{equation}
converges to zero as $n$ tends to infinity.

\begin{Definition} A $P_n$-weakly stationary sequence $(X_n)_{\nin}$
is called asymptotically $P_n$-deterministic if the prediction
error $\delta_n$ satisfies
\begin{equation}\label{asympdet}
    \lim_{n \rightarrow \infty} \delta_n =0.
\end{equation}
\end{Definition}

From Hilbert space theory we know that $\hat{X}_{n+1}$ can be
characterized by the property
\begin{equation}\label{15}
E\left( (X_{n+1} - \hat{X}_{n+1}) \overline{Y} \right) \;= \;
\langle X_{n+1} - \hat{X}_{n+1},Y \rangle \;=0
\end{equation}
for all $Y \in H_n.\;$ From \ref{15} we get for
$b=(b_{n,0},...,b_{n,n})^T $ the linear equation
\begin{equation}
\Phi^T b \;=\; \varphi,
\end{equation}
where $\varphi \;= \left( E(X_{n+1} \overline{X_0}),E(X_{n+1}
\overline{X}_1),...,E(X_{n+1}\overline{X_n}) \right)^T \;$ and
$\Phi$ is the $(n+1) \times (n+1)$-matrix $ \Phi \;= \left(
E(X_i\overline{X_j}) \right)_{0 \leq i,j \leq n}.\;$ Matrices
structured like $\Phi$ can be seen as generalizations of
Toeplitz-matrices with regard to the orthogonal polynomial
sequence $(P_n)_{\nin}.\;$ The dependency of the $i,j$- entry on
$i-j \;$ in the Toeplitz case is substituted by the condition $\;
\Phi_{i,j}=\varepsilon_{i}\ast\varepsilon_{j}
(d)=\sum\limits_{k=|i-j|}^{i+j} g(i,j;k) \; d(k)\;$ with $d(k) \;=
E(X_n\overline{X_0}).\;$

A Methods of fast inversion of such matrices resembling the
Durbin-Levinson algorithm (compare Brockwell and Davies
\cite{Bro}), will be presented in the last section of the chapter.

\subsection{Best linear Predictors}

Using the isometric isomorphism between $H$ and $L^2(D_s,\mu)$,
the problem of finding the coefficients  $b_{nk},\; k=0,
\ldots,n$, in the representation
\[
\hat X_{n+1} = \sum_{k=0}^n b_{nk} X_k \] can be transferred to
orthogonal polynomials:

Let $\Pi^1_{n+1}$ be the set of all polynomials $Q(x)$ of degree
$n+1$ and leading coefficients $1$ with respect to $P_{n+1}\;$,
that is $Q(x) = P_{n+1}(x) - \sum_{k=0}^n \beta_{nk} P_k (x)$. We
determine the polynomial $Q_{n+1}^*(x) \in \Pi_{n+1}^1$ with the
property:
\[
|| Q_{n+1}^*||_2 = \min_{Q \in \Pi^1_{n+1}} \{ || Q||_2 \}.\]

 The
coefficients $b_{nk}$ of the predictor $\hat X_{n+1}$ and the
prediction error $\delta_n$ can be derived from the representation

\[Q_{n+1}^* (x) = P_{n+1}(x) - \sum_{k=0}^n b_{nk} P_k(x).\]

\begin{Theorem} \label{s311}

Let $(q_n)_{n \in \nz_0}$ be the orthonormal polynomial sequence
with respect to the spectral measure $\mu$ and denote $\varrho_n
=\varrho_n (\mu)$ the positive leading coefficients of $q_n(x) $,
i.e. $q_n(x)= \varrho_n x^n+ ...$. Further, denote $\sigma_n
=\sigma_n(\pi)$ the leading coefficients of $P_n(x)$.

\begin{enumerate}
\item The minimizing polynomial $Q_{n+1}^* (x) \in \Pi^1_{n+1}$ of
the prediction is given by

\begin{equation}
  Q_{n+1}^* (x) = \frac{\sigma_{n+1}}{\varrho_{n+1}} q_{n+1}(x);
  \end{equation}
\item The best linear one step predictor has the representation:
 \[ \hat X_{n+1} = \int_{D_s}( P_{n+1} (x) -
    \frac{\sigma_{n+1}}{\varrho_{n+1}} q_{n+1}(x)) dZ(x);\]
\item The prediction error satisfies
 \[
  \delta_n = \frac{\sigma_{n+1}}{\varrho_{n+1}};\]
\item $\hat X_{n+1}$ has the innovation representation:
  \[
   \hat X_{n+1} = \sum_{k=0}^n < P_{n+1}, q_k > \frac{\varrho_k}{\sigma_k}
    (X_k-\hat X_k)\]
  with the scalar product $<\; , \;>$ in $L^2(D_s,\mu)$.\\
\end{enumerate}
\end{Theorem}

\begin{Proof}
All polynomials $Q \in \Pi^1_{n+1}$ can be written as
\[
 Q(x) = \sum_{k=0}^n \beta_k q_k (x) + \frac{\sigma_{n+1}}{\varrho_{n+1}}
 q_{n+1}(x).\]

Hence
\[
||Q||_2^2 = \int_{D_s} | Q(x)|^2 d \mu(x) = | \beta_0|^2+ ... +
|\beta_n|^2
 + | \frac{\sigma_{n+1}}{\varrho_{n+1}}|^2.\]
$||Q||_2^2$ becomes minimum for $\beta_0 = ... = \beta_n = 0$ and
a), b), c) follow:
\[
Q_{n+1}^* (x) = \frac{\sigma_{n+1}}{\varrho_{n+1}} q_{n+1} (x);\]
\[\delta_n=|| Q_{n+1}^* ||_2 = \frac{\sigma_{n+1}}{\varrho_{n+1}};\]
\[\hat X_{n+1} = X_{n+1} - (X_{n+1} - \hat X_{n+1}) =
\int_{D_s} (P_{n+1} (x) - \frac{\sigma_{n+1}}{\varrho_{n+1}}
q_{n+1}(x)) d Z(x).\]

(d) results from the Fourier expansion of $P_{n+1} (x) - Q_{n+1}^*
(x)$ with respect to $(q_n)_{n \in \nz_0}$:
\begin{eqnarray*}
 \hat X_{n+1} &=& \int_{D_s} \sum_{k=0}^n < P_{n+1} - Q_{n+1}^*, q_k> q_k(x)
    d Z (x)\\[0.2cm]
 &=& \sum_{k=0}^n \int_{D_s} < P_{n+1}, q_k> q_k (x) d Z (x) \\[0.2cm]
 &=& \sum_{k=0}^n < P_{n+1}, q_k> \frac{\varrho_k}{\sigma_k} (X_k - \hat X_k).
 \end{eqnarray*}
 \end{Proof}

Linear n-step prediction can be discussed likewise. The
representation of the optimal predictor $\hat X_{n+m}^{(m)} \;
(m>1)$ comprises besides the leading coefficients further
coefficients of the involved polynomial systems. Accordingly
complicated are investigations of the asymptotic of the prediction
error $\delta_n^{(m)}$. We give the representation of $\hat
X_{n+m}^{(m)}$ and $\delta_n^{(m)}$, but focus in the sequel on
the one step prediction error.

\newpage
\begin{Theorem}
With the notations as in Theorem \ref{s311} one has:
\begin{enumerate}
\item The best linear m step predictor $\hat X_{n+m}^{(m)}$
satisfies:
\begin{equation} \label{gl335}
\hat X_{n+m}^{(m)} = \int_{D_s} (P_{n+m} (x) - \sum_{k=n+1}^{n+m}
c_{Pq} (n+m,k) \: q_k (x)) d Z(x)
\end{equation}
with the connection coefficients $c_{Pq}(n+m,k)$ providing
$P_{n+m}$ from $q_k$.\\
\item The m-step prediction error $\delta_n^{(m)} := ||\hat
X_{n+m}^{(m)} - X_{n+m} ||_2$ is given by:
\begin{equation}
 \delta_n^{(m)} = \sqrt{\sum_{k=n+1}^{n+m} c_{Pq} (n+m,k)^2}.
 \end{equation}
\end{enumerate}

\end{Theorem}
\begin{Proof}
We consider the class of
polynomials $\Pi_{n,n+m}^1$ having the form\\
$\displaystyle Q(x) = P_{n+m}(x)-\sum_{k=0}^n \beta_{nk} P_k(x)$
for $\beta_{nk} \in \rz$ ($n \! \in \! \nz_0, \: k=0,...,n$). Now,
we look for the polynomial minimizing $\displaystyle \int_{D_s}
|Q(x)|^2 d \mu(x)$. Every $Q(X)$ can be written as
\[
 Q(x) = \sum_{k=n+1}^{n+m} c_{Pq}(n+m,k) q_k(x)+\sum_{k=0}^n
 \alpha_{nk} q_k(x)\]
with some $\alpha_{nk} \in \rz,\: n \in \nz_0, \; k=0,...,n$.\\[0.3cm]
Due to $\displaystyle \int_{D_s} | Q(x)|^2 d \mu(x)=
\sum_{k=n+1}^{n+m} | c_{Pq} (n+m,k)|^2+ \sum_{k=0}^n |
\alpha_{nk}|^2$ the minimum is taken for $\alpha_{nk}=0, \;
k=0,...,n$, which implies (a) and (b).
\end{Proof}

The rest of the section is devoted a thorough investigation of the
asymptotic of the one step predictor
\begin{equation}\label{17}
    \delta_n = \frac{\sigma_{n+1}}{\varrho_{n+1}}.
\end{equation}
If the orthogonalization measure $\pi$ satisfies the
Kolmogorov-Szeg\"o property, we can give a rather complete
characterization when $\delta_n \rightarrow 0 $ holds.

\subsection{Kolmogorov-Szeg\"o class}

To each polynomial sequence $(p_n)_{\nin}$ orthonormal with
respect to a probability measure $\nu \in M^1 ([-1,1])$ one can
associate a unique polynomial sequence $(\psi_n)_{\nin}$ on $
[-\pi,\pi]$ orthonormal with respect to a measure $\alpha$ given
by
\begin{equation}\label{18}
d \alpha (t) \;=\; |\sin t | \; d \nu (\cos t),\qquad \qquad t \in
\; [-\pi,\pi].
\end{equation}
Denote the positive leading coefficients of $\psi_n(t) $ by
$\varrho_n (\alpha) $ and those of $p_n(t)$ by $\varrho_n
(\nu).\;$ We apply some important results on orthogonal
polynomials on the unit circle. The original references are
Geronimus \cite{Ger} and Szeg\"{o} \cite{Szeg}, see also Lubinsky
\cite{Lub}.

If the Radon Nikodym derivative $\alpha'$ of $\alpha$ fulfils the
Kolmogorov-Szeg\"o property, i.e. $\ln (\alpha') \;\in L^1
(\,[-\pi,\pi]\,) \;$ or equivalently
\begin{equation}
\int_{-\pi}^\pi \ln (\alpha'(t)) \; dt \; > -  \infty,
\end{equation}
then $\varrho_n (\alpha) $ converges monotonically increasingly
towards the geometric mean
\begin{equation}
\varrho(\alpha) := \; \exp \left( \;- \frac{1}{4 \pi}
\;\int_{-\pi}^\pi \ln (\alpha'(t)) \; dt\right),
\end{equation}
see (\cite{Lub} Theorem 3.4). If the Kolmogorov-Szeg\"o property
is not valid we have $\lim\limits_{n \rightarrow \infty} \varrho_n
(\alpha) \;= \infty.$ Transferring this result to $[-1,1]$ we
have, see (\cite{Ger}, Theorem 9.2)

\begin{Proposition} \label{Prop4.1}
Let $(p_n)_{\nin}$ be an orthonormal polynomial sequence with
respect to a measure $\nu \in M^1 ([-1,1]), \; \; (|\supp \nu | \;
= \infty).\;$ Then alternatively we have:
\begin{enumerate}
\item[(i)] If $\displaystyle\int_{-1}^1 \displaystyle\frac{\ln
(\nu' (x))}{\sqrt{1-x^2}} \; dx \; > - \infty,\;$ then there exist
positive constants $C_1, C_2$ such that for all $\nin$
$$ C_1 \; \leq \; \frac{\varrho_n (\nu)}{2^n} \;\leq \;C_2.$$
\item[(ii)] If $ \displaystyle\int_{-1}^1 \displaystyle\frac{\ln
(\nu' (x))}{\sqrt{1-x^2}} \; dx \; = - \infty, \;$ then $
\lim\limits_{n \rightarrow \infty} \displaystyle\frac{\varrho_n
(\nu)}{2^n} \; = \infty.$
\end{enumerate}
\end{Proposition}

\begin{proof}
Consider the orthonormal polynomial sequence $(\psi_n)_{\nin}$ on
$[-\pi,\pi]$ with leading coefficients $\varrho_n (\alpha),\;$
where $\alpha$ is the measure as in \ref{18}. The coefficients
$\varrho_n(\alpha)$ and $\varrho_n (\nu)$ satisfy the
inequalities, see (\cite{Ger}, (9.9)):
\begin{equation}
\frac{\varrho_{2n-1} (\alpha)}{2 \sqrt{\pi}} \; \leq \;
\frac{\varrho_n(\nu)}{2^n} \; \leq \; \frac{\varrho_{2n}
(\alpha)}{\sqrt{\pi}}.
\end{equation}
Moreover we have
$$ \int_{-\pi}^\pi \ln (\alpha' (t)) \; dt \; = \; 2 \int_{-1}^1 \frac{\ln (\nu' (x))}{\sqrt{1-x^2}} \; dx,
$$
and the assertions (i) and (ii) follow by the convergence of
$\varrho_n (\alpha)$ towards $\varrho (\alpha)$ or towards
infinity, respectively.
\end{proof}

For the orthonormal version $p_n(x)$ of $P_n(x)$ we have $p_n(x)
\;= \sqrt{h(n)}\; P_n(x).\;$ Hence $ \sigma_n(\pi) \;=
\displaystyle\frac{\varrho_n(\pi)}{\sqrt{h(n)}} \;$ and by
(\ref{17}) we obtain immediately:

\begin{Theorem} \label{4.2}
Let $(X_n)_{\nin}$ be a weakly stationary random sequence on the
polynomial hypergroup ${\Bbb N}_0$ induced by $(P_n)_{\nin}.\;$
Suppose that the orthogonalization measure $\pi$ of $(P_n)_{\nin}$
has support contained in $[-1,1]$ and fulfils the
Kolmogorov-Szeg\"o property on $[-1,1]$, i.e.
\begin{equation}
\int_{-1}^1 \frac{\ln (\pi' (x))}{\sqrt{1-x^2}} \; dx \; > \; -
\infty.
\end{equation}
Then $(X_n)_{\nin}$ is asymptotically $P_n$-deterministic if and
only if
\begin{equation}
\lim_{n \rightarrow \infty} \frac{\sqrt{h(n)} \; \varrho_n
(\mu)}{2^n} \; = \;\infty.
\end{equation}
In particular we have:
\begin{enumerate}
\item[(i)] If $\lim\limits_{n \rightarrow \infty} h(n) \; = \infty
\;$, then $(X_n)_{\nin}$ is asymptotically $P_n$-deterministic
without any assumption on the spectral measure $\mu.\;$ For the
prediction errors we have
$$ \delta_n \; = \; O \left( \frac{1}{\sqrt{h(n+1)}} \right) \qquad \qquad \mbox{as } n \rightarrow \infty. $$
\item[(ii)] If $\{ h(n): \; \nin \} \;$ is bounded, then
$(X_n)_{\nin}$ is asymptotically $P_n$-deterministic if and only
if $\mu$ does not fulfil the Kolmogorov-Szeg\"o property on
$[-1,1].$ \item[(iii)] If $\{ h(n):\; \nin \}\;$ is unbounded,
then $(X_n)_{\nin}$ is asymptotically $P_n$-deterministic provided
$\mu$ does not fulfil the Kolmogorov-Szeg\"o property on $[-1,1].$
\end{enumerate}
\end{Theorem}

\begin{proof}
We know that $\delta_n \; = \displaystyle\frac{\sigma_{n+1}
(\pi)}{\varrho_{n+1} (\mu)} \; = \displaystyle\frac{\varrho_{n+1}
(\pi)}{\sqrt{h(n+1)} \; \varrho_{n+1} (\mu)}.\;$

By the Kolmogorov-Szeg\"o property for $\pi$ and Proposition
\ref{Prop4.1} there exist positive constants $C_1,C_2$ such that
$$ C_1 \; \leq \; \frac{\varrho_{n+1} (\pi)}{2^{n+1}} \;\leq \; C_2 \qquad \quad \mbox{ for } \nin,$$
and hence
$$ C_1 \; \leq \; \sqrt{h(n+1)} \; \frac{ \varrho_{n+1} (\mu)}{2^{n+1} } \;
\delta_n \; \leq \; C_2 \qquad \mbox{ for } \nin. $$ Now applying
\ref{Prop4.1} again every statement follows.
\end{proof}

From the preceding theorem we know that $(X_n)_{\nin}$ shows a
prediction behaviour similar to classical weakly stationary
processes provided $\{ h(n): \; \nin \} \;$ is bounded and the
orthogonalization measure $\pi$ fulfils the Kolmogorov-Szeg\"o
condition. We now prove for the case $ D_s = [-1,1]\;$ that if
$(h(n))_{\nin}$ does not converge to infinity, the
Kolmogorov-Szeg\"o property of $\pi$ is valid.

\begin{Proposition}\label{Prop4.3}
Let $(P_n)_{\nin}$ induce a polynomial hypergroup on ${\Bbb N}_0$,
and suppose that $ D_s = [-1,1].\;$ If $\pi$ does not fulfil the
Kolmogorov-Szeg\"o property, then $\lim\limits_{n \rightarrow
\infty} h(n) \; = \infty.$
\end{Proposition}

\begin{proof}
Consider the representation of $P_n(x)$ by Chebyshev polynomials
of the first kind
\begin{equation}\label{24}
P_n(x) \;=\; \sum_{k=0}^n c_{PT}(n,k) T_k(x),
\end{equation}
where $T_k(x) \; = \cos ( k \arccos x) \;= 2^{k-1} x^k + \ldots
\,.\;$ Since $ |P_n(x)| \; \leq 1 \;$ for every $x \in D_s, \;
\nin, \;$ we obtain by applying the orthogonalization measure $d
\pi_T (x) \;= \displaystyle\frac{1}{\pi} \;
\displaystyle\frac{dx}{\sqrt{1-x^2}} \;$ of $(T_k)_{k \in {\Bbb
N}_0}$
$$ 1 \; \geq \int_{-1}^1 P_n^2 (x) \; d\pi_T (x) \; \geq a_{n,n}^2 \;
\int_{-1}^1 T_n^2 (x) \; d \pi (x) \;= \; \frac{a_{n,n}^2}{2}\; .
$$ Comparing the leading coefficients in \ref{24}, we get
$$ a_{n,n} \;= \; \frac{\sigma_n(\pi) }{ 2^{n-1}} \; =\;
\frac{\varrho_n(\pi)}{\sqrt{h(n)} \; 2^{n-1}},$$ and hence $2 h(n)
\; \geq \left( \varrho_n(\pi) / 2^{n-1}\right)^2.$

By Proposition \ref{Prop4.1} we have $\lim\limits_{n \rightarrow
\infty} \varrho_n (\pi)/2^n \; = \infty \;$ and hence
$\lim\limits_{n \rightarrow \infty} h(n) \; = \infty.$
\end{proof}

We investigate some examples to determine the asymptotic behavior
of the prediction error.
\begin{enumerate}
\item[(a)] Let $P_n^{(\alpha, \beta)}$ be the Jacobi polynomials
with $\alpha \geq \beta > -1 \;$ and $ \alpha + \beta + 1 \geq 0.
\;$ The Haar weights are
$$ h(n) \;=\; \frac{(2n+ \alpha + \beta +1) \, (\alpha + \beta +1)_n (\alpha+1)_n}
{(\alpha + \beta +1) \; n! \; (\beta +1)_n}\; .$$ Evidently $\pi$
satisfies the Kolmogorov-Szeg\"o property. If $\alpha \neq -
\frac{1}{2} $ we have $\lim\limits_{n \rightarrow \infty} h(n) \;
= \infty. \;$ Hence every weakly stationary random sequence on the
polynomial hypergroup ${\Bbb N}_0$ induced by
$(P_n^{(\alpha,\beta)})_{\nin}$ is asymptotically
$P_n^{(\alpha,\beta)}$-deterministic.

For the prediction error we have
$$ \delta_n \; = \; O(n^{-\alpha - \frac{1}{2}}) \qquad\quad  \mbox{ as } n \rightarrow \infty.$$
If $\alpha =- \frac{1}{2},$ then $P_n^{(-1/2,-1/2)} (x) \;=T_n(x)
\;$ are the Chebyshev polynomials of the first kind. The Haar
weights are $h(0)=1, \; h(n)=2 \;$ for $n \in {\Bbb N},\;$ and
Theorem \ref{4.2} (ii) can be utilized. \item[(b)] $P_n^{(\nu)}
(x;\alpha) \;$ the associated ultraspherical polynomials with
$\alpha > -\frac{1}{2},\;\;\nu \geq 0.\;$ These polynomials are
studied in detail in (\cite{la02}, (3)). The Haar weights are
$$ h(n) \;=\; \frac{(2n+2\alpha+2\nu+1)}{4 \alpha^2 (2 \alpha + 2\nu +1)(\nu+1)_n(2\alpha + \nu +1)_n} \;
\left( (2 \alpha + \nu)_{n+1} - \; (\nu)_{n+1} \right)^2$$ and the
orthogonalization measure $\pi$ on $[-1,1]$
$$d\pi (x) \;=\; c_{\alpha \nu}\;  g(x) \; (1-x^2)^\alpha dx \qquad \mbox{ with}$$
$$ g(\cos t) \;=\; \left| _2F_1 (\frac{1}{2}-\alpha, \nu; \nu+\alpha+ \frac{1}{2}; e^{2it})
\right|^{-2},\qquad x = \cos t. $$
To show that $\pi$ fulfils the
Kolmogorov-Szeg\"o property it is convenient to look at the
leading coefficients. From (\cite{la02}, (3.10)) we have
$$ \varrho_n(\pi) \;=\; 2^n \; \left( \frac{(\nu + \alpha + \frac{3}{2})_n
(\nu + \alpha + \frac{1}{2})_n }{(\nu + 1)_n (\nu+2\alpha +1)_n}
\right)^{1/2} $$ and the asymptotic properties of the Gamma
function yield
$$ \lim_{n \rightarrow \infty} \frac{\varrho_n(\pi)}{2^n} \;=\;
\left( \frac{\Gamma(\nu+1) \; \Gamma (\nu+2\alpha+1)}{\Gamma
(\nu+\alpha+\frac{3}{2}) \; \Gamma(\nu+\alpha+\frac{1}{2})}
\right)^{1/2}. $$ Proposition \ref{Prop4.1} yields the
Kolmogorov-Szeg\"o property of $\pi$. Using once more the
asymptotic properties of the Gamma function we see that $h(n)
=O(n^{2 \alpha +1}) \;$ and hence we have for the prediction error
$$\delta_n \; = \; O(n^{-\alpha-\frac{1}{2}}) \qquad \quad \mbox{ as } n \rightarrow \infty. $$
Note that we enlarge the domain of $\alpha$ compared to H\"osel
and Lasser (\cite{hoe}, Corollary 2). \item[(c)] $P_n(x ;
\nu,\kappa) \,$ the Bernstein-Szeg\"o polynomials with $\nu,\,
\kappa \geq 0,\;\kappa-1< \nu <1.$ These polynomials under
consideration are orthogonal with respect to the measure on
$[-1,1]$
$$ d \pi (x) \; = \; c_{\nu \kappa} \; \frac{dx}{g(x) \; \sqrt{ 1-x^2}},$$
where $g(x) \; = |\nu e^{2it} + \kappa e^{it} +1|^2,\;\; x=\cos
t,\;$ is a polynomial with $g(x)>0$ for all $x \in [-1,1].\;$ By
Szeg\"o \cite{Szeg} these polynomials can be represented
explicitly by Chebyshev polynomials of the first kind
$$ \begin{array}{lll}
P_n(x;\nu,\kappa) &=& \displaystyle\frac{1}{\nu+\kappa+1} \;
\left( T_n(x) + \kappa T_{n-1}(x) + \nu T_{n-2}(x)\right),\quad n \geq 2,\\
&&\\
P_1(x;\nu,\kappa) &=& \displaystyle\frac{1}{\nu+\kappa+1} \;
\left( (\nu+1)T_1(x) + \kappa T_0(x) \right), \\
&&\\
P_0(x;\nu,\kappa) &=&1.
\end{array} $$
An easy calculation shows
$$ P_1(x;\nu,\kappa) \; P_n(x;\nu,\kappa) \; = \; \frac{\nu+1}{2(\nu+\kappa+1)} \;
P_{n+1} (x;\nu,\kappa) \qquad $$
\begin{equation}\label{25}
+ \; \frac{\kappa}{\nu+\kappa+1} \; P_n(x;\nu,\kappa) \; + \;
\frac{\nu+1}{2(\nu+\kappa+1)} \; P_{n-1}(x;\nu,\kappa)
\end{equation}
for $n \geq 3.\;$ It is straightforward to calculate the
linearization coefficients and to check directly that
$\left(P_n(x;\nu,\kappa)\right)_{\nin}$ induces a polynomial
hypergroup on ${\Bbb N}_0$ provided $\nu,\kappa \geq 0$ and
$\kappa-1 < \nu<1.\;$ Since the recurrence coefficients in
\ref{25} are constant for $n \geq 3$ the Haar weights are bounded,
and Theorem \ref{4.2}(ii) can be applied.
\end{enumerate}

\subsection{The general case}

There are important orthogonal polynomial sequences inducing a
polynomial hypergroup but not belonging to the Kolmogorov-Szeg\"o
class, e.g. the Cartier-Dunau polynomials discussed before. Often
it is not possible or at least very difficult to decide the
membership to the Kolmogorov-Szeg\"o class. Throughout this
subsection we suppose on $\pi$ only that the corresponding
orthogonal polynomials induce a polynomial hypergroup on $ {\Bbb
N}_0.\;\; (X_n)_{\nin} $ will be a weakly stationary random
sequence on this polynomial hypergroup with spectral measure
$\mu.\;$

A standard procedure in the Hilbert space $L^2(D_s,\mu)$ yields
\begin{equation}
\delta_n^2 \; = \;
\frac{\Delta(X_0,...,X_{n+1})}{\Delta(X_0,..,X_n)},
\end{equation}
where
$$ \Delta(X_n,...,X_{n+m}) \; = \; \mbox{ det } \; \left( \begin{array}{ccc}
\langle X_n,X_n \rangle & \cdots & \langle X_n,X_{n+m }\rangle \\
\vdots && \vdots \\
\langle X_{n+m},X_n \rangle &\cdots & \langle X_{n+m},X_{n+m}
\rangle
\end{array} \right) $$
Notice that $\Delta (X_0,...,X_n) \; >0 \;$ for all $n \in {\Bbb
N}_0 $ if and only if $\supp \mu$ is infinite. We assume further
on that $|\supp \mu | \; = \infty.$

The following sequence of upper bounds holds for $\delta_n.$
$$ \delta_n^2 \; = \; \frac{\Delta(X_0,...,X_{n+1})}{\Delta(X_0,...,X_n)} \; \leq \;
\frac{\Delta(X_1,...,X_{n+1})}{\Delta(X_1,...,X_n)} \; \leq $$
\begin{equation}\label{27}
\leq \; ... \; \leq \; \frac{\Delta(X_n,X_{n+1})}{\Delta(X_n)} \;
\leq \; \Delta(X_{n+1}),
\end{equation}
see {\cite{Mitr} p.46).  Because of Proposition \ref{Prop4.3} we
concentrate on the case $h(n) \rightarrow \infty.\;$ The following
Lemma is needed.

\begin{Lemma} \label{Lemma4.4}
Let $(P_n)_{\nin}$ induce a polynomial hypergroup on ${\Bbb N}_0,$
and assume that $h(n) \rightarrow \infty.\;$ Then for every $k,l
\in {\Bbb N}_0$ holds $g(n,n+l,k) \rightarrow 0 \;$ as $n
\rightarrow \infty.$
\end{Lemma}

\begin{proof}
For $n\geq k$ we have with the Cauchy-Schwarz inequality
$$ g(n,n+l,k) \; = \; h(k) \; \int_{D_s} P_n(x) \; P_{n+l}(x) \; P_k(x) \; d\pi(x) $$
$$ \leq \; h(k) \; g(n,n,0)^{1/2} \; \left( \int_{D_s} (P_{n+l}(x) \; P_k(x))^2 \; d\pi(x) \right)^{1/2}. $$
Since $\; |P_{n+l}(x)P_k(x)| \; \leq 1 \; $ for $x \in D_s$ we see
that $g(n,n+l,k) \rightarrow 0$ as $n \rightarrow \infty.$
\end{proof}

\begin{Theorem} \label{4.5}
Let $(X_n)_{\nin}$ be a weakly stationary random sequence on the
polynomial hypergroup ${\Bbb N}_0$ induced by $(P_n)_{\nin}.\;$
Assume that $h(n) \rightarrow \infty$ as $n \rightarrow \infty.\;$
If $d(n) \; = E(X_n \overline{X_0}) $ tends to zero, the random
sequence $(X_n)_{\nin}$ is asymptotically $P_n$-deterministic.
\end{Theorem}

\begin{proof}
We have $ E(X_n \overline{X_n}) \; = \displaystyle\sum_{k=0}^{2n}
g(n,n,k) \; d(k).\;$ Since $\displaystyle\sum_{k=0}^{2n} g(n,n,k)
\; =1 \;$ and $g(n,n,k) \rightarrow 0$ as $n \rightarrow \infty$
by Lemma \ref{Lemma4.4}, Toeplitz's Lemma (see \cite{Kno}, p.377)
yields $\Delta(X_n) \; = E(X_n\overline{X_n}) \; \rightarrow 0,\;$
and by (\ref{27}) we have $\delta_n \rightarrow 0.$
\end{proof}

\begin{Remark} \label{Remark1.2}
Theorem \ref{4.5} improves Theorem 2 of \cite{hoe}, where we had
to assume that the recurrence coefficients $a_n, b_n$ and $c_n$
are convergent.

If the spectral measure $\mu$ is absolutely continuous with
respect to the measure $\pi$, we can deduce that $d(n) \rightarrow
0,\;$ see (\cite{Bloo}, Theorem 2.2.32(vi)).
\end{Remark}

We apply (\ref{27}) further to investigate also spectral measures
$\mu$ that contain discrete or singular parts.

\begin{Theorem} \label{4.6}
Let $(X_n)_{\nin}$ be a weakly stationary random sequence on the
polynomial hypergroup ${\Bbb N}_0$ induced by $(P_n)_{\nin},$
where $h(n) \rightarrow \infty$ as $n \rightarrow \infty.\;$
Assume that the spectral measure has the form
$$ \mu \; = \; f \pi \; + \mu_0 \; + \mu_1,$$
where $f \in L^1(D_s,\pi),\;\; \supp \mu_0 \subseteq D_{s,0} := \{
x \in D_s: \; P_n(x) \rightarrow 0$ for $n \rightarrow \infty \}
\;$ and $\supp \mu_1 \subseteq D_s \backslash D_{s,0} \;$ is a
finite set or empty. If $ m = \; |\supp \mu_1 |>0 \;$ and there is
some $n_0 \in {\Bbb N}_0$ such that
\begin{equation}
\inf_{n \geq n_0} \; \det \;\left( \begin{array}{ccc}
\langle P_n,P_n \rangle_{\mu_1} & \cdots& \langle P_n,P_{n+m-1} \rangle_{\mu_1} \\
\vdots && \vdots\\
\langle P_{n+m-1},P_n\rangle_{\mu_1} & \cdots & \langle
P_{n+m-1},P_{n+m-1} \rangle_{\mu_1}
\end{array} \right)
\quad > 0,
\end{equation}
then the random sequence $(X_n)_{\nin}$ is asymptotically
$P_n$-deterministic.
\end{Theorem}

\begin{proof}
If $\mu_1 =0$ we have $\Delta(X_n) \; = E(X_n \overline{X_n}) \; =
\displaystyle\int_{D_s} P_n^2(x) f(x) d\pi(x) \; +
\displaystyle\int_{D_s} P_n^2(x) d\mu_0(x).\;\;$ By Theorem
\ref{4.5} (and Lemma \ref{Lemma4.4} above) and the theorem of
dominated convergence we have $\Delta(X_n) \rightarrow 0$ and the
statement follows.

If $m = \; |\supp \mu_1|>0, \;$ we consider the Gramian
determinants $\Delta(X_n,...,X_{n+m-1})$ and
$\Delta(X_n,...,X_{n+m}).\;$ Expanding the determinants we obtain
$$ \Delta(X_n,...,X_{n+k}) \; = \; \det \;\left( \begin{array}{ccc}
\langle P_n,P_n \rangle_{\mu_1} &\cdots& \langle P_n,P_{n+k} \rangle_{\mu_1} \\
\vdots&&\vdots \\
\langle P_{n+k},P_n \rangle_{\mu_1} &\cdots & \langle
P_{n+k},P_{n+k}\rangle_{\mu_1}
\end{array} \right) \quad + \; \sum_{\sigma \in I_k} \varphi_\sigma,$$
where $\; |I_k|\;$ is finite and independent of $n$, whereas the
$\varphi_\sigma$'s are products of\\ $\langle P_{n+i_1},P_{n+j_1}
\rangle_{f \pi + \mu_0} $
and $\langle P_{n+i_2},P_{n+j_2} \rangle_{\mu_1}$ containing at least one factor\\
$\langle P_{n+i},P_{n+j} \rangle_{f \pi + \mu_0}.\;$ Since all $
\langle P_{n+i},P_{n+j} \rangle_{f\pi + \mu_0} $ tend to zero as
$n \rightarrow \infty,$ we conclude from the assumption that
$1/\Delta(X_n,...,X_{n+m-1})\;$ is bounded for $n \geq n_0.\;$
Since $\; |\supp \mu_1| \; =m \;$ the $(m+1) \times
(m+1)$-determinant
$$ \det \; \left( \begin{array}{ccc}
\langle P_n,P_n \rangle_{\mu_1}& \cdots& \langle P_n,P_{n+m} \rangle_{\mu_1} \\
\vdots&&\vdots\\
\langle P_{n+m},P_n \rangle_{\mu_1} &\cdots& \langle
P_{n+m},P_{n+m} \rangle_{\mu_1}
\end{array} \right) $$
has to be zero. Hence $\Delta(X_n,...,X_{n+m})\;$ tends to zero
with $n \rightarrow \infty.\;$ By (\ref{27}) we get $\delta_n
\rightarrow 0.$
\end{proof}

A typical situation to use Theorem \ref{4.6} is, when $\pi$ is
even and $\supp \mu_1 \subseteq D_s \backslash D_{s,0} =\{
-1,1\}.\;$ Since $P_n(1) =1 $ and $P_n(-1) \; =(-1)^n\;$ we obtain
for $\mu_1 = \alpha \epsilon_1 + \beta \epsilon_{-1}$
$$ \det \; \left( \begin{array}{cc}
\langle P_n,P_n\rangle_{\mu_1}& \langle P_n,P_{n+1} \rangle_{\mu_1} \\
\langle P_{n+1},P_n \rangle_{\mu_1} & \langle P_{n+1},P_{n+1}
\rangle_{\mu_1}
\end{array} \right) \;\; = \; 4 \alpha \beta.$$

The part of the spectral measure contained in $D_{s,0} = \{ x \in
D_s: \; P_n(x) \rightarrow 0 \} \;$ is easy to deal with, as we
have already seen in the proof of Theorem \ref{4.6}. Hence we want
to derive results on the size of $D_{s,0}.$

We restrict ourselves from now on to the case that $\pi$ is even,
i.e. $b_n=0$ for all $\nin.\,$ At first we consider the case of
$a_n \rightarrow a$ as $n \rightarrow \infty$, where $\frac{1}{2}
<a<1.\;$ The Cartier-Dunau polynomials $P_n(x)$ with $q \geq 2$
belong to this class.

We know then that $\supp \pi = [-\gamma,\gamma] \,$ is a proper
subset of $D_s = [-1,1]$ with \mbox{$\gamma = 2 \sqrt{a(1-a)}.\;$}
Note that $\gamma =1$ exactly when $a = \frac{1}{2}.\;$ In fact,
using theorems of Blumenthal and Poincar\'{e} this is proven in
Lasser (\cite{la02}, Theorem 2.2). In addition we can obtain the
following result. Compare also Voit (\cite{Voit2}, Theorem
8.2(4)).

\begin{Proposition} \label{Prop4.7}
Assume that $(P_n)_{\nin}$ is even and induces a polynomial
hypergroup on ${\Bbb N}_0.$ Furthermore let $a_n \rightarrow a$
with $\frac{1}{2} <a<1.\;$ Then $P_n(x) \rightarrow 0 $ as $n
\rightarrow \infty$ for each $x \in \; ]-1,1[\;.$
\end{Proposition}

\begin{proof}
Fix $\alpha_0 \in \; ]\gamma,1[\;. \;$ By the separating property
of zeros of $P_n(x),$ see \cite{Szeg}, we get $P_n(\alpha_0)
>0.\;$ Now $Q_n(x) \; = P_n(\alpha_0x) /P_n(\alpha_0) \,$ defines
a polynomial hypergroup on ${\Bbb N}_0.$ The recurrence relation
of the $Q_n(x)$ is given by
$$ x Q_n(x) \; = \; \tilde{a}_n Q_{n+1}(x) \; + \; \tilde{c}_n Q_{n-1}(x),$$
where
$$ \tilde{a}_n \; = \; a_n \; \frac{P_{n+1}(\alpha_0)}{P_n(\alpha_0) \, \alpha_0},
\qquad \tilde{c}_n \; = \; c_n \;
\frac{P_{n-1}(\alpha_0)}{P_n(\alpha_0) \,\alpha_0}.$$ By the
Poincar\'{e} theorem, compare the proof of Theorem 2.2 in
\cite{la02}, we have
$$ \lim_{n \rightarrow \infty} \frac{P_n(\alpha_0)}{P_{n-1}(\alpha_0)} \; <1.$$
Since $\tilde{a}_n$ is convergent, the dual space of the
polynomial hypergroup induced by $(Q_n)_{\nin}$ is $[-1,1].\;$ In
particular
$$ |P_n(\alpha_0x) /P_n(\alpha_0)| \; \leq 1$$
for all $x \in [-1,1],\;\; \nin.\;$ That means $\lim\limits_{n
\rightarrow \infty} \; |P_n(\alpha_0x) | \; \leq \lim\limits_{n
\rightarrow \infty} P_n(\alpha_0) \; =0 \;$ for all $x \in
[-1,1].\;$ Since $\alpha_0 \in\;]\gamma,1[ \,$ can be chosen
arbitrarily we have $\lim\limits_{n \rightarrow \infty} P_n(x)
=0\;$ for all $ x \in \; ]-1,1[\;.\;$
\end{proof}

\begin{Corollary} \label{Cor4.8}
Let $(X_n)_{\nin}$ be a weakly stationary random sequence on the
polynomial hypergroup ${\Bbb N}_0$ induced by $(P_n)_{\nin}.\;$
Suppose that $\pi$ is even and $a_n \rightarrow a$ with
$\frac{1}{2} <a<1.\;$ Then $(X_n)_{\nin}$ is asymptotically
$P_n$-deterministic.
\end{Corollary}

\begin{proof}
By Proposition \ref{Prop4.3} we get $\lim\limits_{n \rightarrow
\infty} h(n) \;=\infty\;$ and the statement follows from
Proposition \ref{Prop4.7}.
\end{proof}

If $a_n \rightarrow \frac{1}{2}$ the problem is much more
involved. A sufficient condition for $P_n(x) \rightarrow 0$ for
every $x \in \;]-1,1[\;\;$ can be derived by making use of the
Turan determinant. Given the orthogonal polynomials $(P_n)_{\nin}$
fulfilling (\ref{threeTerm}), the Turan determinant is defined by
\begin{equation}
\theta_n(x) \;=\; h(n) \; \left( P_n^2(x) \;-\;\frac{a_n}{a_{n-1}}
\;P_{n-1}(x) \; P_{n+1}(x) \right)
\end{equation}

\begin{Proposition}\label{Prop4.9}
Assume that $(P_n)_{\nin}$ is even and induces a polynomial
hypergroup on ${\Bbb N}_0.\;$ Further let $a_n \rightarrow
\frac{1}{2}$ as $ n \rightarrow \infty.\;$ The following
inequalities for $\theta_n(x)$ are valid:
\begin{enumerate}
\item[(i)] $ \theta_n(x) /h(n) \;\leq C_1\;\left( P_{n-1}^2(x) +
P_n^2(x) + P_{n+1}^2(x) \right) \;\;$ for all $x \in {\Bbb P},\;
\nin,\;$ where $C_1>0$ is a constant independent of $x$ and $n$.
\item[(ii)] $\;|\theta_n(x)-\theta_{n-1}(x) |\, /h(n) \;\leq C_2
\;|a_{n-1} c_n -a_{n-2}c_{n-1} | \; \left( P_{n-1}^2(x) + P_n^2(x)
\right) \;\;$ for all $x \in [-1,1], \;\nin,\;$ where $C_2 >0$ is
a constant independent of $x$ and $n$. \item[(iii)] Given some
$\delta \in \;]0,1[\;\;$ there is $N \in {\Bbb N}$ such that
\\ $\theta_n(x) /h(n) \; \geq C_3 \;\left( P_{n-1}^2(x) + P_n^2(x) \right) \;   \;$
for all $x \in [-1+\delta,1-\delta],\;n\geq N,\;\;$ where $C_3>0$
is a constant independent of $x$ and $n$.
\end{enumerate}
\end{Proposition}

\begin{proof}
\begin{enumerate}
\item[(i)] Since $\lim\limits_{n \rightarrow \infty}
\displaystyle\frac{a_n}{a_{n-1}} \;=1 \;$ and $ \;2 \;|P_{n-1}(x)
P_{n+1}(x)| \;\leq P_{n-1}^2(x) + P_{n+1}^2(x),\;$ the inequality
of (i) follows immediately. \item[(ii)] By using the recurrence
relation (\ref{threeTerm}) we get
\begin{equation} \label{30}
\theta_n(x) \;=\; h(n) \; P_n^2(x) \;+\; h(n-1) \,P_{n-1}^2(x)
\;-\; \frac{x}{a_{n-1}} \; h(n) \,P_{n-1}(x) P_n(x)
\end{equation}
and
\begin{equation} \label{31}
\theta_n(x) \;=\; h(n) \,P_n^2(x) \;+\; h(n+1)\;
\frac{a_nc_{n+1}}{a_{n-1} c_n} \; P_{n+1}^2(x) \;- h(n)
\;\frac{a_n x}{a_{n-1} c_n} \; P_n(x) P_{n+1}(x).
\end{equation}
Employing (\ref{30}) to $\theta_n$ and (\ref{31}) to $\theta_{n-1}
$ gives
$$ \theta_n(x) - \theta_{n-1}(x) \;=\; h(n) \; \left( 1\;-\;
\frac{a_{n-1}c_n}{a_{n-2} c_{n-1}} \right) \;P_n^2(x) $$
$$+\; \left( h(n-1) \; \frac{a_{n-1}}{a_{n-2} c_{n-1}} \;-\; h(n) \;
\frac{1}{a_{n-1}} \right) \; x \,P_{n-1}(x) P_n(x) $$
$$ = \; h(n) \;\left( \frac{a_{n-2}c_{n-1} - a_{n-1}c_n}{a_{n-2}c_{n-1}} \;P_n^2(x) \;+\;
\frac{a_{n-1}c_n -c_{n-1} a_{n-2}}{c_{n-1}a_{n-2}a_{n-1}} \; x \,
P_{n-1}(x) P_n(x) \right). $$ Since $a_n \rightarrow \frac{1}{2},
\; c_n \rightarrow \frac{1}{2}\;\;$ we obtain for $|x|\leq 1 $
$$|\theta_n(x) -\theta_{n-1}(x)|\, /h(n) \;\leq \; C_2 \;|a_{n-1}c_n \;-\; a_{n-2}c_{n-1} | \;
\left( P_{n-1}^2(x) \;+\; P_n^2(x) \right).$$ \item[(iii)]
Applying (\ref{30}) we have
$$ \theta_n(x) \;=\; h(n) \;\left( P_n(x) \;-\; \frac{x}{2a_{n-1}} \;
P_{n-1}(x) \right)^2 \;+\; \left( 1 \;-\; \frac{x^2}{4c_na_{n-1}}
\right) \; h(n-1) \;P_{n-1}^2(x)$$ and
$$ \theta_n(x) \;=\;h(n-1) \;\left( P_{n-1}(x) \;-\;\frac{x}{2c_n} \;P_n(x)
\right)^2 \;+\;\left( 1 \;-\; \frac{x^2}{4 c_na_{n-1}} \right) \;
h(n) \;P_n^2(x).$$ In particular it follows
$$ \theta_n(x) \;\geq \; \left( 1 \;-\; \frac{x^2}{4c_na_{n-1}} \right) \; h(n-1) \;P_{n-1}^2(x) $$
and $$ \theta_n(x) \;\geq \; \left( 1 \;-\;
\frac{x^2}{4c_na_{n-1}} \right) \;h(n) \;P_n^2(x).$$ Since $x \in
[-1+\delta,1-\delta] \;$ there is a constant $C_3$ and $N \in
{\Bbb N}$ such that
$$ \theta_n(x) /h(n) \;\geq \;C_3 \; \left( P_{n-1}^2(x) \;+\; P_n^2(x)\right)
\qquad \mbox{ for all } n \geq N.$$
\end{enumerate}
\end{proof}

\begin{Theorem} \label{4.10}
Assume that $(P_n)_{\nin}$ is even and induces a polynomial
hypergroup on ${\Bbb N}_0.\;$ Let $a_n \rightarrow \frac{1}{2}$
and $h(n) \rightarrow \infty$ as $n \rightarrow \infty.\;$ Further
suppose that
\begin{equation}
\sum_{n=1}^\infty |a_nc_{n+1} \;-\; a_{n-1}c_n| \;< \infty.
\end{equation}
Then every weakly stationary random sequence $(X_n)_{\nin}$ on the
polynomial hypergroup ${\Bbb N_0}$ induced by $(P_n)_{\nin}$ is
asymptotically $P_n$-deterministic.
\end{Theorem}

\begin{proof}
It suffices to prove that $P_n(x) \rightarrow 0$ as $n \rightarrow
\infty \;$ for every $x \in \;]-1,1[\;.\;$ Proposition
\ref{Prop4.9}(ii) and (iii) imply
$$ |\theta_n(x) -\theta_{n-1}(x)| \;\leq \; h(n) \;C_2\;|a_{n-1}c_n \;-\; a_{n-2}c_{n-1} | \;
\left( P_{n-1}^2(x) \;+\; P_n^2(x) \right) $$
$$ \leq \; \frac{C_2}{C_3} \;| a_{n-1}c_n \;-\; a_{n-2}c_{n-1} | \; \theta_n(x) \;
=\; \epsilon_n \theta_n(x) $$ for all $n \geq N,\;$ where
$\epsilon_n \;= \displaystyle\frac{C_2}{C_3}\;|a_{n-1}c_n
\;-a_{n-2}c_{n-1}|.\;$ Hence
$$\frac{1}{1+\epsilon_n} \; \theta_{n-1}(x) \; \leq \; \theta_n(x) \;
\leq \frac{1}{1-\epsilon_n} \; \theta_{n-1}(x)$$ for every $n \geq
N.\;$ Since $\displaystyle\sum_{n=2}^\infty \epsilon_n \;$ is
convergent, $\theta_n(x)$ is convergent, too. Applying once more
Proposition \ref{Prop4.9}(iii) and $h(n) \rightarrow \infty,$ we
get $P_n(x) \rightarrow 0.$
\end{proof}

We conclude this section with some further examples.
\begin{enumerate}
\item[(d)] The Cartier-Dunau polynomials that are essential for
the investigation of stationary radial stochastic processes on
homogeneous trees have the property that $a_n =
\displaystyle\frac{q}{q+1}.\;$ Hence for $q \geq 2$ we can use
Corollary \ref{Cor4.8} and get that $\delta_n \rightarrow 0$ for
the corresponding random sequences $(X_n)_{\nin}$. \item[(e)]
$P_n(x; \beta|q) \;$ the (continuous) $q$-ultraspherical
polynomials with $-1 < \beta <1 \;$ and $0<q<1.\;$ Their
hypergroup structure is studied in \cite{la01}, compare also
Bressoud \cite{Bres}. The recurrence coefficients $a_n,c_n$ are
not given in explicit form, but we know the asymptotic behavior,
see (\cite{Bloo}, p.168). In fact
$$\alpha_n = \frac{1}{2}+\frac{1}{n}+o(\frac{1}{n})$$ and
hence $a_n \rightarrow \frac{1}{2}$ as $n \rightarrow \infty$. It
now follows with elementary computations from (\ref{HaarProd})
that $h_n \rightarrow \infty$ with growing n.

To check properties of Theorem \ref{4.10} we investigate the
monotonicity of the sequence $(\alpha_n)_\nin.\;$ A direct
calculation shows that $(\alpha_n)_\nin$ is increasing if $ \beta
\leq q,\;$ and decreasing if $q \leq \beta.\;$ In both cases we
have $\sum\limits_{n=1}^\infty |\alpha_{n+1} - \alpha_n | \; = \;
|\frac{1}{4} - \alpha_1|.\;$ Now Theorem \ref{4.10} implies that
every corresponding random sequence is asymptotically
$P_n(\cdot\;\;;\beta|q)$-deterministic.

\end{enumerate}

\section{Prediction with additional information}
Moving average processes play a prominent role in the classical
time series theory and have a lot of applications. Their
comparatively simple structure allows efficient algorithms for
estimating required parameters (for example \cite{Bro}).

We recall from equation (\ref{MovAverage}) that MA(q)-processes
with respect to polynomial hypergroups have the form
$$X_n  = \sum_{k=0}^q a_k {\bf T_n} Z_k h(k),$$
where $(Z_k)_{k \in \nz_0}$ is a white noise sequence with respect
to the polynomial system and ${\bf T_n}$ the associated shift
operators. The spectral measure $\mu$ is, like in the classical
case, completely determined by the coefficients $a_k \in \cz, \:
k=0,...,q$:
$$\mu = |\hat a|^2 \pi= | \sum_{k=0}^q a_k P_k h(k)|^2 \pi.$$
Therefore, the prediction of such processes is more easy: No
estimation of $\mu$ has to be done. In special cases, the
orthogonal polynomial system to $\mu$, which is required for the
optimal linear predictor, is explicitly known:

\begin{Example}

Let $P_n^{(\alpha, \alpha)} (x)$ be the ultrasherical polynomials
with orthogonalization measure $\pi_{\alpha,\alpha}$ and
\begin{eqnarray*}
X_n& =& {\bf T_n} Z_0 + {\bf T_n} Z_1\\
 &=&g(n,1,n-1) Z_{n-1}
+ \Bigl(1 + g(n,1,n)\Bigr) Z_n + g(n,1,n+1) Z_{n+1}.
\end{eqnarray*}
the corresponding MA(1) process. Then, the spectral measure is
$$d\mu(x) = (1+x)^2 d\pi_{\alpha,\alpha}(x),$$
and the Jacobian polynomials $P_n^{(\alpha,\alpha+2)}$ the
associated orthogonal polynomials.

\end{Example}

For MA(q) processes with respect to a polynomial system we
introduce an alternative concept:

Assume, that besides the realizations of $X_1,\ldots,X_n$ the
first $q$ initial realizations of the associated white noise
$Z_0,...,Z_{q-1}$ are also observable. Then, the optimal linear
predictor has a representation which does not resort to the
spectral measure.

In the special case of $MA(1)$ processes with respect to Chebyshev
polynomials of the first kind we give an explicit representation
of the predictor involving Chebyshev polynomials of the second
kind. We are not aware if this pleasing result has an analogy in
other polynomial systems.

To make formulas more readable, we subsume $h(k)$ in $a(k)$ and
write until the end of the chapter, somewhat carelessly, for $a(k)
h(k)$ only $a(k)$.

Now, let $(X_n)_{n \in \nz_0}$ be a $MA(q)$-process with respect
to $(P_n)_{n \in \nz_0}$ and
$$\tilde H_N:=<Z_q,...,Z_{q+N}>,\quad N=0,1,2,...,\;$$
the Hilbert space generated from the white noise sequence
$Z_q,...,Z_{q+N}$. $\tilde H_N$ is a sub-Hilbert space of $\tilde
H:=\overline {<(Z_n)_{n \in \nz_0}>}$. Denote as $\tilde X_N$ the
orthogonal projection from $X_N$ to $\tilde H_N$.

\begin{Remark} \label{bem331}
The representation
\[X_n = \sum_{k=0}^q a_k \sum_{s=|k-n|}^{n+k} g(n,k,s)\; Z_s =
  \sum_{s=0}^{q+n} ( \sum_{k=0}^q a_k \: g(n,k,s))\; Z_s\]
induces with $a_q =1$ for $n \in \nz_0$:

\begin{eqnarray*}
\tilde X_n &=& \sum_{s=q}^{n+q} ( \sum_{k=0}^q a_k \: g(n,k,s)) Z_s\\
& = & g(n,q,n\!+\!q) Z_{n+q} + [ g(n,q,n\!+\!q\!-\!1) + g(n,q\!-\!1,n\!+\!q\!-\!1) a_{q-1} ]\\
&& \cdot  Z_{n+q-1}+ \ldots + \sum_{k=0}^q a_k \: g(n,k,q)\; Z_q.\\
\end{eqnarray*}
As $g(n,q,n\!+\!q) \neq 0$ is always satisfied, the recursion \\[0.2cm]

\hspace*{3.2 cm}$Z_q = \tilde X_0$\\
and for $n=1, \ldots , N:$
\begin{equation} \label{gl331}
 Z_{q+n} = \frac{1}{g(n,q,n\!+\!q)} ( \tilde X_n - \sum_{s=q}^{q+n-1}
   ( \sum_{k=0}^q a_k \: g(n,k,s))\; Z_s)
\end{equation}
is defined.

Thus, $Z_{q+n}$ has the representation:\\[0.2cm]
\[ Z_{q+n} = \sum_{r=0}^n b_{n r} \tilde X_r
 = \sum_{r=0}^n b_{nr} (X_r- \sum_{s=0}^{q-1} (\sum_{k=0}^q a_k \: g(r,k,s))\; Z_s).
\]
This shows especially that $ \tilde H_N = < \tilde X_0,...,\tilde
X_N> \subseteq <X_0,...,X_N, Z_0,...,Z_{q-1}>$ is valid.
\end{Remark}

\begin{Definition}
The orthogonal projection $\stackrel{\approx}{X}_{N+1}$ of
$X_{N+1}$ to \mbox{$\stackrel{\approx}{H}_N := <
X_0,...,X_N,Z_0,...,Z_{q-1}>$} is called \bf{predictor of
$X_{N+1}$ given the knowledge $Z_0,...,Z_{q-1}$}.
\end{Definition}

\begin{Theorem}
Let $\displaystyle X_n = \sum_{k=0}^q a_k {\bf T_n} Z_k,$ be a
$MA(q)$-process with respect to $(P_n)_{n \in \nz_0}$ and assume
$a_q=1$. Then it is true that:
{
\begin{enumerate}
\item $\stackrel{\approx}{X}_{N+1}$ has the representation:
\[
  \stackrel{\approx}{X}_{N+1} = {g(N\!+\!1,q,N\!+\!q\!+\!1)} ( \sum_{s=0}^{q-1}
   \tilde b_{N+1 s} \;Z_s - \sum_{r=0}^N b_{N+1 \: r}\; X_r )\]
with \quad $\displaystyle \tilde b_{N\!+\!1 s} := \sum_{r=0}^{N+1}
\sum_{k=0}^q b_{N+1 \: r} \: a_k \: g(r,k,s),$ where the
coefficients \mbox{$b_{N\!+\!1 \: r}, \; r=0,...,N\!+\!1$}, are
given by the recursion (\ref{gl331}) for $Z_{q+N+1}$;.\\[-0.3cm]
\item The prediction error is:
\[ || X_{N+1} - \stackrel{\approx}{X}_{N+1} ||_2
  = \frac{g(N\!+\!1,q,N\!+\!q\!+\!1)}{\sqrt{h(q\!+\!N\!+\!1)}}.\]
\end{enumerate}
}
\end{Theorem}

\begin{Proof}
Remark \ref{bem331} shows:
\begin{eqnarray}
 Z_{q+N+1} &= & \sum_{r=0}^{N+1} b_{N+1 \: r} \tilde X_r\\
         &=& \sum_{r=0}^{N+1} b_{N+1 \: r} X_r - \sum_{r=0}^{N+1} b_{N+1 \:r}
   \sum_{s=0}^{q-1} ( \sum_{k=0}^q a_k \: g(r,k,s))\; Z_s.
\end{eqnarray}
This equation can be solved for $X_{N+1}$:
\begin{eqnarray} \label{gl333}
 X_{N+1} & = & {g(N\!+\!1,q,N\!+\!q\!+\!1)} ( Z_{q+N+1} - \sum_{r=0}^N b_{N+1 \: r}
   X_r \nonumber \\
   && + \sum_{s=0}^{q-1} ( \sum_{r=0}^{N+1}  \sum_{k=0}^q b_{N+1 \: r}
   a_k \: g(r,k,s))\; Z_s).
\end{eqnarray}

Since
\[ ( Z_{q+N+1}, Y) = 0 \mbox{ \quad f\"{u}r } Y \in
  \stackrel{\approx}{H}_N\]
is true, the orthogonal projection from $X_{N+1}$ to $
\stackrel{\approx}{H}_N$ is given by leaving out
$${g(N\!+\!1,q,N\!+\!q\!+\!1)} \; Z_{q+N+1}$$
in the representation (\ref{gl333}). Therefore the prediction
error is:
\begin{eqnarray*}
 || X_{N+1} -  \stackrel{\approx}{X}_{N+1} ||_2 & =& || {g(N\!+\!1,q,N\!+\!q\!+\!1)} \; {Z_{q+N+1}}||_2\\[0.2cm]
   &=&\frac{g(N\!+\!1,q,N\!+\!q\!+\!1)}{\sqrt{h(q\!+\!N\!+\!1)}}.
\end{eqnarray*}
\end{Proof}

\begin{Remark}
If there is a $C>0$ with $ C\leq g(N\!+\!1,q,N\!+\!q\!+\!1),$ the
prediction error tends with growing $N$ towards zero if and only
if the Haar measure $h(N)$ growth unboundedly.
The order of the convergence is $\sqrt{h(N)}$. \\
\end{Remark}
\newpage
\begin{Theorem}
Let $(X_n)_{n \in \nz_0}$ be a $MA(1)$-process with respect to
Chebyshev polynomials of the first kind:\\[-0.2cm]
\[ X_n = \sum_{k=0}^1 a_k {\bf T_n} Z_k \mbox{\qquad and $a_1:=1, \;a_0:=a$}\]
and $U_n(x)$ the Chebyshev polynomials of the second kind. Then,
the predictor of $X_{N+1}$ given the knowledge $Z_0$
is:\\[-0.2cm]
\[ \stackrel{\approx}{X}_{n+1} = \frac{1}{2} ([ a U_{n+1} (-a) + U_n (-a) ] Z_0 - 2 \sum_{t=0}^n
   U_{n+1-t} (-a) X_t + U_{n+1} (-a) X_0 ).\]
\end{Theorem}

\begin{proof}
With \hspace{3cm} $X_0=a Z_0 +Z_1$;
\[ X_n = a Z_n + \frac{1}{2} Z_{n+1} + \frac{1}{2} Z_{n-1} \quad
\mbox{ for } n\geq 1\] one has \quad
$\tilde X_0 = Z_1, \quad \tilde X_1 = a Z_1 + \frac{1}{2} Z_2$ \quad and\\[0.1cm]
\[
 \tilde X_n = X_n = a Z_n + \frac{1}{2} Z_{n+1}+ \frac{1}{2} Z_{n-1}
 \qquad \mbox{ f\"{u}r $n \geq 2$}.\]

The matrix form of the system is:\\
\[
 \left( \begin{array}{lllll}
           1\\
           a           &   \frac{1}{2} \\
           \frac{1}{2} &       a       & \frac{1}{2}  \\
           0           &  \frac{1}{2}  &      a       &  \frac{1}{2}\\
           \vdots \ddots &   \ddots    &  \ddots      &  \ddots     \\
           0 \cdots     &      0       & \frac{1}{2}  &     a  &\frac{1}{2}\\
        \end{array}
 \right)
 \left(
 \begin{array}{c}
   Z_1\\
   \vdots\\
   \vdots\\
   \vdots\\
   Z_{n+1}\\
  \end{array}
  \right)
  =
  \left(
  \begin{array}{c}
   \tilde X_0\\
   \vdots\\
   \vdots\\
   \vdots\\
   \tilde X_n\\
  \end{array}
  \right).
\]
\vspace{0.3cm}

The first column of the inverse of this triangular matrix has to
satisfy the equations
\[c_0 =1 ; \quad a c_0 + \frac{1}{2} c_1 =0;\]
\[\frac{1}{2} c_{n-2} + a c_{n-1} + \frac{1}{2} c_n =0, \quad n\geq 2 .\]
Comparison with the three term recursion of Chebyshev polynomials
of the second kind (for example: \cite{la01}, 3(a))
\[ U_{-1} (x)=0 , \qquad U_0 (x)=1,\]
\[2 x U_{n-1}(x)= U_n (x) + U_{n-2} (x) \mbox{\qquad for $ n=1,2, \ldots$}\]
shows:\\[0.3cm]
 $  c_k =  U_k (-a)\; $ is solution of the equations for the first column of the inverse matrix.\\[0.2cm]
The 2nd to $(n+1)st$-columns yield the same equations besides the
initial value $c_0=2$.
In this case the solutions are:\quad $c_k=2 U_k (-a)$.\\[0.4cm]
So, we get as inverse matrix
\[
 \left( \begin{array}{ccc}
  U_0 (-a) \\
  U_1 (-a) & 2 U_0 (-a) \\
  \vdots   & \vdots & \ddots\\
  U_n(-a)  & 2 U_{n-1} (-a) \cdots & 2 U_0 (-a)\\
  \end{array}
 \right),
\]
and further
\[
\stackrel{\approx}{X}_{n+1} = {g(n\!+\!1,1,n\!+\!2)} ( \tilde
b_{n+1 \: 0} \; Z_0 - \sum_{t=0}^n
    b_{n+1 \: t} \; X_t) \]
with
\begin{eqnarray*}
b_{n+1 \:0} & = & U_{n+1} (-a); \\[0.2cm]
b_{n+1 \:k} & = & 2 U_{n+1-k} (-a), \qquad k=1, \ldots, n+1;\\[0.2cm]
\tilde b_{n+1 \: 0} & = & \sum_{t=0}^{n+1} \sum_{k=0}^1 b_{n+1 \:
t}
a_k \: g(t,k,0)\\[0.2cm]
& = & a U_{n+1} (-a) + U_n (-a).
\end{eqnarray*}
This finishes the proof.
\end{proof}
\newpage

\section{A Levinson type algorithm}
In this last section we present an algorithm which is suited to
calculate the coefficients of the best linear predictor for
$P_n$-weakly stationary sequences, given before. The algorithm
parallels the Levinson algorithm used in the prediction theory of
classical weakly stationary processes \cite{Del}. It is a fast
inversion algorithms working for a large class of symmetric
positive definite matrices. These matrices are given by orthogonal
polynomials not necessarily inducing a hypergroup. For the
applications we have in mind, the matrices are the covariance
matrices of $P_n$-weakly stationary processes.

\subsection{Coefficients of the best linear one-step predictors}

From Theorem  \ref{s311} we know that the coefficients $b_{n k}$
of the best linear one-step predictor
\[
\hat X_{n+1} = \sum_{k=0}^n b_{n k} X_k \] are determined by the
equation
\[Q_{n+1}^* (x) = P_{n+1}(x) - \sum_{k=0}^n b_{nk} P_k(x)\]
with
\[
  Q_{n+1}^* (x) = \frac{\sigma_{n+1}}{\varrho_{n+1}} q_{n+1}(x);\]

where $q_n, n=0,1, \ldots$ are the orthonormal polynomials with
respect to the spectral measure $\mu$ and have leading
coefficients $\varrho_{n+1}=\varrho_{n+1}(\mu)$. As before
$\sigma_{n}=\sigma_{n}(\pi)$ are the leading coefficients of the
hypergroup polynomials $P_n$.

The algorithm, which we present in the following, provides from
the covariance matrix
$$\mathbf{M}_{n+1}=(d(k,l))_{0 \leq k,l \leq n+1}$$
the connection coefficients $c_{\phi P}(k,l),\quad 0 \leq l \leq k
\leq n+1$ between $P_k$ and the monic polynomials $\phi_l$ being
orthogonal with respect to the spectral measure $\mu$.

Assuming that we know the leading coefficients of the $P_n$,
equation
\[
 Q^*_{n+1} (x)= \sigma_{n+1} (\pi) \phi_{n+1}(x)= P_{n+1}
  (x) + \sigma_{n+1} (\pi) \sum_{k=0}^n c_{\phi P} (n+1,k) P_k (x)\]
 shows how we get the coefficients of the predictor:
  \[
   \hat X_{n+1} = - \sigma_{n+1} (\pi) \sum_{k=0}^n c_{\phi P} (n+1,k) X_k.\]

We already mentioned that the algorithm does not require that
$P_n$ defines a polynomial hypergroup. The strategy of fitting
together a modified Chebyshev algorithm \cite{Gau} and an
algorithm presented in \cite{AsII} works for a large class of
matrices.

Techniques which allow efficient triangular decompositions and
inversions of structured matrices are important for many
applications and have been studied for a long time (for example
\cite{Olshevsky}, \cite{Schur}). Involved algorithms are usually
considered as fast if they do the job for n-th order matrices in
$O(n^2)$ steps or less. So the classical Cholesky decomposition of
general positive definite symmetric matrices is not fast as it
requires $O(n^3)$ steps \cite{Golup}. Fast algorithms for non
sparse matrices exist for certain structured types. These types
include Toeplitz and Hankel matrices and more general matrices
with a displacement structure \cite{Pan}, \cite{Kailath},
\cite{Heinig}.

\subsection{$(P_s,\mu)$-structured matrices}

We start from a fixed system of real polynomials $(P_n)_{\nin }$
with $deg(P_n)=n$, having positive leading coefficients and being
orthogonal with respect to a positive Borel measures $\pi$ on the
real axis with infinite support. Such polynomials satisfy a
three-term recursion
$$x P_k (x)= a_k P_{k+1}(x)+ b_k P_k(x)+c_k P_{k-1}(x) \qquad k=0,1,\ldots$$
$$P_{-1}(x)=0, \qquad P_0(x)=P_0 >0$$
with recursion coefficients $a_k>0$, $c_k>0$ and $b_k \in
\mathbb{R}$  \cite{AsI}. We assume that the recursion coefficients
as well as the linearization coefficients $g(k,l,s)$ in the
equation
$$P_k(x) P_l(x) = \sum_{s=|k-l|}^{k+l} g(k,l,s) P_s(x)$$
are explicitly known. For a lot of polynomial systems such
explicit representations can be found in the literature (for
example \cite{Szeg}, \cite{Ismail}).

Next, we choose another Borel measure $\mu$ on $\mathbb{R}$ with
infinite support and finite moments $\mu_k=\int_\mathbb{R}x^k
d\mu(x)$ of all orders $k=0,1,\ldots$.

With these settings we call an $(n+1)\times(n+1)$-matrix
$\mathbf{M_n}$ \emph{$(P_s,\mu)$-structured} if it has the form

$$\left(
\begin{array}{cccccc}
  \int_\mathbb{R} P_0(x) P_0(x) d\mu(x) &  & \cdots & \int_\mathbb{R} P_0(x) P_n(x) d\mu(x)
  \\\\
  \vdots &  & \ddots  & \vdots \\\\
  \int_\mathbb{R} P_n(x) P_0(x) d\mu(x) & & \cdots &  \int_\mathbb{R} P_n(x) P_n(x) d\mu(x)
\end{array}
\right).$$

We denote the entries of $\mathbf{M_n}$ as $d(k,l)$ and
 set $d(k):= \int_\mathbb{R} P_k(x) d\mu(x)$.

 As the orthogonal polynomials $\{P_k\}_{k=0}^{n}$ and the defining measure $\mu$
can be chosen arbitrarily except for the stated restrictions, this
definition allows for a large amount of structured matrices.

A simple example is provided by the Chebyshev polynomials of the
first kind $T_k(x)$. The recursion formula of the polynomials
implies
$$d(k,l)=1/2 (d(|k-l|)+ d(k+l)),$$
showing that $( T_s,\mu)$-structured matrices have a specific
Toeplitz-plus-Hankel structure.

In the following we will present an algorithm which yields a
decomposition of the inverse
$$\mathbf{M}_n^{-1}= \mathbf{L}_n^{T} \mathbf{D}_n \mathbf{L}_n $$
into a product of a lower triangular matrix $\mathbf{L}_n$, its
transposed $\mathbf{L}_n^{T}$ and a diagonal matrix with positive
entries. The calculation required is $O(n^2)$ and therefore linear
equations
$$\mathbf{M_n}\mathbf{x_{n+1}}=\mathbf{b_{n+1}}$$
with given $(n+1)$-vector $\mathbf{b_{n+1}}$ can be solved fast
for $\mathbf{x_{n+1}}$.

Let  $\phi_k$  be the monic polynomials associated  to the measure
$\mu$ and let $\alpha_k$, $\beta_k$  be the recursion coefficients
in the three term recursion
$$x \phi_k= \phi_{k+1}+\alpha_k \phi_{k}+\beta_k \phi_{k-1} \qquad k=0,1,\ldots$$
  $$\phi_{-1}(x)=0, \qquad \phi_0(x)= 1.$$
Further, let $c_{\phi P}(k,l)$ be the
connection coefficients  of the expansion
$$\phi_k(x)=\sum_{l=0}^k c_{\phi P}(k,l)  P_l(x)$$
and
$$\sigma_{k,l}=\int_{\mathbb{R}} \phi_k(x) P_l(x)d\mu(x)$$
the mixed moments. Note, that one defines $c_{\phi P}(k,l)=0$ for
$k<l$ and that $\sigma(k,k)$ and $c_{\phi P}(k,k)$ are not zero
for $k=0,1,\ldots$.

\begin{Proposition}
 The inverse $\mathbf{M}_n^{-1}$ of a $(P_s,\mu)$-structured
 matrix $\mathbf{M}_n$ always exists and has a decomposition
 $$\mathbf{M_n}^{-1}= \mathbf{L}_n^{T}\mathbf{D}_n \mathbf{L}_n $$
 with
$$\mathbf{L}_n=
\left(
\begin{array}{cccccc}
   c_{\phi P}(0,0) & 0            & \cdots &  0  \\
   c_{\phi P}(1,0) & c_{\phi P}(1,1)&  \ddots &\vdots   \\
  \vdots & &\ddots &   0 \\
c_{\phi P}(n,0) & \cdots&\cdots  &  c_{\phi P}(n,n)
\end{array}
\right)
$$

and a diagonal matrix $\mathbf{D}_n$ with k-th diagonal element
$1/(\sigma_{k,k}\cdot c_{\phi P}(k,k))$.
\end{Proposition}

\begin{proof}
The inverse $\mathbf{M}_n^{-1}$ exists since the columns of
$\mathbf{M}_n$ are linear independent: $\sum_{i=0}^n \lambda_i
d(s,i)=0$ for all $s=0,\ldots,n$ means that for those $s$
$$\int_\mathbb{R} P_s(x)Q(x) d\mu(x)=0$$
with the polynomial $Q(x)= \sum_{i=0}^n \lambda_i P_i(x)$ is
valid. But this is only possible if $Q(x)$ is the zero polynomial
and $\lambda_i=0$ for $i=0,\ldots,n$.

For the decomposition first observe that
$$\int_\mathbb{R} (\phi_k(x))^2 d\mu(x)=\int_\mathbb{R} \phi_k(x)\sum_{l=0}^k c_{\phi P}(k,l)  P_l(x) d\mu(x)=
 c_{\phi P}(k,k) \sigma_{k,k}.$$
The last equality holds since the orthogonality of the $\phi_k(x)$
with respect to $\mu$ implies $\sigma_{k,l}=0$ for $l<k$.

Now, again using orthogonality one gets for $0 \leq k,l \leq n$
$$  c_{\phi P}(k,k) \sigma_{k,k}\cdot \delta_{k,l} =\int_\mathbb{R} \phi_k(x)\phi_l(x)
d\mu(x)= \sum_{s=0}^k \sum_{t=0}^l c_{\phi P}(k,s)c_{\phi
P}(l,t)d(s,t).$$

With the convention $c_{\phi P}(k,l)=0$ for $k<l$ this is just the
matrix equation
$$\mathbf{D}_n^{-1}= \mathbf{L}_n \mathbf{M_n} \mathbf{L}_n^T $$
which implies $\mathbf{D}_n= (\mathbf{L}_n^T)^{-1}
\mathbf{M}_n^{-1} \mathbf{L}_n^{-1} $ and thus
$$\mathbf{M}_n^{-1}= \mathbf{L}_n^{T}\mathbf{D}_n \mathbf{L}_n .$$

\end{proof}

The task is now to calculate the coefficients $c_{\phi P}(k,l)$
and $\sigma_{k,l}$ for $0 \leq k,l \leq n$ from the matrix
$\mathbf{M}_n$. Actually our algorithm needs the moments $d(k)$
for $k=0,\ldots,2n-1$. $\mathbf{M}_n$ shows $d(k)$ only up to
$k=n$. But this does not spoil anything:

\begin{Proposition}\label{moments}
The $(n+1)\times(n+1)$ entries of a $(P_s,\mu)$-structured matrix
are uniquely given by the moments  $d(k), k=0,\ldots,2n$. Vice
versa these $(2n+1)$ moments can be calculated from the matrix in
$O(n^2)$ steps.
\end{Proposition}

\begin{proof}
The linearization of $P_k(x)P_l(x)$ yields
$$d(k,l)= \sum_{s=|k-l|}^{k+l}g(k,l,s)d(s),$$
showing that $d(k,l), 0 \leq k,l \leq n$ can be calculated from
$d(k), 0 \leq k \leq 2n$.

Now, let the matrix and thus $d(k,l), 0 \leq k,l \leq n$ be given.
Then $d(k,0)=P_0 d(k)$ is already at hand for $k=0,\ldots,n$.
Since $g(k,l,k+l)\neq 0$ is always valid, one obtains the missing
moments $d(n+k), 1 \leq k \leq n$ successively from
$$d(n+k)=\frac{1}{g(n,k,n+k)}\left(d(n,k)-\sum_{s=n-k}^{n+k-1}g(n,k,s)d(s)\right).$$
The required calculation is obviously $O(n^2)$.
\end{proof}

\begin{Remark}
With respect to this observation one could also have the viewpoint
that $d(k), k=0,\ldots,2n$ are given beforehand and $\mathbf{M}_n$
is derived. If the $(P_n)_{\nin }$ define a hypergroup on
$\mathbb{N}_0$ one might start from a sequence being positive
definite with respect to the hypergroup. In this case, the measure
$\mu$ needs not to be stated explicitely. Its existence follows
from a Bochner's theorem (\ref{spec1}).
\end{Remark}

\subsection{The algorithms} The first ingredient of our algorithm
is a modified Chebyshev algorithm (\cite{Gau},p.76). It provides
the recursion coefficients $\alpha_k$, $\beta_k$ from the mixed
moments $\sigma_{k,l}$. As we do not require that the $P_n$ are
monic we get slightly different formulas compared to \cite{Gau}.

\begin{Proposition}
The recursion coefficients $\{\alpha_k,\; \beta_k\}_{k=1}^{n-1}$,
can be calculated from the

values $\{d(k)\}_{k=0}^{2n-1}$ and $\{a_k, b_k,
c_k\}_{k=0}^{2n-2}$ as follows:

Start with the initialization:
$$\sigma_{-1,-1}=1, \sigma_{-1,l}=\sigma_{l,-1}=0, \quad l=0,1, \ldots, 2n-2;$$
$$\sigma_{0,l}= d(l), \quad l=0,1, \ldots, 2n-1; \qquad \beta_0 = a_{-1}=c_0 = 0$$
and continue for $k=1,2, \ldots,
n-1; \quad l=k,k+1, \ldots, 2n-k-1$ with the recursion for the
mixed momentums
\begin{equation}
\sigma_{k,l}= a_l \sigma_{k-1,l+1} + [b_l  - \alpha_{k-1}] \sigma_{k-1, l}+
 c_l \sigma_{k-1,l-1}-\beta_{k-1}a_0 \sigma_{k-2,l} \label{eq:mixed moments}
 \end{equation}
and the coefficients
\begin{equation}
\beta_k= a_{k-1} \frac{\sigma_{k,k}}{\sigma_{k-1,k-1}}; \quad
\alpha_k = b_k + a_k \frac{\sigma_{k,k+1}}{\sigma_{k,k}} - a_{k-1}
\frac{\sigma_{k-1,k}}{\sigma_{k-1,k-1}}.\label{eq:beta}
\end{equation}
\end{Proposition}

\begin{proof}
For the mixed moments one finds
\begin{eqnarray*}
 \sigma_{k,l}  &=&  \int_{\mathbb{R}} \phi_k(x) P_l(x)d\mu(x)\\
  &=& \int_{\mathbb{R}} [x \phi_{k-1}(x)-\alpha_{k-1}\phi_{k-1}(x)-\beta_{k-1}\phi_{k-2}(x)] P_l(x)d\mu(x)\\
  &=& \int_{\mathbb{R}} \phi_{k-1}(x)[a_l P_{l+1}(x) + b_l P_l(x)+ c_l P_{l-1}(x)]d\mu(x)-
  \alpha_{k-1}\sigma_{k-1,l}-\beta_{k-1}\sigma_{k-2,l}\\
  &=& a_l \sigma_{k-1,l+1}+[b_l - \alpha_{k-1}]\sigma_{k-1,l}+
  c_l \sigma_{k-1,l-1}-\beta_{k-1}\sigma_{k-2,l}.
\end{eqnarray*}

The orthogonality of $\phi_k$ with respect to $\mu$ gives
$\sigma_{k,l}=0$ for all $k>l$. Especially one has for all $k \geq
1$
$$0=\sigma_{k+1,k-1}=a_{k-1}\sigma_{k,k}-\beta_k
\sigma_{k-1,k-1}$$ and thus
$$ \beta_k= a_{k-1} \frac{\sigma_{k,k}}{\sigma_{k-1,k-1}}.$$

Likewise one gets
$$0=\sigma_{k+1,k}=a_k \sigma_{k,k+1} + [b_k -\alpha_{k}] \sigma_{k,k}-\beta_k \sigma_{k-1,k}$$
giving
\begin{eqnarray*}
\alpha_k &=& b_k +  a_k \frac{\sigma_{k,k+1}}{\sigma_{k,k}} -
\beta_k \frac{\sigma_{k-1,k}}{\sigma_{k,k}} \\
&=& b_k +  a_k \frac{\sigma_{k,k+1}}{\sigma_{k,k}} - a_{k-1}
\frac{\sigma_{k-1,k}}{\sigma_{k-1,k-1}} .
\end{eqnarray*}

\end{proof}

With the received recursion coefficients $\alpha_k$ and $\beta_k$
we calculate as next step the connection coefficients $c_{\phi
P}(k,l)$ of the expansion. This can be done with an algorithm
given by Askey \cite{AsII} for connecting monic polynomials. We
again have to rewrite it, as the polynomials $P_s$ need not be
monic.

\begin{Proposition}\label{Askey}
The connection coefficients $\{c_{\phi P}(k,l)\}_{l=0}^{k}$ for
$k=1,\ldots,n$ can be

calculated from $\{\alpha_k,\beta_k\}_{k=0}^{n-1}$ and  $\{a_k,
b_k, c_k\}_{k=0}^{n}$ as follows:

Start with
$$c_{\phi P}(0,0)=\frac{1}{P_0}$$
and the conventions $\beta_0 = a_{-1}=c_0 = 0$, $ c_{\phi
P}(k,l)=0$ for $l>k$ or $k<0$ or $l<0$.

Continue for $1\leq k \leq n-1$ and $l=0,\ldots,k+1$ with

\begin{eqnarray*}
c_{\phi P}(k+1,l)&=&a_{l-1}c_{\phi P}(k,l-1)+ (b_l-\alpha_k)
c_{\phi P}(k,l)\\\\
&&+ c_{l+1} c_{\phi P}(k,l+1) -\beta_k c_{\phi P}(k-1,l).
\end{eqnarray*}

\end{Proposition}

\begin{proof}
$c_{\phi P}(0,0)=\frac{1}{P_0}$ follows from $\quad \phi_0(x)=1$.

The three term recursions of $\phi_k$ and $P_k$ show for $k=0,
\ldots,n-1$
\begin{eqnarray*}
\phi_{k+1}(x)&=& x \phi_{k}(x)-\alpha_k \phi_{k}(x)-\beta_k \phi_{k-1}(x)\\
&=& x \sum_{l=0}^k c_{\phi P}(k,l)P_l(x)-\alpha_k \sum_{l=0}^k
c_{\phi P}(k,l)P_l(x)-\beta_k \sum_{l=0}^{k-1}
c_{\phi P}(k-1,l)P_l(x)\\
&=& \sum_{l=0}^k c_{\phi
P}(k,l)[a_l P_{l+1}(x)+b_l P_l(x)+c_l P_{l-1}(x)]\\
& & -\alpha_k \sum_{l=0}^k c_{\phi P}(k,l)P_l(x)-\beta_k
\sum_{l=0}^{k-1} c_{\phi P}(k-1,l)P_l(x).
\end{eqnarray*}

Now, we reorder according to $P_l$. Hereby we shift indices and
use the above conventions.
\begin{eqnarray*}
\phi_{k+1}(x)&=& \sum_{l=0}^{k+1}\big[a_{l-1} c_{\phi P}(k,l-1)+(b_l -\alpha_k) c_{\phi P}(k,l)\\
&&+ c_{l+1}c_{\phi P}(k,l+1)- \beta_k c_{\phi P}(k-1,l)]P_l(x)
\end{eqnarray*}

These equations show an explicit representation of the connection
coefficients.
\end{proof}

The technique for a fast decomposition of our structured matrices
now follows from combining the given algorithms. We thereby get
rid of the coefficients $\alpha_k$ and $\beta_k$.

\begin{Theorem}
The coefficients $\{c_{\phi P}(k,l)\}_{0 \leq k,l \leq n}$ and the
mixed moment $\{\sigma_{k,k}\}_{k=0}^n$ required for the
decomposition
$$\mathbf{M_n}^{-1}= \mathbf{L}_n^{T}\mathbf{D}_n \mathbf{L}_n $$

can be found from the entries of $\mathbf{M_n}$, and the three
term recursion coefficients $\{a_k, b_k, c_k\}_{k=0}^{2n-2}$, with
the following algorithm having $O(n^2)$ complexity:

{\bf First Step: Calculation of $d(s)$}

For $k=0,\ldots,n$
$$d(k)=d(k,0)/P_0 .$$
For $k=1,\ldots,n$:
$$d(n+k)=\frac{1}{g(n,k,n+k)}\left(d(n,k)-\sum_{s=n-k}^{n+k-1}g(n,k,s)d(s)\right)$$

{\bf Second Step: Calculation of $\sigma_{k,l}$}

$\sigma_{-2,-2}=\sigma_{-1,-1}=1$, $a_{-1}=c_0 = 0$.

For $l=0,1, \ldots, 2n-2$
$$\sigma_{-1,l}=\sigma_{l,-1}=0.$$
For $0 \leq l<k \leq n$
$$\sigma_{k,l}=0.$$
For $\quad l=0,1, \ldots, 2n-1$
$$\sigma_{0,l}=d(l).$$
For $k= 1,\ldots,n; \quad l=k,\ldots 2n-k-1$
\begin{eqnarray*}
\sigma_{k,l}&= &a_l \sigma_{k-1,l+1} + [b_l - b_{k-1} -
a_{k-1}\frac{\sigma_{k-1,k}}{\sigma_{k-1,k-1}} +
a_{k-2}\frac{\sigma_{k-2,k-1}}{\sigma_{k-2,k-2}}]
\sigma_{k-1,l}\\
&&+ \; c_l
\sigma_{k-1,l-1}-a_{k-2}\frac{\sigma_{k-1,k-1}}{\sigma_{k-2,k-2}}\sigma_{k-2,l}.
\end{eqnarray*}

{\bf Third Step: Calculation of $c_{\phi P}(k,l)$}

$c_{\phi P}(0,0)=1/P_0$,  $c_{\phi P}(-1,-1)= c_{\phi
P}(-2,-2)=0$.

For $0 \leq l \leq n$

$c_{\phi P}(-1,l)= c_{\phi P}(l,-1)= c_{\phi P}(-2,l)= c_{\phi
P}(l,-2)=0$.

For $0 \leq k<l \leq n$
$$c_{\phi P}(k,l)=0.$$

For  $k=1,\ldots,n-1, \quad l=0,\ldots,k+1$
\begin{eqnarray*}
c_{\phi P}(k+1,l)&= &a_{l-1} c_{\phi P}(k,l-1) + [b_l -b_k - a_k
\frac{\sigma_{k,k+1}}{\sigma_{k,k}} +
a_{k-1}\frac{\sigma_{k-1,k}}{\sigma_{k-1,k-1}}]c_{\phi P}(k,l)\\
&&+ \; c_{l+1} c_{\phi P}(k,l+1) -
a_{k-1}\frac{\sigma_{k,k}}{\sigma_{k-1,k-1}} c_{\phi P}(k-1,l).
\end{eqnarray*}

\end{Theorem}

\begin{proof}
The calculation of $d(k)$ is given in Proposition \ref{moments}.
The representation of the mixed moments results from the modified
Chebyshev algorithm by inserting (\ref{eq:beta}) into
(\ref{eq:mixed moments}). Plugging the representations of
$\alpha_k$ and $\beta_k$ in the algorithm presented in Proposition
\ref{Askey} gives the stated recursions for $c_{\phi P}$. All
steps have complexity $O(n^2)$.
\end{proof}


\begin{thebibliography}{99}

\bibitem{andl}
 And\v{e}l, J.: Statistische Analyse von Zeitreihen.
 Akademie Verlag, Berlin, 1984.

\bibitem{Arn}
 Arnaud, J.P.: Stationary processes indexed by a homogeneous
tree. Ann. of Prob. 22 (1994), 195-218.

\bibitem{Aron}
Aronszajn, N.: Theory of reproducable kernels. Trans. Amer. Math.
Soc. 68 (1950), 337-404.

\bibitem{AsI}
Askey R.: Orthogonal polynomials, Mem. Amer. Math.
Soc. 18 (213), Providence, 1979.

\bibitem{Ask}
 Askey R.: Orthogonal Polynomials and Special Functions.
 SIAM, Philadelphia, 1975.

\bibitem{AsII}
 Askey R.: Orthogonal expansions with positive coefficients II,
 Siam J. Math. 2 (1971), 340-346.

\bibitem{Ressel}
Berg, C., Christensen, J. P. R., Ressel, P.: Harmonic Analysis on
Semigroups Theory of Positive Definite and Related Functions,
Springer, Berlin, 1984.

\bibitem{Bloo}
 Bloom, W.R., Heyer, H.: Harmonic Analysis of Probability Measures
on Hypergroups. De Gruyter, Berlin, 1995.

\bibitem{Blow}
Blower, G.: Stationary processes for translation operators. Proc.
London Math. Soc. 72 (1996), 697-720.

\bibitem{Bochner}
Bochner, S.: Vorlesung \"{u}ber Fouriersche Integrale, Leipzig,
1932.

\bibitem{Bo}
Bochner, S.:  Stationarity, boundedness, almost periodicity of
random valued functions. in: Proc. Third Berkeley Symp. Math.
Statist. Probab. Vol. 2, Berkeley, pp. 7-27, 1956.


\bibitem{Bres}
Bressoud, D.M.: Linearization and related formulas for
$q$-ultraspherical polynomials. SIAM J.Math. Anal. 12 (1981),
161-168.

\bibitem{Bro}
Brockwell, P.J., Davies, R.A.: Time Series: Theory and Methods.
Springer, New York, 1991.

\bibitem{ChaRao}
Chang, D.K., Rao, M.M.: Bimeasures and non stationary processes in
Real and Stochastic Analysis, Rao, M.M., Ed., Wiley, New York,
1986.

\bibitem{Cha}
Chang, D.K., Rao, M.M.: Special representations of weakly
harmonizable processes. Stochastic Anal. Appl. 6 (1988), 169-189.

\bibitem{Chi}
Chihara, T.S.:  An Introduction to Orthogonal Polynomials, Gordon
and Breach, New York, 1978.

\bibitem{Cow}
Cowling, M., Meda, St., Setti, A.G.: An overview of harmonic
analysis on the group of isometries of a homogeneous tree. Expo.
Math. 16 (1998), 385-424.

\bibitem{Cra}
Cram\'{e}r, H.: On the theory of stationary random processes. Ann.
Math. 41 (1940), 215-230.

\bibitem{CraHilb}
Cram\'{e}r, H.: Stochastic processes as curves in Hilbert spaces.
Theory. Probability Appl. 9 (1964), 193-204.

\bibitem{Deh}
Dehay , D.: Spectral analysis of the covariance of the almost
periodically correlated processes. Stochastic Processes Appl. 50
(1994), 315-330.


\bibitem{Del}
Delsarte, P., Genin, Y.: On the role of orthogonal polynomials on
the unit circle in digital signal processing applications
(115-133), in: Orthogonal polynomials, theory and practice, P.
Nevai (ed.), NATO Advanced Science Institutes Series C:
Mathematical and Physical Sciences, 294. Kluwer, 1990.


\bibitem{Dev}
 Devroye, L., Gy\"orfi, L.: Nonparametric Density
Estimation: The $L^1$-View. Wiley, New York, 1985.


\bibitem{Gard}
Gardner, W.A., Napolitano, A., Paura L.: Cyclostationarity: Half a
century of research. Signal Processing 86 (4) (2006), 639-697.

\bibitem{Gasper}
Gasper,G.: Linearizations of the product of Jacobi polynomials II.
Canad. J. Math. 22 (1970), 582-593.

\bibitem{Gau}
Gautschi W.: Orthogonal polynomials, computation and
approximation. Numerical Mathematics and Scientific Computation.
Oxford University Press, 2004.

\bibitem{Genossar}
Genossar M.J., Lev-Ari H., Kailath T.: Consistent estimation of
 the cyclic autocorrelation, IEEE Trans. Sig. Proc. 42 (1994),
 565-603.

\bibitem{Ger}
Geronimus, Y.L.: Polynomials Orthogonal on a Circle and Interval.
Pergamon, Oxford, 1960.


\bibitem{Geto}
Getoor, R.K.: The shift operator for non-stationary stochastic
processes. Duke Math. J. 23 (1956), 175-187.

\bibitem{Glad}
Gladyshev, E.G.: Periodically correlated random sequences. Soviet
Math. 2 (1961), 385-388.

\bibitem{Glady}
Gladyshev, E.G.: On multidimensional stationary random processes.
Theory Prob. and Appl. 3( 1958), 425-428.

\bibitem{Golup}
Golub, G.H., van Loan C.F.: Matrix computations, 2nd edn., The
John Hopkins University Press, 1989.


\bibitem{Graham}
Graham, C.C., McGehee, O.C.: Essays in Commutative Harmonic
Analysis, Springer, New York, 1979.


\bibitem{Han}
Hanin, L.G., Schreiber, B.S.:  Discrete spectrum of
nonstationary stochastic processes on groups. J. Theoretical
Probability 11 (1998), 1111-1133.


\bibitem{Heinig}
Heinig G., Rost K.: New fast algorithms for Toeplitz-plus-Hankel
matrices, SIAM Journal Matrix Anal. Appl. 25(3) (2003), 842-857.


\bibitem{Hey1}
 Heyer, H.:  Stationary random fields over hypergroups.
Gaussian Random Fields, World Scientific, pp. 197-213, 1991.

\bibitem{Hey2}
Heyer, H.: The covariance distribution of a generalized random
field over a commutative hypergroup. Contemporary Math. AMS 261,
Probability on Algebraic Structures, pp. 73-82, 2000.


\bibitem{wie}
 H\"{o}sel V., Lasser R.: A Wiener theorem for orthogonal polynomials,
 J. Funct. Anal. 133(2) (1995), 395-401.

\bibitem{hoe}
H\"{o}sel V., Lasser R.: One-step prediction for $P_n$-weakly
stationary processes, Monatshefte f\"{u}r Mathematik 113 (1992),
199-212.

\bibitem{hoecov}
 H\"osel, V.: On the estimation of covariance functions on
$P_n$-weakly stationary processes. Stoch. Anal. and Appl. 16
(1998), 607-629.

\bibitem{hoelaP}
H\"osel, V., Lasser, R.: Prediction of weakly stationary sequences
on polynomial hypergroups. Annals of Probability 31(1) (2003),
93-114.

\bibitem{Hof}
Hofmann, M.: Some Classes of Nonstationary Processes on Groups and
Hypergroups, Diploma thesis, TUM, 2006.

\bibitem{Hou}
Houdr\'{e}, C.: Harmonizability, V-boundedness, (2,p)-boundedness
of stochastic processes. Probab. Theory Rel. Fields 84 (1990),
39-54.

\bibitem{Houdre}
Houdr\'{e}, C.: Linear Fourier and stochastic analysis, Probab.
Theory Rel. Fields 87 (1990), 167-188.


\bibitem{HuRep}
Hurd, H.L.: Representation of strongly harmonizable periodically
correlated processes and their covariances. J. Multivar. Anal. 29
(1989), 53-67.

\bibitem{HuCor}
Hurd, H.L.: Correlation theory of almost periodically correlated
processes. J. Multivar. Anal. 37 (1991), 24-45.

\bibitem{Hurd}
Hurd, H.L., Miamee, A.: Periodically correlated random sequences.
Spectral theory and practice. Wiley Series in Probability and
statistics, 2007.


\bibitem{Ismail}
Ismail, M.E.H.: Classical and quantum orthogonal polynomials in
one variable. Cambridge University Press, Cambridge, 2005.

\bibitem{Isserlis}
Isserlis, L.: On a formula for the product-moment coefficient of
any order of a normal frequency distribution in any number of
variables. Biometrika 12 (1918), 134-139.

\bibitem{jew}
 Jewett, R.I.: Spaces with an abstract convolution of measures.
 Adv. in Math. 18 (1975), 1-101.

\bibitem{Kailath}
 Kailath, T., Sayed  A.H.: Displacement structure, theory and
applications. SIAM Review 37 (1995), 297-386.

\bibitem{Kak}
 Kakihara, Y.: Multidimensional Second Order Stochastic Processes.
Series on Multivariate Analysis 2, World Scientific Publ. Co.,
Inc., River Edge, NJ, 1997.

\bibitem{Kallia}
Kallianpur, G., Mandrekar, V.: Spectral theory of stationary
H-valued processes. J. Multivariate Anal. 1(1) (1971), 1-16.

\bibitem{Kam}
Kamp\'{e} de F\'{e}riet, J.: Correlation and spectrum of
asymptotically stationary random functions. Mathematics Student 30
(1962), 55-67.

\bibitem{KampFer}
Kamp\'{e} de F\'{e}riet, J., Frenkiel, F.N.: Correlation and
spectra for nonstationary random functions. Math. Comput. 16,
(1962), 1-21.

\bibitem{Khin}
Khintchine, A.Ya.: Korrelationstheorie der station\"{a}ren
stochastischen Prozesse. Math. Ann. 109 (1934), 605-615.


\bibitem{Kno}
Knopp, K.:  Theorie und Anwendung der unendlichen Reihen.
Springer, Berlin, 1922.

\bibitem{Kol}
Kolmogorov, A.N.: Stationary sequences in Hilbert space. Bull.
Moskov. Gos. Univ. Matematika 2 (1941), 1-40.


\bibitem{Lasser}
 Lasser, R.: A modification of stationarity for stochastic processes
 induced by orthogonal polynomials. H. Heyer, ed., Probability
 Measures  on Groups  IX , Lecture Notes Math. 1379, Springer, Berlin, pp. 185-191, 1989.

\bibitem{LaMM}
 Lasser, R.: On the modified momentum problem with solution
 carried by $[-1,1]$. H. Heyer, ed.,Probability Measures on Groups and related
 Structures, World Scientific Publications, 1995.


\bibitem{la01}
 Lasser, R.: Orthogonal polynomials and hypergroups.
 Rend. di. Mat. 3 (1983), 185-209.

\bibitem{la02}
 Lasser, R.: Orthogonal polynomials and hypergroups II - the symmetric
 case. Trans. Amer. Math. Soc. 341 (1994), 749-770.


\bibitem{laleo}
 Lasser, R., Leitner, M.: On the estimation of the mean of weakly
 stationary and polynomial weakly stationary sequences
 J. Mult. Anal. 35 (1990), 31-47.

\bibitem{lalest}
 Lasser, R., Leitner, M.: Stochastic processes indexed by hypergroups I,
 J. Theor. Prob. 2 (1989), 301-311.

\bibitem{laob}
 Lasser, R., Obermaier, J., Strasser, W.: On the
consistency of weighted orthogonal series density estimators with
respect to $L^1$-norm. Nonparametric Statistics 3 (1993), 71-80.

 \bibitem{laro}
 Lasser, R., R\"{o}sler, M.: Linear mean estimation of weakly stationary
 stochastic processes under the aspect of optimality and asymptotic
 optimality. Stoch. Proc. Appl. 38 (1991), 279-293.

\bibitem{lei2}
 Leitner, M.: Stochastic processes indexed by hypergroups II.
 J. Theor. Prob. 4 (1991) 321-332.

\bibitem{Leihyp}
 Leitner, M.: Hyper-weakly harmonizable processes and operator families.
 Stochastic Anal.Appl. 13 (1995), 471-485.

\bibitem{Lindahl}
Lindahl; R.J., Maserick; P.H.: Positive-definite functions on
involutionsemigroups. Duke. Math. J. 38 (1971), 771-782.

\bibitem{loeve}
Lo\`{e}ve, M.: Probability Theory. Van Nostrand,Princeton, N.J.,
1955.

\bibitem{Lub}
Lubinsky, D.S.: A survey of general orthogonal polynomials for
weights on finite and infinite intervals. Acta Appl. Math. 10
(1987), 237-296.

\bibitem{Man}
Mandrekar, V., Salehi, H.: Prediction Theory and Harmonic
Analysis. North-Holland, Amsterdam, 1983.

\bibitem{Marty}
Marty, F., R\^{o}le de la notation de hypergroupe dans l'\'{e}tude
de groupes non ab\'{e}liens. C.R. Acad. Sci. Paris 201 (1935),
636-638.

\bibitem{Mate}
M\'{a}t\'{e}, A., Nevai,P.: Bernstein's inequality in $L^p$ for $0
< p <1$ and $(C,1)$ bounds for orthogonal polynomials. Ann. of
Math. 111 (1980), 145-154.

\bibitem{Mercer}
Mercer, J.: Functions of positive and negative type, and their
connection with the theory of integral equations. Proc. Royal Soc.
Series A, 83 (559) (1909), 69-70.


\bibitem{Mia}
Miamee, A. G.: Periodically correlated processes and their
stationary dilations. SIAM J. Appl. Math. 50 (1990), 1194-1199.


\bibitem{MiaSal}
Miamee, A. G.,  Salehi H.: On the prediction of periodically
correlated stochastic processes. in: Multivariate Analysis V,
North-Holland, Amsterdam, pp.167-179, 1980.

\bibitem{Mitr}
Mitrinovi\'{c}, D.S.: Analytic Inequalities. Springer, Berlin,
1970.

\bibitem{nev}
Nevai, P.: Orthogonal Polynomials. Memoirs Amer. Math. Soc. 213,
Providence, R.I., 1979.

\bibitem{Nie}
Niemi, H.: On stationary dilations and the linear
prediction of certain stochastic processes.
 Soc. Sci. Fenn. Comment. Phy.-Math. 45 (1975), 111-130.

\bibitem{Niemi}
Niemi, H.: Stochastic processes as Fourier transforms of
stochastic measures. Annal. Acad. Sci. Fenn. AI Math. 591 (1975),
1-47.

\bibitem{Olshevsky}
Olshevsky, V., (ed.): Fast algorithms for structured matrices,
theory and applications. AMS-IMS-SIAM Joint summer research
conference on fast algorithms in mathematics, computer science and
engineering, 2001. Amer. Math. Soc., Contemp. math. 323, 2003.

\bibitem{Pan}
Pan, V. Y.: Structured Matrices and Polynomials. Unified superfast
algorithms. Birkh\"{a}user, Boston, 2001.

\bibitem{Part}
Parthasarathy, K.R.: Positive-definite kernels, continuous tensor
products and central limit theorems of probability theory. Lecture
Notes in Mathematics 272, 1972.

\bibitem{Parzen}
Parzen, E.: Spectral analysis of asymptotically stationary time
 series. Bull. Internat. Statist. Inst. 39 (1962) 87-103.

\bibitem{Pri}
Priestley, M.B.: Non-linear and Non-stationary Time Series
Analysis. Academic Press, New York, 1988.


\bibitem{Rao}
 Rao M.M.: Real and Stochastic Analysis. Wiley, New York, 1986.

\bibitem{RaoBiH}
Rao, M.M.: Bimeasures and harmonizable processes. In: Lecture
Notes in Math. Vol. 1379, Probability Measures on Groups IX, ed.
H.Heyer, Springer, pp.254-298, 1989.

\bibitem{RaoHa}
Rao, M.M.: Harmonizable processes: structure theory. L'Enseign.
Math. 28 (1982), 295-351.

\bibitem{RaoView}
Rao, M.M.: A view of harmonizable processes. in: Statistical Data
Analysis and Inference, North-Holland, Amsterdam, pp. 597-615,
1989.

\bibitem{RiSz}
Riesz, F., Sz.-Nagy, B.: Vorlesungen \"{u}ber Funktionalanalysis.
VEB Deutscher Verlag der Wissenschaften, Berlin, 1956.

\bibitem{Roes}
R\"{o}sler, M.: On the dual of a commutative signed hypergroup.
manuscripta math. 88 (1995), 147-163.

\bibitem{Roza}
Rozanov, Y.A.: Spectral theory of multidimensional stochastic
random processes with discrete time. Uspekhi Mat. Nauk. XIII(2)
(1958), 93-142.

\bibitem{Ros}
Rosenblatt, M.: Stationary Sequences and Random Fields.
Birkh\"{a}user, Boston, 1985.

\bibitem{Rozanov}
Rozanov, Y.A.: Spectral analysis of abstract functions. Theor.
Prob. Appl. 4 (1959), 271-287.

\bibitem{rudre}
Rudin, W.: Real and Complex Analysis. Third edition. Mc Graw Hill,
New York, 1983.

\bibitem{Sas}
Sasv\'{a}ri, Z.: Positive Definite and Definitizable Functions.
Akademie Verlag, Berlin, 1994.

\bibitem{SchAsRel}
Schreiber, B.S.:Asymptotically stationary and related processes.
Stochastic Processes and Functional Analysis. A Volume of Recent
Advances in Honor of M.M. Rao, A.C. Krinik and R.J. Swift, eds.,
Lect. Notes in Pure and Appl. Math. Vol. 238, pp. 363-397, Marcel
Dekker, New York, 2004


\bibitem{SchAs}
Schreiber, B.S.: Asymptotically stationary processes on amenable
groups. Stochastic Analysis and Appl. 22 (2004), 1525-1551.

\bibitem{Schur}
Schur, I.: \"{U}ber Potenzreihen, die im Innern des
Einheitskreises beschr\"{a}nkt sind. J. Reine Angew. Math. 147
(1917), 205-232.

\bibitem{schu}
 Schuster, A.: On the investigation of hidden periodicities with
 application to a supposed 26-day period of meteorological
 phenomena.Ter. Mag. Atmos. Elect. 3 (1898), 13-41.


\bibitem{Shir}
Shiryayev, A.N.:Probability. Springer, New York, 1984.

\bibitem{Stewart}
Stewart J.: Positive definite functions and generalizations. An
historical survey. Rocky Mountain J.Math 6(3) (1976), 409-429.

\bibitem{Stone}
Stone, M.H.: On one parameter unitary groups in Hilbert spaces.
Annals of Math. 33 (1932), 643-648.

\bibitem{Szeg}
Szeg\"o, G.: Orthogonal Polynomials. Amer. Math. Soc., Providence,
1975.

\bibitem{Sz}
Szeg\"o, G.: Orthogonal Polynomials. Amer. Math. Soc., New York,
1959.

\bibitem{Nagy}
Sz.-Nagy, B.: Transformation de l'espace de Hilbert, fonctions de
type positif sur un groupe. Acta Sci. Math 15 Szeged (1954),
104-114.

\bibitem{Czech}
Sz.-Nagy, B.: Contibution en Hongrie \`{a} la th\'{e}orie
spectrale des transformations lin\'{e}aires. Czech. Math. J. 6(2)
(1956). 166-176.

\bibitem{Voit1}
Voit, M.: Central limit theorems for random walks on ${\Bbb N}_0$
that are associated with orthogonal polynomials. J.Multivar. Anal.
34 (1990), 290-322.

\bibitem{Voit2}
Voit, M.: Factorization of probability measures on symmetric
hypergroups. J.Austral. Math. Soc. Series A 50 (1991), 417-467.

\bibitem{Wall}
Wall, H.S.: Hypergroups. Amer.J.Math. 59 (1937), 77-98.

\bibitem{Wien}
Wiener, N.: Generalized harmonic analysis. Acta Mathematika, V 55,
pp. 117-258.1930.

\bibitem{WiMa1}
Wiener, N., Masani, P.: The prediction theory
of multivariate stochastic processes, Part 1. Acta Math. 98
(1957), 111-150.

\bibitem{WiMa2}
Wiener, N., Masani, P.: The prediction theory of
multivariate stochastic processes, Part II. Acta Math. 99 (1958),
93-137.


\bibitem{Yag}
Yaglom, A.M.: Correlation Theory of Stationary and Related Random
Functions. Vol.s I+II, Springer, New York, 1987.



\end{thebibliography}
\end{document}